\input amstex 
\documentstyle{amsppt} 
\loadbold
\let\bk\boldkey

\magnification=1200
\hsize=5.75truein
\vsize=8.75truein 
\hcorrection{.25truein}
\loadeusm 
\let\scr\eusm
\loadeurm 

\font\Rm=cmr12 
\font\rrm=cmr12 scaled \magstep1

\font\fScr=eusb10 \define\Scr#1{\text{\fScr #1}} 
 
\define\Ind{\text{\rm Ind}}
\define\Aut#1#2{\text{\rm Aut}_{#1}(#2)}
\define\End#1#2{\text{\rm End}_{#1}(#2)}
\define\Hom#1#2#3{\text{\rm Hom}_{#1}({#2},{#3})} 
\define\GL#1#2{\roman{GL}_{#1}(#2)}
\define\M#1#2{\roman M_{#1}(#2)}
\define\Gal#1#2{\text{\rm Gal\hskip.5pt}(#1/#2)}
\define\Ao#1#2{\scr A^0_{#1}(#2)} 
\define\Go#1#2{\scr G_{#1}(#2)} 
\define\upr#1#2{{}^{#1\!}{#2}} 
\define\uL#1{{}^{\,\text{\rm L}}{#1}} 
\define\pre#1#2{{}_{#1\!}{#2}} 
\define\N#1#2{\text{\rm N}_{#1/#2}} 
\define\sw{\text{\rm sw}} 
 
\define\Mid{\,\big|\,} 
\define\ewc#1{\Scr E^{\text{\rm C}}(#1)} 
\define\biP#1#2#3{\upr{2\,}\Psi_{(#1/#2,#3)}} 
\define\bip#1{\upr{2\,}\Psi_#1} 
\let\ge\geqslant
\let\le\leqslant
\let\ups\upsilon
\let\vD\varDelta 
\let\ve\varepsilon 
\let\eps\epsilon 
\let\vf\varphi
\let\vF\varPhi 
\let\vG\varGamma

\let\vS\varSigma 
\let\vs\varsigma 
\let\vt\vartheta
\let\vT\varTheta
\let\vU\varUpsilon 
\let\vX\varXi 
\let\wt\widetilde 
\define\wP#1{\widehat{\scr P}_{#1}} 
\define\wW#1{\widehat{\scr W}_{#1}} 
\define\wss#1{\widehat{\scr W}^{\text{\rm ss}}_{#1}} 
\define\wwr#1{\widehat{\scr W}^{\text{\rm wr}}_{#1}} 
\define\awr#1{\widehat{\scr W}^{\text{\rm awr}}_{#1}} 
\define\wG#1{\widehat{\text{\rm GL}}_{#1}} 
\define\wA#1{{\widehat{\scr A\,}\!}_{#1}} 
\define\wR#1#2{\widehat{\scr R}_{#1}({#2})} 
\define\wRp#1#2{\widehat{\scr R}^+_{#1}({#2})} 
\topmatter \nologo \nopagenumbers
\title\nofrills \rrm 
Local Langlands correspondence and ramification for Carayol representations
\endtitle 
\rightheadtext{Carayol representations} 
\author 
Colin J. Bushnell and Guy Henniart 
\endauthor 
\leftheadtext{C.J. Bushnell and G. Henniart}
\affil 
King's College London and Universit\'e de Paris-Sud 
\endaffil 
\address 
King's College London, Department of Mathematics, Strand, London WC2R 2LS, UK. 
\endaddress
\email 
colin.bushnell\@kcl.ac.uk 
\endemail
\address 
Laboratoire de Math\'ematiques d'Orsay, Univ Paris-Sud, CNRS, Universit\'e
Paris-Saclay, 91405 Orsay, France.
\endaddress 
\email 
Guy.Henniart\@math.u-psud.fr 
\endemail 
\date May 2017, revised September 2018 \enddate 
\abstract 
Let $F$ be a non-Archimedean locally compact field of residual characteristic $p$ with Weil group $\Cal W_F$. Let $\sigma$ be an irreducible smooth complex representation of $\Cal W_F$, realized as the Langlands parameter of an irreducible cuspidal representation $\pi$ of a general linear group over $F$. In an earlier paper, we showed that the ramification structure of $\sigma$ is determined by the fine structure of the endo-class $\varTheta$ of the simple character contained in $\pi$, in the sense of Bushnell-Kutzko. The connection is made via the {\it Herbrand function} $\Psi_\varTheta$ of $\varTheta$. In this paper, we concentrate on the fundamental Carayol case in which $\sigma$ is totally wildly ramified with Swan exponent not divisible by $p$. We show that, for such $\sigma$, the associated Herbrand function satisfies a certain functional equation, and that this property essentially characterizes this class of representations. We calculate $\Psi_\varTheta$ explicitly, in terms of a classical Herbrand function arising naturally from the classification of simple characters. We describe exactly the class of functions arising as Herbrand functions $\Psi_\varXi$, as $\varXi$ varies over the set of totally wild endo-classes of Carayol type. In a separate argument, we derive a complete description of the restriction of $\sigma$ to any ramification subgroup and hence a detailed  interpretation of the Herbrand function. This gives concrete information concerning the Langlands correspondence. 
\endabstract 
\toc{} 
\subhead 1. Classical Herbrand functions \endsubhead 
\subhead 2. Certain simple characters \endsubhead 
\subhead 3. Functional equation \endsubhead 
\subhead 4. Symmetry and the bi-Herbrand function \endsubhead 
\subhead 5. Characters of restricted level \endsubhead 
\subhead 6. Variation of parameters \endsubhead 
\subhead 7. The Herbrand function \endsubhead 
\subhead 8. Representations with a single jump \endsubhead 
\subhead 9. Ramification structure \endsubhead 
\subhead 10. Parameter fields \endsubhead 
\endtoc 
\endtopmatter 
\document \baselineskip=14pt \parskip=4pt plus 1pt minus 1pt
\subhead 
1 
\endsubhead 
Let $F$ be a non-Archimedean, locally compact field with residual characteristic $p$. Let $\scr W_F$ be the Weil group of a separable closure $\bar F/F$. For a real variable $x\ge0$, let $\scr R_F(x) = \scr W_F^x$ be the corresponding ramification subgroup of $\scr W_F$ and $\scr R^+_F(x)$ the closure of $\bigcup_{y>x} \scr R_F(y)$. We use the conventions of \cite{38} here, so that $\scr R_F(0)$ is the inertia group $\scr I_F$ and $\scr R_F^+(0)$ is the wild inertia group $\scr P_F$ in $\scr W_F$. If $\scr G$ is any of this list of locally profinite groups, $\widehat{\scr G}$ will denote the set of equivalence classes of irreducible, smooth, complex representations of $\scr G$. We shall be concerned with the {\it ramification structure\/} of certain $\sigma\in \wW F$, that is, the structure of the restricted representations $\sigma\Mid \scr R_F(x)$ and $\sigma\Mid \scr R_F^+(x)$, for $x>0$. 
\par 
On the other side, let $\Ao nF$ denote the set of equivalence classes of irreducible, cuspidal, complex representations of the general linear group $\GL nF$, $n\ge1$, and set $\wG F = \bigcup_{n\ge1} \Ao nF$. For $\pi \in \wG F$, write $\roman{gr}(\pi) = n$ to indicate $\pi \in \Ao nF$. Such a representation $\pi$ contains a {\it simple character\/} $\theta_\pi$ in $\GL nF$ \cite{15} and, up to conjugation, only one \cite{10}. The {\it endo-class\/} $\vT_\pi$ of $\theta_\pi$ is therefore uniquely determined by $\pi$. Let $\Scr E(F)$ denote the set of endo-classes of simple characters over $F$. (For the notion of endo-class, see \cite{3} or the summary in any of \cite{2}, \cite{5}, \cite{10}.) 
\par 
Denote by $\pi \mapsto \uL\pi$ the Langlands correspondence $\wG F \to \wW F$ \cite{20}, \cite{25}, \cite{35}, \cite{37}. Writing $\sigma = \uL\pi$, the fine structure of the endo-class $\vT_\pi$ and the ramification structure of $\sigma$ determine each other \cite{13} 6.5 Corollary. The relationship is expressed via a certain {\it Herbrand function\/} $\Psi_{\vT_\pi}$ attached to the endo-class $\vT_\pi$. In this paper, we consider a particularly interesting class of representations, comprising what we call {\it Carayol representations.} We compute the associated Herbrand functions. We list the functions which arise as Herbrand functions. We interpret the results in terms of the ramification structure of the associated Galois representations, from which we extract information about the Langlands correspondence. 
\subhead 
2 
\endsubhead 
We review the background from \cite{13} with as little formality as possible. If $\pi \in \wG F$ and $\sigma = \uL\pi\in \wW F$, the endo-class $\vT_\pi$ determines the restriction $\sigma\Mid \scr P_F$. More precisely, $\sigma$ defines an element $[\sigma]_0^+$ of the orbit space $\scr W_F\backslash \wP F$, namely the orbit of irreducible components of $\sigma\Mid \scr P_F$. The Langlands correspondence induces a canonical bijection (\cite{5} 8.2 Theorem, \cite{12} 6.1) 
$$ 
\aligned 
\Scr E(F) &\longrightarrow \scr W_F\backslash \wP F, \\ \vT &\longmapsto \uL\vT, 
\endaligned 
\tag A 
$$ 
by 
$$ 
[\uL\pi]_0^+ = \uL\vT_\pi, \quad \pi \in \wG F. 
$$ 
Results developed through \cite{3}, \cite{4}, \cite{5}, \cite{6}, \cite{7}, \cite{9} and particularly \cite{12} show that the map (A) is central to understanding of the Langlands correspondence. 
\subhead 
3 
\endsubhead  
The starting point of \cite{13} is that each of the sets $\Scr E(F)$, $\scr W_F\backslash \wP F$ carries a canonical {\it ultrametric.} That on $\Scr E(F)$, denoted by $\Bbb A$, is built on the fact that simple characters are characters of compact groups carrying canonical filtrations, and those filtrations provide a medium via which the characters may be compared. The ultrametric $\Bbb A$ relates to Swan exponents of pairs of representations, as defined from the local constants of \cite{29}, \cite{39}. Let $\vT \in \Scr E(F)$ and choose $\pi\in \wG F$ such that $\vT_\pi = \vT$. There is a unique continuous function $\vF_\vT(x)$, $x\ge0$, such that 
$$ 
\vF_\vT\big(\Bbb A(\vT,\vT_\rho)\big) = \frac{\sw(\check\pi\times \rho)}{\roman{gr}(\pi)\,\roman{gr}(\rho)}\,, 
$$ 
for any $\rho\in \wG F$. The function $\vF_\vT$ is piecewise linear, strictly increasing and convex. It is given by an explicit formula \cite{13} (4.4.1) derived from the conductor formula of \cite{14} 6.5 Theorem. We call $\vF_\vT$ the {\it structure function\/} of $\vT$. 
\par 
The ultrametric on $\scr W_F\backslash \wP F$, denoted by $\Delta$, is defined by comparing representations via the canonical filtration of $\scr P_F$ by ramification groups: for $\sigma,\tau\in \wW F$, 
$$ 
\Delta([\sigma]_0^+,[\tau]_0^+) = \inf \{x>0: \Hom{\scr R_F(x)}\sigma\tau \neq 0\}. 
$$ 
The ultrametric $\Delta$ likewise relates to Swan exponents of tensor products of pairs of representations of $\scr W_F$ \cite{21}. For $\sigma\in \wW F$, there is a unique continuous function $\vS_\sigma(x)$, $x\ge 0$, such that 
$$ 
\vS_\sigma\big(\Delta([\sigma]_0^+,[\tau]_0^+)\big) = \frac{\sw(\check\sigma\otimes \tau)}{\dim\sigma \cdot \dim\tau}\,, 
$$ 
for all $\tau\in \wW F$. The function $\vS_\sigma$ is piecewise linear, strictly increasing and convex. It is given by a formula derived from the ramification structure of $\sigma$, reproduced in (2.2.2) below. If $\vS_\sigma$ is smooth at $x$, its derivative satisfies 
$$ 
\vS_\sigma'(x) = \dim\End{\scr R_F(x)}\sigma\big/(\dim\sigma)^2\big.. 
$$ 
We call $\vS_\sigma$ the {\it decomposition function\/} of $\sigma$: it depends only on the orbit $[\sigma]_0^+$.   
\par 
If $\vT \in \Scr E(F)$, set $\Psi_\vT = \vF_\vT^{-1}\circ\vS_\sigma$, for any $\sigma\in \wW F$ such that $[\sigma]_0^+ = \uL\vT$. The Langlands correspondence respects Swan exponents of pairs and $\dim(\uL\pi) = \roman{gr}(\pi)$, $\pi\in \wG F$, so  
$$ 
\Psi_\vT\big(\Delta(\uL\vT, \uL\vX)\big) = \Bbb A(\vT,\vX), \quad \vX \in \Scr E(F).  
$$ 
The function $\Psi_\vT$ is called the {\it Herbrand function\/} of $\vT$. It is continuous, strictly increasing and piecewise linear. 
\par 
If we take the view that $\vT \in \Scr E(F)$ has been given, in terms of the standard classification from \cite{15}, it is a simple matter to write down the function $\vF_\vT$. The Interpolation Theorem \cite{13} 7.5 shows, in principle, how to compute $\Psi_\vT$ {\it directly from $\vT$,} without reference to $\uL\vT$. It yields the decomposition function $\vS_\sigma$ and therefore a numerical account of the ramification structure of $\sigma$, just in terms of $\vT$. The Interpolation Theorem is not easy to apply directly, but it is the foundation of much of what we do here. 
\subhead 
4 
\endsubhead 
We specify the classes of representation on which we focus. 
\par 
Let $\vT \in \Scr E(F)$. Assuming, as we invariably do, that $\vT$ is non-trivial, it is the endo-class of a simple character $\theta\in \scr C(\frak a,\beta)$ attached to a simple stratum $[\frak a,m,0,\beta]$ in some matrix algebra $\M nF$ (following the conventions of \cite{15}). In particular, $\beta \in \GL nF$ and the algebra $F[\beta]$ is a field: one says that $F[\beta]$ is a {\it parameter field\/} for $\vT$. The positive integers $\deg\vT = [F[\beta]{:}F]$ and $e_\vT = e(F[\beta]|F)$ are invariants of $\vT$. The {\it slope\/} $\vs_\vT$ of $\vT$, defined by $\vs_\vT = m/e_\frak a$, where $e_\frak a$ is the period of the hereditary $\frak o_F$-order $\frak a$, is likewise an invariant of $\vT$. If $\pi \in \wG F$ satisfies $\vT_\pi = \vT$, then $\vs_\vT = \sw(\pi)/\roman{gr}(\pi)$. However, neither $\theta$ nor $\vT$ determines the parameter field $F[\beta]$: see the later parts of section 6. 
\par 
Say that $\vT \in \Scr E(F)$ is {\it totally wild\/} if $\deg\vT = e_\vT = p^r$, for an integer $r\ge0$. If $\vT$ is totally wild, say it is {\it of Carayol type\/} if $\deg\vT > 1$ and the integer $e_\vT\vs_\vT$ is not divisible by $p$. Let $\ewc F$ denote the set of endo-classes $\vT \in \Scr E(F)$ that are totally wild of Carayol type. We aim to calculate $\Psi_\vT$ for all $\vT \in \ewc F$. 
\par
We concentrate on this case for two reasons. First, 7.1 Proposition of \cite{13} reduces the problem of calculating Herbrand functions to the totally wild case. Second, we have to work with simple characters. The definition of simple character in \cite{15} is rigidly hierarchical in nature and proofs are almost always inductive along this hierarchy. The first inductive step concerns the case where the element $\beta$ (as above) is {\it minimal over\/} $F$ \cite{15} (1.4.14). For totally wild endo-classes, this is the Carayol case.  
\par 
On the other side, say that $\sigma \in \wW F$ is {\it totally wild\/} if the restriction $\sigma\Mid \scr P_F$ is irreducible. In particular, $\dim\sigma = p^r$, for some $r\ge0$. Denote by $\wwr F$ the set of totally wild elements of $\wW F$. An endo-class $\vT \in \Scr E(F)$ is then totally wild if and only if there exists $\sigma\in \wwr F$ such that $[\sigma]_0^+ = \uL\vT$ ({\it cf\.} section 6 of \cite{12}). Say that $\sigma \in \wwr F$ is {\it of Carayol type\/} if $\dim\sigma\neq 1$ and $p$ does not divide $\sw(\sigma)$. Thus $\sigma\in \wwr  F$ is of Carayol type if and only if $[\sigma]_0^+ = \uL\vT$, for some $\vT \in \ewc F$. We shall see that these representations $\sigma$ exhibit a family of quite singular properties, reflecting the special nature of the endo-classes $\vT\in \ewc F$.
\subhead 
5 
\endsubhead 
We review our main results. They are organized into three principal theorems, that compl\-ement and support each other, followed by a substantial application. 
\par 
For any $\vT \in \Scr E(F)$, the Herbrand function $\Psi_\vT(x)$ satisfies $\Psi_\vT(0) = 0$ and $\Psi_\vT(x) = x$ for $x\ge \vs_\vT$ \cite{13} 6.2 Proposition. The derivative $\Psi'_\vT(x)$ has only finitely many discontinuities in the interesting region $0<x<\vs_\vT$: we call them the {\it jumps of\/} $\Psi_\vT$. When $\vT \in \ewc F$, the function $\Psi_\vT(x)$ is {\it convex\/} in the region $0\le x\le \vs_\vT$. The reasons for this are simple (2.4), but the property is very useful. 
\proclaim{Theorem 1} 
Let $\vT \in \ewc F$. The graph $y= \Psi_\vT(x)$, $0\le x\le \vs_\vT$, is symmetric with respect to the line $x{+}y = \vs_\vT$. That is, 
$$ 
\vs_\vT-x = \Psi_\vT\big(\vs_\vT - \Psi_\vT(x)\big), \quad 0\le x\le \vs_\vT. 
\tag B
$$ 
\endproclaim 
Theorem 1 has a satisfying converse. The group of characters of $U^1_F$ acts on the set $\Scr E(F)$ following the natural twisting action of characters of $F^\times$ or $\scr W_F$ on $\wG F$ or $\wW F$. We denote this action by $(\chi,\vT) \mapsto \chi\vT$. It has the property $\Psi_{\chi\vT} = \Psi_\vT$ \cite{13} 7.4 Proposition. We get: 
\proclaim{Corollary} 
Let $\vT \in \Scr E(F)$ be totally wild of degree $p^r$, for some $r\ge 1$, and suppose that $\vs_\vT \le \vs_{\chi\vT}$ for all characters $\chi$ of $U^1_F$. The function $\Psi_\vT$ then has the symmetry property \rom{(B)} if and only if $\vT \in \ewc F$.  
\endproclaim 
Theorem 1, together with some preliminary calculations, suggests the definition of a family of elementary functions. Let $r\ge 1$ and let $E/F$ be a totally ramified field extension of degree $p^r$. Let $m$ be a positive integer not divisible by $p$ and set $\vs = m/p^r$. Let $\psi_{E/F}$ be the classical Herbrand function of $E/F$ \cite{18}, \cite{38}. Define $c$ by the equation $c + p^{-r}\psi_{E/F}(c) = \vs$. There is then a unique function $\biP EF\vs(x)$, defined for $0\le x\le \vs$, such that the graph $y = \biP EF\vs(x)$ is symmetric with respect to the line $x{+}y=\vs$ and $\biP EF\vs(x) = p^{-r}\psi_{E/F}(x)$, for $0\le x\le c$. Functions of this form will be called {\it bi-Herbrand functions.} 
\par 
Our strategy is to identify $\Psi_\vT$, $\vT \in \ewc F$, as a specific bi-Herbrand function. Let $\deg\vT = p^r$. There is a simple stratum $[\frak a,m,0,\alpha]$ in $\M{p^r}F$ such that $\vT$ is the endo-class of some $\theta\in \scr C(\frak a,\alpha)$. Thus $F[\alpha]/F$ is totally ramified of degree $p^r$ and $p$ does not divide $m = -\ups_{F[\alpha]}(\alpha)$. In this notation, $\vs_\vT = m/p^r$.  If $\|\scr C(\frak a,\alpha)\|$ denotes the set of endo-classes of elements of $\scr C(\frak a,\alpha)$, then $\|\scr C(\frak a,\alpha)\| \i \ewc F$. 
\par 
The set  $\|\scr C(\frak a,\alpha)\|$ is not well-adapted to our purposes, because the function $\vT \mapsto \Psi_\vT$ is not constant there. Indeed, it may vary widely: see 7.2 Theorem 1. To overcome this problem, we specify a non-empty subset $\scr C^\star(\frak a,\alpha)$ of $\scr C(\frak a,\alpha)$, using an explicit formula given in 7.1 below: we say that $\theta$ {\it conforms to\/} $\alpha$ to indicate $\theta\in \scr C^\star(\frak a,\alpha)$. Let $\|\scr C^\star(\frak a,\alpha)\|$ denote the set of endo-classes of characters $\theta \in \scr C^\star(\frak a,\alpha)$. 
\proclaim{Theorem 2} 
Let $\vT \in \ewc F$ have degree $p^r$ and $\vs_\vT = m/p^r$. There is a simple stratum $[\frak a,m,0,\alpha]$ in $\M{p^r}F$ such that $\vT\in \|\scr C^\star(\frak a,\alpha)\|$. For any such stratum, 
$$ 
\Psi_\vT(x) = \biP{F[\alpha]}F{\vs_\vT}(x), \quad 0\le x\le \vs_\vT. 
\tag C
$$ 
\endproclaim 
Theorem 2 has the following consequence. 
\proclaim{Corollary} 
Let $E/F$ be a totally ramified field extension of degree $p^r$, $r\ge 1$, and let $m$ be a positive integer not divisible by $p$. There exists $\vT \in \ewc F$, with parameter field $E/F$, such that 
$$ 
\Psi_\vT(x) = \biP EF{m/p^r}(x), \quad 0\le x \le m/p^r = \vs_\vT. 
$$ 
\endproclaim 
The corollary is an effective tool for constructing representations of $\scr W_F$ with specified ramification properties. An application of the technique is given in 9.7. 
\subhead 
6 
\endsubhead 
Our third result looks at the problem from the Galois side. Let $\sigma \in \wwr F$ be of Carayol type and dimension $p^r$. Define $\vT \in \ewc F$ by $[\sigma]_0^+ = \uL\vT$. As $r\ge1$, the function $\Psi_\vT$ has at least one jump \cite{13} 7.7. If $\Psi_\vT$ has {\it exactly\/} one jump, we say that $\sigma$ is {\it H-singular.} In section 8, we analyze the structure of such representations in some detail: they belong to a rather special class of ``Heisenberg representations'' (as one says). 
\par 
Without restriction on the number of jumps, define a number $c_\vT$ by the equation 
$$
c_\vT + \Psi_\vT(c_\vT) = \vs_\vT, \quad \vT \in \ewc F. 
$$ 
By the symmetry of Theorem 1, $c_\vT$ is a jump of $\Psi_\vT$ if and only if $\Psi_\vT$ has an odd number of jumps and, in that case, $c_\vT$ is the middle one. 
\proclaim{Theorem 3} 
Let $\sigma \in \wwr F$ be of Carayol type and dimension $p^r$. Let $\vT \in \ewc F$ satisfy $\uL\vT = [\sigma]_0^+$. 
\roster 
\item The restriction $\sigma\Mid \scr R_F^+(c_\vT)$ is a direct sum of characters. 
\item 
Let $\xi$ be a character of $\scr R_F^+(c_\vT)$ occurring in $\sigma$, let $\scr W_{L_\xi}$ be the $\scr W_F$-stabilizer of $\xi$, and let $\sigma_\xi$ be the natural representation of $\scr W_{L_\xi}$ on the $\xi$-isotypic subspace of $\sigma\Mid \scr R_F^+(c_\vT)$. The field extension $L_\xi/F$ is totally ramified of degree dividing $p^r$ and $\sigma = \Ind_{L_\xi/F}\,\sigma_\xi$. Moreover, 
$$ 
\Psi_\vT(x) = p^{-r} \psi_{L_\xi/F}(x), \quad 0\le x\le c_\vT. 
\tag {\rm D} 
$$ 
\item 
If $\Psi_\vT$ has an odd number of jumps, then $\sigma_\xi $ is irreducible, totally wild, H-singular, of Carayol type and dimension $p^r\big/[L_\xi{:}F] \big.\neq 1$. 
\item   
If $\Psi_\vT$ has an even number of jumps, then $\sigma_\xi$ is a character and $[L_\xi{:}F] = p^r$. 
\endroster 
\endproclaim 
By symmetry, the relation (D) determines $\Psi_\vT$ completely. Any two choices of the character $\xi$ are $\scr W_F$-conjugate, so the same applies to the field $L_\xi$. The field extension $L_\xi/F$ is not usually Galois but, after a suitable tamely ramified base field extension, it has a canonical presentation as a tower of elementary abelian extensions faithfully reflecting the ramification structure of $\sigma$. 
\par 
The canonical presentation of $\sigma$ as an induced representation, 
$$ 
\sigma =  \Ind_{L_\xi/F}\,\sigma_\xi = \Ind_{\scr W_{L_\xi}}^{\scr W_F}\,\sigma_\xi, 
$$ 
is derived from arithmetic considerations. It can claim to be more natural than anything provided by a purely group-theoretic approach. 
\par 
The restrictions $\sigma\Mid \scr R_F(x)$, $\sigma\Mid \scr R_F^+(x)$ follow a clear pattern, underlying the symmetry property of Theorem 1. To give the flavour, suppose there are at least two jumps. Let $j$ be the least and $\bar\jmath$ the greatest. The restriction $\sigma\Mid \scr R_F(j)$ is irreducible, while $\sigma\Mid \scr R_F^+(\bar\jmath)$ is a multiple of a character. The  restriction $\sigma\Mid \scr R^+_F(j)$ is a multiplicity-free direct sum of irreducible representations while $\sigma\Mid \scr R_F(\bar\jmath)$ is a direct sum of characters, its isotypic components being the restrictions of the irreducible components of $\sigma\Mid \scr R^+_F(j)$. The pattern repeats for the second and penultimate jump, and so on. 
\subhead 
7 
\endsubhead  
We now have two expressions, (C) and (D), for the Herbrand function $\Psi_\vT$ of $\vT \in \ewc F$. Together they show how to read the {\it algebraic\/} structure of the decompositions $\sigma\Mid \scr R_F(x)$, $x>0$, directly from the presentation $\vT \in \|\scr C^\star(\frak a,\alpha)\|$. Our final tranche of results treats this in some detail. 
\par 
In the same context, the number $c_\vT$ (as in part 6) and the function $\Psi_\vT$, as $\vT$ ranges over $\|\scr C^\star(\frak a,\alpha)\|$, depend only on  $\alpha$. We therefore denote them by $c_\alpha$ and $\Psi_\alpha$ respectively. Let $j_\infty(\alpha) = j_\infty(F[\alpha]|F)$ be the largest jump of the classical Herbrand function $\psi_{F[\alpha]/F}$. The definition of $\biP{F[\alpha]}F{\vs_\vT}$ and Theorem 2 show that $\Psi_\alpha$ has an even number of jumps if and only if $j_\infty(\alpha) < c_\alpha$. 
\par 
Let $\scr G^\star(\alpha)$ be the set of $\sigma \in \wwr F$ such that $[\sigma]_0^+ \in \uL\|\scr C^\star(\frak a,\alpha)\|$. 
\proclaim{Theorem 4A} 
If $\sigma, \tau \in \scr G^\star(\alpha)$, the representations $\sigma\Mid \scr R_F^+(c_\alpha)$, $\tau\Mid \scr R_F^+(c_\alpha)$ are equivalent. In particular, any character $\xi$ of $\scr R_F^+(c_\alpha)$ occurring in $\sigma\Mid \scr R_F^+(c_\alpha)$ also occurs in $\tau\Mid \scr R_F^+(c_\alpha)$. 
\endproclaim 
All representations $\sigma \in \scr G^\star(\alpha)$ therefore give rise to the same conjugacy class of field extensions $L_\xi/F$ and the associated representations $\sigma_\xi$ all have the same dimension $p^r/[L_\xi{:}F]$. 
\par 
To go further, there is a second field extension to be taken into account. If $\rho\in \wW F$ has dimension $n$, let $\bar\rho: \scr W_F \to \roman{PGL}_n(\Bbb C)$ be the associated projective representation. The kernel of $\bar\rho$ is of the form $\scr W_E$, where $E/F$ is finite and Galois. One calls $E/F$ the {\it centric field\/} of $\rho$. Returning to the main topic, let $\widetilde L_{\sigma,\xi}/L_\xi$ be the centric field of the H-singular representation $\sigma_\xi \in \wwr{L_\xi}$. The extension $\widetilde L_{\sigma,\xi}/L_\xi$ is Galois. It is non-trivial if and only if $\dim\sigma_\xi > 1$, that is, $\biP {F[\alpha]}F{\vs_\vT}$ has an odd number of jumps.  
\par 
Let $w_\alpha = w_{F[\alpha]/F}$ be the wild exponent (1.6.1) of the field extension $F[\alpha]/F$. We divide into two cases. Say that $\alpha$ is {\it $\star$-exceptional\/} if $j_\infty(\alpha) = c_\alpha$ and the integer $l_\alpha = m{-}w_\alpha$ is even and positive. Otherwise, say that $\alpha$ is {\it $\star$-ordinary.} (This terminology is suggested by the usage of \cite{33}, but is not equivalent to it.) 
\par 
For our next result, we fix a character $\xi$ of $\scr R_F^+(c_\alpha)$ occurring in $\sigma\in \scr G^\star(\alpha)$ and abbreviate $L = L_\xi$, $\wt L_\sigma = \wt L_{\sigma,\xi}$. Let $T_\sigma/F$ be the maximal tame sub-extension of $\wt L_\sigma/L$. Let $d_\sigma$ be the number of characters $\chi$ of $\scr W_L$ such that $\phi\otimes \sigma_\xi \cong \sigma_\xi$.  
\proclaim{Theorem 4B} 
\roster 
\item 
If $\alpha$ is $\star$-ordinary, then $\wt L_\sigma = \wt L_\tau$, for all $\sigma,\tau\in \scr G^\star(\alpha)$. 
\item 
If $\alpha$ is $\star$-exceptional, then $T_\sigma = T_\tau$ and $d_\sigma = d_\tau$, for all $\sigma,\tau\in \scr G^\star(\alpha)$. There are, at most, $d_\sigma$ Galois extensions of the form $\wt L_\tau/L$, $\tau \in \scr G^\star(\alpha)$. 
\endroster 
\endproclaim 
The bound in part (2) is achieved when $[T_\sigma{:}L]$ is not divisible by $p$. (In general, we don't know what happens here, but $p$ can divide $[T_\sigma{:}L]$: see 9.6 Example.) In part (1), the set $\scr G^\star(\alpha)$ bears a canonical structure as principal homogeneous space over an easily described group of characters of $L^\times$. 
\subhead 
8 
\endsubhead 
We give an overview of our methods and the layout of the paper. 
\par 
Section 1 is a free-standing account of the classical Herbrand functions $\psi_{E/F}$, $\vf_{E/F}$ of a finite field extension $E/F$. For Galois extensions $E/F$, much of what we need can be deduced from the standard account in \cite{38}. We develop the same level of detail for non-Galois extensions, starting from Deligne's notes \cite{18}. 
\par 
The development proper starts with section 2. We introduce the main players and fix the basic notation. We take a simple stratum $[\frak a,m,0,\alpha]$ in the matrix algebra $\M{p^r}F$, $r\ge 1$, as in part 4 above, and a simple character $\theta\in \scr C(\frak a,\alpha)$ of endo-class $\vT$. Thus $\vT \in \ewc F$ and $\vs_\vT = m/p^r$.  The Interpolation Theorem of \cite{13} readily yields $\Psi_\vT(x) = p^{-r} \psi_{F[\alpha]/F}(x)$ in the range $0\le x\le \vs_\vT/2$. In the region $\vs_\vT/2 < \Psi_\vT(x) \le \vs_\vT$, it interprets the value $\Psi_\vT(x)$ in terms of intertwining properties of certain simple strata. 
\par 
Section 3 is devoted to the proof of Theorem 1. The argument is couched almost entirely in terms of Galois representations. Take $\sigma \in \wwr F$ of dimension greater than $1$. After a tame base field extension, 8.3 Theorem of \cite{13} gives a sufficiently canonical presentation $\sigma = \Ind_{K/F}\,\tau$, where $K/F$ is cyclic of degree $p$. After an elementary change of variables, the jumps of $\vS_\tau$ are among those of $\vS_\sigma$ but one or two of them are ``flattened'', in an obvious sense. One of these is invariably the first. If $\sigma$ is of Carayol type, the other is the last: this follows from an application of the conductor formula of \cite{14} 6.5 Theorem, which also gives a relation between the first and last jumps. One may then assume that $\tau$ has the symmetry property and proceed by induction on dimension. 
\par 
Section 4 makes a transition back to the GL-side. The combination of convexity and symmetry imposes significant restrictions on the piecewise linear graph $y=\Psi_\vT(x)$ in the relevant region $0\le x\le \vs_\vT$. We abstract these properties in the definition of the bi-Herbrand function $\biP EF\vs$. Much of the section is devoted to listing elementary, but useful, geometric properties of the graphs of $\Psi_\vT$ and $\biP EF\vs$. Our strategy is to identify $\Psi_\vT$ as a bi-Herbrand function. In many cases, one can do that straightaway: see 4.6 Example. This simple case also has a role in the more complicated arguments that follow. 
\par 
Sections 5 and 6 are highly technical in nature, preparing the way for the arguments of section 7. In section 5, we use the Interpolation Theorem to identify, via some delicate intertwining and conjugacy arguments, a subset of $\|\scr C(\frak a,\alpha)\|$ on which the Herbrand function $\Psi_\vT$ takes the expected value $\biP{F[\alpha]}F{\vs_\vT}$. The specification of this set, which we temporarily call $\scr L_\alpha$, is quite subtle. There is nothing canonical or natural about $\scr L_\alpha$, but it is a vital computational device. 
\par 
The set $\scr C(\frak a,\alpha)$ does not determine $\alpha$, although it does determine $\frak a$ and the integer $m$. Let $\roman P(\frak a,\alpha)$ be the set of $\beta\in \GL{p^r}F$ for which $[\frak a,m,0,\beta]$ is a simple stratum satisfying $\scr C(\frak a,\beta) = \scr C(\frak a,\alpha)$. In section 6, we examine various ways in which one can construct elements $\beta$ of $\roman P(\frak a,\alpha)$ while keeping track of the relation between the sets $\scr L_\alpha$ and $\scr L_\beta$. 
\par 
In section 7, we first define the subset $\scr C^\star(\frak a,\alpha)$ of simple characters $\theta\in \scr C(\frak a,\alpha)$ that {\it conform to\/} $\alpha$. We show that, if $\theta'\in \scr C(\frak a,\alpha)$, there exists $\alpha' \in \roman P(\frak a,\alpha')$ to which $\theta'$ conforms. The calculations in sections 5 and 6 give a first result (7.2 Theorem 1) from which Theorem 2 follows. 
\par
With section 8, we return to the Galois side. We first re-cast the general theory of representations of, loosely speaking, Heisenberg type and so identify the class of representations with Herbrand function having a single jump. This is in preparation for section 9, where we prove Theorem 3. That result is given in two tranches. In the first (9.2), we assume that $\sigma$ is ``absolutely wild'', in the sense that its centric field extension is totally wildly ramified. The argument there develops the method of section 3. 
\par 
The general case is presented separately as 9.5 Corollary. The transition to the general case is, we found, surprising in both its simplicity and its exactness. It marks a change in direction in the paper. Until the end of section 7, we rely on the fact that, when using the Interpolation Theorem to compute the Herbrand function, one can impose an arbitrary finite, tamely ramified, base field extension while losing no control: the method is illustrated in the proof of 2.6 Proposition and then used repeatedly until the end of the proof of 9.2 Theorem. From 9.5 Corollary on, we have to take account of the tame structures destroyed by such a process. Theorems 4A and 4B follow in section 10, where we combine and compare the main results of sections 7 and 9. 
\par 
Some parts of Theorems 4 are foreshadowed, often in more detail, in the classical literature of dimension $p$ \cite{23}, \cite{32}, \cite{33}, \cite{36}. There is a device from \cite{36} that allows us to remove the distinction between ordinary and exceptional elements $\alpha$, {\it provided $p\neq2$.} We summarize this in 10.6, and then briefly review the historical context. 
\remark{Acknowledgement} 
Both the transparency and accuracy of the exposition were greatly improved by detailed criticism from referees. We are extremely grateful for their efforts. 
\endremark 
\head 
Background and notation 
\endhead 
General notations are quite familiar: $\frak o_F$ is the discrete valuation ring in $F$, $\frak p_F$ is the maximal ideal of $\frak o_F$ and $\ups_F$ is the normalized additive valuation. For $k\ge 1$, $U^k_F$ is the congruence unit group $1{+}\frak p_F^k$. Similarly, if $\frak a$ is a hereditary $\frak o_F$-order in some matrix algebra, then $U^k_\frak a = 1{+}\frak p^k$, where $\frak p$ is the Jacobson radical $\roman{rad}\,\frak a$ of $\frak a$. For real $x$, $x\mapsto [x]$ is the greatest integer function. 
\par 
If $E/F$ is a finite field extension, $\psi_{E/F}$, $\vf_{E/F}$ are the classical Herbrand functions discussed in section 1. If $E/F$ is Galois and $\vG = \Gal EF$, then $\vG_a$, $\vG^a$, $a\ge 0$, are the ramification subgroups of $\vG$ in the lower, upper numbering conventions of \cite{38}. 
The symbols $\scr W_F$, $\wW F$, $\scr P_F$, $\wP F$, $\wG F$, $\Scr E(F)$, $\uL\vT$, $[\sigma]_0^+$, $\scr R_F(x)$, $\scr R_F^+(x)$ all retain the meaning given them in the introduction. Notation concerned with simple characters is all taken from \cite{15} and \cite{3}. For the special cases considered here, full definitions are given in 2.1--3. The broader summary in \cite{2} may be found helpful. Certain special notations recur sporadically. Their definitions may be found as follows: $\vs_\vT$ (2.1), $\vs_\sigma$ (2.2), $\wwr F$ (3.2), $\awr F$ (3.2), $\ewc F$ (2.3), $j_\infty(E|F)$ (1.5), $w_{E/F}$ (1.6),  $ \scr C^\star$ (7.1). 
\head\Rm 
1. Classical Herbrand functions 
\endhead 
Let $E/F$ be a finite, separable field extension. As we go through the paper, we rely on properties of the classical Herbrand function $\psi_{E/F}$ and its inverse $\vf_{E/F}$. For Galois extensions $E/F$, many of these are to be found in \cite{38}. In the general case, we develop them from the outline in \cite{18}. Beyond that, we need estimates of the {\it jumps of\/} $\psi_{E/F}$, that is, the discontinuities of the derivative $\psi'_{E/F}(x)$, $x>0$. With only minor changes, the formalism applies equally well to inseparable extensions $E/F$: we indicate how this is done in 1.7. 
\par 
We conclude the section with what seems to be a novel result on the structure of a broad class of totally ramified extensions. We do not need this until near the end of the paper but it fits well in the present context. The reader may wish to skip that, or even the entire section, referring back to it as needed. 
\subhead 
1.1 
\endsubhead 
Let $E/F$ be a finite Galois extension. The Herbrand function $\psi_{E/F}(x)$ is defined, for $x\ge -1$, in \cite{38} IV \S3 but we shall always assume $x\ge 0$. If $K/F$ is a Galois extension contained in $E$, the fundamental {\it transitivity property\/} $\psi_{E/F} = \psi_{E/K}\circ \psi_{K/F}$ holds. If the finite separable extension $E/F$ is not Galois, we follow \cite{18}. Let $E'/F$ be a finite Galois extension containing $E$. The function $\psi_{E'/E}$ is positive and strictly increasing, so we may set 
$$ 
\psi_{E/F} = \psi_{E'/E}^{-1}\circ \psi_{E'/F}. 
\tag 1.1.1 
$$ 
Because of the transitivity property for Galois extensions, this definition does not depend on the choice of $E'/F$. The relation 
$$ 
\psi_{E/F} = \psi_{E/K} \circ \psi_{K/F} 
\tag 1.1.2 
$$ 
then holds for any tower $F\subset K\subset E$ of finite separable extensions. In all cases, $\vf_{E/F}$ shall be the inverse function for $\psi_{E/F}$, 
$$ 
\vf_{E/F}\circ \psi_{E/F}(x) = x = \psi_{E/F} \circ \vf_{E/F}(x), \quad x\ge 0. 
\tag 1.1.3 
$$ 
\proclaim{Lemma} 
\roster 
\item 
If $K/F$ is finite and tamely ramified, then $\psi_{K/F}(x) = ex$, where $e = e(K|F)$. 
\item 
If $E/F$ is finite separable and $K/F$ is  finite and tamely ramified, with $e(K|F) = e$, then $\psi_{EK/K}(x) = e(EK|E)\,\psi_{E/F}(x/e)$. If $E/F$ is totally wildly ramified, then $\psi_{EK/K}(x) = e\,\psi_{E/F}(x/e)$. 
\endroster 
\endproclaim  
\demo{Proof} 
Part (1) follows immediately from the definitions here and in \cite{38}. By (1.1.2) and part (1),  $\psi_{EK/F}(x) = \psi_{EK/K}\circ \psi_{K/F}(x) = \psi_{EK/K}(ex)$. On the other hand, $\psi_{EK/F}(x) = \psi_{EK/E}\circ \psi_{E/F}(x) = e(EK|E)\psi_{E/F}(x)$, whence part (2) follows. \qed 
\enddemo 
The lemma reduces most questions to the totally wildly ramified case.  
\subhead 
1.2 
\endsubhead 
We list some properties of the graph $y = \psi_{E/F}(x)$, $x\ge 0$. 
\proclaim{Proposition 1} 
Let $E/F$ be a finite separable extension and write $e = e(E|F) = e_0p^r$, where $e_0$ is an integer not divisible by $p$. 
\roster 
\item 
The function $\psi_{E/F}$ is continuous, piece-wise linear, strictly increasing and convex. 
\item If $x$ is sufficiently large, then $\psi'_{E/F}(x) = e$.  
\item There exists $\eps > 0$ such that $\psi_{E/F}(x) = e_0x$, for $0 \le x<\eps$. 
\item The derivative $\psi'_{E/F}$ is continuous except at a finite number of points. 
\endroster 
\endproclaim 
\demo{Proof} 
All assertions are standard when $E/F$ is Galois, and (2)--(4) then follow from (1.1.2) in general. In (1), the first two properties are clear while, by (3), $\psi'_{E/F}(x) = e_0 \ge 1$ for $x$ positive and sufficiently small. It is enough, therefore to show that $\psi_{E/F}$ is convex. By 1.1 Lemma (2), we need only prove that $\psi_{EK/K}$ is convex for some finite tame extension $K/F$. We choose $K/F$ to be the maximal tame sub-extension of the normal closure $E'/F$ of $E/F$. This reduces us to the case in which $E'/F$ is totally wildly ramified. If $E = F$, there is nothing to prove, so assume otherwise. The proper subgroup $\Gal{E'}E$ of the finite $p$-group $\Gal{E'}F$ is contained in a normal subgroup of index $p$. That is, there is a Galois sub-extension $F'/F$ of $E/F$ of degree $p$. In the relation $\psi_{E/F} = \psi_{E/F'} \circ \psi_{F'/F}$, the function $\psi_{F'/F}$ is convex since $F'/F$ is Galois. By induction on degree, $\psi_{E/F'}$ is convex, whence so is $\psi_{E/F}$. \qed 
\enddemo 
This technique of the proof of the proposition will be used again, so we make a formal definition. 
\definition{Definition} 
Let $E/F$ be a finite separable extension, with normal closure $E'/F$. Say that $E/F$ is {\it absolutely wildly ramified\/} if $E'/F$ is totally wildly ramified. 
\enddefinition 
In the notation of the definition, let $K/F$ be the maximal tame sub-extension of $E'/F$. The extension $EK/K$ is then absolutely wildly ramified. From the proof of Proposition 1, we extract a useful property. 
\proclaim{Gloss} 
If $E/F$ is absolutely wildly ramified, there exists a Galois extension $F'/F$, of degree $p$, such that $F'\subset E$. 
\endproclaim 
We give a second application. 
\proclaim{Proposition 2} 
Let $E/F$ be finite, separable and totally wildly ramified. If $\psi_{E/F}$ is smooth at $x$, then the value $\psi'_{E/F}(x)$ is a non-negative power of $p$. 
\endproclaim 
\demo{Proof} 
The result is standard when $E/F$ is Galois. Otherwise, let $K/F$ be finite and tamely ramified. Part (2) of 1.1 Lemma implies that the result holds for $E/F$ if and only if it holds for $EK/K$. It is therefore enough to treat the case of $E/F$ absolutely wild. As in the Gloss, let $F'/F$ be a sub-extension of $E/F$ that is Galois of degree $p$. The extension $F'/F$ has the desired property since it is Galois. By induction on the degree, we may assume that it holds equally for $E/F'$. The proposition then follows from the transitivity relation $\psi_{E/F} = \psi_{E/F'}\circ \psi_{F'/F}$. \qed 
\enddemo 
\subhead 
1.3 
\endsubhead 
As in the Galois case, the function $\psi_{E/F}$ reflects properties of the norm map $\N EF:E^\times \to F^\times$. 
\proclaim{Proposition} 
Let $E/F$ be a finite separable extension. Let $\chi$ be a character of $F^\times$ such that $\sw(\chi) = k\ge 1$. The character $\chi\circ \N EF$ of $E^\times$ then has the properties 
\roster 
\item $\sw(\chi\circ \N EF) \le \psi_{E/F}(k)$ and,  
\item 
if $\psi_{E/F}'$ is continuous at $k$, then $\sw(\chi\circ \N EF) = \psi_{E/F}(k)$. 
\endroster 
\endproclaim 
\demo{Proof} 
The result is standard when $E/F$ is Galois \cite{38} V Proposition 9. 
\par 
Suppose next that $E/F$ is tamely ramified and set $e = e(E|F)$. Thus $\psi_{E/F}(x) = ex$, $x\ge0$. If $\chi$ is a character of $F^\times$ with $\sw(\chi) = k \ge1$, then $\sw(\chi\circ\N EF) = ek$ and there is nothing to prove. 
\par 
Transitivity now reduces us to the case where $E/F$ is totally wildly ramified. Also, if $K/F$ is a finite tame extension, the result holds for $E/F$ if and only if it holds for $EK/K$. We may therefore assume that $E/F$ is absolutely wildly ramified. Let $F'$ be a field, $F\subset F'\subset E$, such that $F'/F$ is Galois of degree $p$ (as in 1.2 Gloss). The result holds for the extension $F'/F$ and so in general, by induction on $[E{:}F]$. \qed 
\enddemo 
\definition{Definition} 
A {\it jump\/} of $\psi_{E/F}$ is a point $x>0$ at which the derivative $\psi'_{E/F}$ is not continuous. Let $J_{E/F}$ denote the set of jumps of $\psi_{E/F}$. 
\enddefinition 
The set $J_{E/F}$ is finite by 1.2 Proposition 1 (4). 
\proclaim{Corollary} 
Let $E/F$ be totally wildly ramified, and let $K/F$ be a finite tame extension, with $e = e(K|F)$. If $\chi$ is a character of $K^\times$ with $\sw(\chi) = k \ge 1$, such that $e^{-1}k \notin J_{E/F}$, then 
$$ 
\sw(\chi\circ \N{EK}K) = \psi_{EK/K}(k) = e\,\psi_{E/F}(e^{-1}k). 
$$ 
\endproclaim 
\demo{Proof}
The second equality is 1.1 Lemma, whence $J_{EK/K} = eJ_{E/F}$. The result now follows from the proposition. \qed 
\enddemo 
\subhead 
1.4  
\endsubhead 
Another familiar property extends to the general case. 
\proclaim{Proposition} 
Let $E/F$ be a finite separable extension. If $\eps > 0$, then 
$$ 
\align 
\scr R_F(\eps)\cap \scr W_E &= \scr R_E(\psi_{E/F}(\eps)), \\ 
\scr R^+_F(\eps)\cap \scr W_E &= \scr R^+_E(\psi_{E/F}(\eps)). 
\endalign 
$$ 
\endproclaim 
\demo{Proof} 
If $E/F$ is Galois, the result follows from \cite{38} IV Proposition 14. The case of $E/F$ tame follows readily. If $K/F$ is a finite tame extension, the result therefore holds for $E/F$ if and only if it holds for $EK/K$ ({\it cf\.} 1.1 Lemma). Thus we need only treat the case where $E/F$ is absolutely wildly ramified. There is a Galois sub-extension $F'/F$ of $E/F$ of degree $p$. If $F' = E$, there is nothing to do, so we assume otherwise. We have 
$$ 
\align 
\scr R_F(\eps) \cap \scr W_E &= \scr R_F(\eps) \cap \scr W_{F'}\cap  \scr W_E \\
&= \scr R_{F'}(\psi_{F'/F}(\eps)) \cap \scr W_E \\ 
&=  \scr R_E(\psi_{E/F'}(\psi_{F'/F}(\eps)) ) = \scr R_E(\psi_{E/F}(\eps)), 
\endalign 
$$ 
by induction on $[E{:}F]$. The second assertion follows. \qed 
\enddemo 
For a sharper result of this kind, see 1.9 Corollary 2 below. 
\subhead 
1.5 
\endsubhead 
Let $j_\infty(E|F)$ be the largest element of $J_{E/F}$. 
\proclaim{Proposition} 
Let $E/F$ be separable and totally wildly ramified. If $\bar E/F$ is the normal closure of $E/F$, then $j_\infty(\bar E|F) = j_\infty(E|F)$. 
\endproclaim 
\demo{Proof} 
Let $K/F$ be a finite tame extension. The result then holds for $E/F$ if and only if it holds for $EK/K$. We may therefore assume that $E/F$ is absolutely wildly ramified. 
\par 
The relation $\psi_{\bar E/F} = \psi_{\bar E/E}\circ \psi_{E/F}$ implies that 
$$ 
J_{\bar E/F} = J_{E/F} \cup \psi_{E/F}^{-1}(J_{\bar E/E}). 
$$ 
We have to show that $j_\infty(E|F)$ is the largest element of this set. Set $\vG = \Gal{\bar E}F$ and $\vD = \Gal{\bar E}E$. The definition of $\vG_x$ \cite{38} IV \S1 gives $\vD_x = \vG_x\cap \vD$, for all $x\ge 0$. Let $k_\infty$ be the largest jump of $\vG$ in this numbering. Thus $\vG_{k_\infty} \neq \{1\} = \vG_{k_\infty+\ve}$, for all $\ve >0$. As $\bar E/F$ is the least Galois extension containing $E$, so $\bigcap_{\gamma \in \vG} \gamma\vD\gamma^{-1} = 1$. That is, $\vD$ has no non-trivial subgroup normal in $\vG$. Since $\bar E/F$ is totally wildly ramified, $\vG_{k_\infty}$ is central in $\vG$, so $\vD_{k_\infty} = \vG_{k_\infty}\cap \vD$ is normal in $\vG$, whence $\vD_{k_\infty} = 1$. The largest jump of $\vD$ is therefore strictly less than $k_\infty$. Translating back, the largest jump $j_\infty(\bar E|E)$ of $\psi_{\bar E/E}$ is strictly less than $\psi_{E/F}(j_\infty(E|F))$. \qed 
\enddemo 
\subhead 
1.6 
\endsubhead 
Let $E/F$ be a finite separable extension. Denote by $d_{E/F}$ the differental exponent of $E/F$: thus $\frak p_E^{d_{E/F}}$ is the different of $E/F$. Define the {\it wild exponent\/} $w_{E/F}$ of $E/F$ by 
$$ 
w_{E/F} = d_{E/F}+1-e(E|F). 
\tag 1.6.1 
$$ 
We record, for use throughout the paper, some basic facts involving the wild exponent. 
\proclaim{Lemma} 
Let $E/F$ be finite, with $E\subset \bar F$. 
\roster 
\item 
If $F\subset K\subset E$, then 
$$ 
w_{E/F} = e(E|K)w_{K/F} + w_{E/K}. 
$$ 
\item 
If $\tau$ is an irreducible representation of $\scr W_E$, then 
$$ 
\sw(\Ind_{E/F}\,\tau) = \big(\sw(\tau) + w_{E/F}\,\dim\tau\big)\,f(E|F). 
$$ 
In particular, 
$$ 
w_{E/F} = \sw\big(\Ind_{E/F}\,1_E\big)/f(E|F), 
$$ 
where $1_E$ is the trivial character of $\scr W_E$. 
\endroster 
\endproclaim 
\demo{Proof} 
Assertion (1) follows from the multiplicativity property of the different and a short calculation. Part (2) follows from the corresponding properties of the Artin exponent \cite{38} Ch\. VI \S2. \qed 
\enddemo 
The main business of the sub-section concerns estimates relating the wild exponent $w_{E/F}$ to the largest jump $ j_\infty(E|F)$ of $\psi_{E/F}$. 
\proclaim{Proposition} 
If $E/F$ is separable and totally wildly ramified of degree $p^r$, then 
$$ 
\psi_{E/F}(x) = p^rx-w_{E/F}, \quad x\ge j_\infty(E|F). 
$$ 
\endproclaim 
\demo{Proof} 
Let $K/F$ be tamely ramified with $e = e(K|F)$. Thus $w_{EK/K} = e\,w_{E/F}$ by the lemma. The result therefore holds for $E/F$ if and only if it holds for $EK/K$. Taking $K/F$ to be the maximal tame sub-extension of the normal closure of $E/F$, we reduce to the case where $E/F$ is absolutely wildly ramified. Part (2) of 1.2 Proposition 1 implies that there is a constant $c_{E/F}$ such that $\psi_{E/F}(x) = p^rx-c_{E/F}$, for $x\ge j_\infty(E|F)$. We show that $c_{E/F} = w_{E/F}$. 
\par
Let $F'/F$ be a sub-extension of $E/F$ that is Galois of degree $p$. In this case, $j_\infty(F'|F)$ is the only jump of $\psi_{F'/F}$, and it equals  $w_{F'/F}/(p{-}1)$ \cite{38} V \S3. The proposition thus holds for $F'/F$. If $E/F$ is Galois, we may assume inductively that $c_{E/F'} = w_{E/F'}$. So, taking $x$ sufficiently large, we get 
$$
\align 
p^rx-c_{E/F} &= \psi_{E/F'}(\psi_{F'/F}(x)) = \psi_{E/F'}(px-w_{F'/F}) \\ 
&= p^rx - p^{r-1}w_{F'/F}-w_{E/F'} = p^rx-w_{E/F},  
\endalign 
$$ 
by the lemma. Thus $c_{E/F} = w_{E/F}$ when $E/F$ is Galois. 
\par 
Suppose that $E/F$ is not Galois. The normal closure $E'/F$ of $E/F$ is totally wildly ramified by hypothesis. So, with $p^s = [E'{:}F]$ and $x$ sufficiently large, we get 
$$ 
\align 
\psi_{E'/F}(x) &= p^sx-w_{E'/F} = \psi_{E'/E}(\psi_{E/F}(x)) \\ 
&= p^{s-r}(p^rx-c_{E/F}) - w_{E'/E}. 
\endalign 
$$ 
Thus $w_{E'/F} = e(E'|E)c_{E/F} - w_{E'/E}$, and the lemma implies $c_{E/F} = w_{E/F}$. \qed 
\enddemo 
\proclaim{Corollary} 
Let $E/F$ be totally wildly ramified of degree $p^r$. If $j_\infty = j_\infty(E|F)$ is the largest jump of $\psi_{E/F}$, then  
$$ 
(p^r{-}1) j_\infty \ge w_{E/F} \ge p^{r-1}(p{-}1) j_\infty \ge p^r j_\infty/2 . 
$$ 
Moreover, $w_{E/F} = (p^r{-}1) j_\infty$ if and only if $j_\infty$ is the only jump of $\psi_{E/F}$. 
\endproclaim 
\demo{Proof} 
Since $\psi_{E/F}(x) \ge x$ for all $x\ge 0$, the first inequality follows directly from the proposition, likewise the final remark. 
\par 
Observe that $\psi_{E/F}'(x) \le p^{r-1}$, for all points $0< x < j_\infty$ at which the derivative is defined (1.2 Proposition 2). The function $\vt(x) = \psi_{E/F}(x){-}p^{r-1}x$ is therefore decreasing on the interval $0<x<j_\infty$. Thus $\vt(j_\infty) \le 0$, or $p^rj_\infty-w_{E/F} \le p^{r-1}j_\infty$, as required. \qed 
\enddemo 
\subhead 
1.7 
\endsubhead 
If $E/F$ is a finite, {\it purely inseparable\/} extension, we set $\psi_{E/F}(x) = x$, $x\ge 0$. If $E/F$ is a finite extension, define 
$$ 
\psi_{E/F} = \psi_{E/E_0}\circ \psi_{E_0/F} = \psi_{E_0/F},  
\tag 1.7.1 
$$ 
where $E_0/F$ is the maximal separable sub-extension of $E/F$. Assuming $E\neq E_0$, the derivative of $\psi_{E/F}$ satisfies $\psi'_{E/F}(x) < [E{:}F]$ for all $x$. We therefore set $j_\infty(E|F) = \infty$ when $E/F$ is not separable. With these definitions, all the results of 1.1--3, 1.5 and 1.6 remain valid. 
\subhead 
1.8 
\endsubhead 
We anticipate a phenomenon arising later on, in sections 5 and 6. 
\par 
Let $E/F$ be totally ramified of degree $p^r$, $r\ge1$. Thus $E = F[\alpha]$, where $\alpha$ is a root of an Eisenstein polynomial $f(X) = X^{p^r} + a_1X^{p^r-1} + \dots + a_{p^r-1}X + a_{p^r} \in \frak o_F[X]$, and one has $d_{E/F} = \ups_E(f'(\alpha))$. 
\par 
Set $a_0 = 1$. If $E/F$ is inseparable, the coefficient $a_j$ is zero unless $j\equiv 0\pmod p$. Each term $(p^r{-}j)a_j\alpha^{j-1}$ in $f'(\alpha)$ vanishes, giving $d_{E/F} = w_{E/F} = \infty$. 
\proclaim{Proposition} 
Suppose $E/F$ is separable and totally ramified of degree $p^r$. There is an integer $k$ such that $0\le k\le p^r{-}1$, and 
$$ 
d_{E/F} = \underset {0\le j \le p^r{-}1} \to {\roman{min}}\, \ups_E\big((p^r{-}j)a_j\alpha^{j-1}\big) \equiv k{-}1 \pmod{p^r}. 
$$ 
In particular, $w_{E/F} \equiv k \pmod p$. If $F$ has characteristic $p$, then $k\not\equiv 0 \pmod p$. 
\endproclaim 
\demo{Proof} 
For $0\le j\le p^r{-}1$, the term $(p^r{-}j)a_j\alpha^{j-1}$ is either zero or 
$$ 
\ups_E\big((p^r{-}j)a_j\alpha^{j-1}\big) \equiv j{-}1 \pmod {p^r}. 
$$
This gives the expression for $d_{E/F}$. If $F$ has characteristic $p$, any term with $j\equiv 0\pmod p$ has valuation $\infty$ and the second assertion follows. \qed 
\enddemo 
If $F$ has characteristic zero, an Eisenstein polynomial $f(X) = X^p{-}a$ gives a field extension $E/F$ of degree $p$ such that $w_{E/F} \equiv 0 \pmod p$. 
\subhead 
1.9 
\endsubhead 
We prove a simple, but under-appreciated, result concerning absolutely wildly ramified extensions $E/F$ (1.2 Definition). It re-appears naturally in the analysis of representations in section 9. 
\par 
Let $E/F$ be a finite separable extension. As before, let $J_{E/F}$ be the set of jumps of the piecewise linear function $\psi_{E/F}$. For $x>0$, define 
$$ 
w_x(E|F) = \underset{\eps\to 0} \to {\text{\rm lim}}\, \psi_{E/F}'(x{+}\eps)\big/\psi'_{E/F}(x{-}\eps). \big. 
$$ 
By 1.2 Proposition 2, $w_x(E|F)$ is a non-negative power of $p$ while $w_x(E|F) > 1$ if and only if $x\in J_{E/F}$. 
\par 
If $E/F$ is a finite Galois extension with $\Gal EF = \vG$, we use the notation $\vG^{y+} = \bigcup_{z>y} \vG^z$, and similarly for the lower numbering.
\proclaim{Proposition} 
Let $E/F$ be separable and absolutely wildly ramified. Let $a$ be the least element of $J_{E/F}$.  
\roster 
\item 
The number $a$ is an integer and there exists a character $\chi$ of $F^\times$ such that $\sw(\chi) = a$ and $\chi\circ \N EF = 1$. 
\item 
Let $D = D_{(1)}(E|F)$ be the group of characters $\chi$ of $F^\times$ such that $\sw(\chi)\le a$ and $\chi\circ\N EF = 1$. All non-trivial elements of $D$ have Swan exponent $a$, and $D$ is elementary abelian of order $w_a(E|F)$. 
\item 
If $E_1/F$ is class field to the group $D$, then $F\subset E_1\subset E$, $\psi_{E_1/F}(a) = a$ and 
$$ 
J_{E/E_1} = \psi_{E_1/F}(J_{E/F}) \smallsetminus \{a\}. 
$$ 
\endroster 
\endproclaim   
\demo{Proof} 
We proceed by induction on $[E{:}F]$. If $[E{:}F] = p$ then, since $E/F$ is absolutely wild, it is Galois and there is nothing to do. Assume, therefore, that $[E{:}F] \ge p^2$. Since $E/F$ is absolutely wild, there is a Galois extension $F'/F$, of degree $p$, contained in $E$ (1.2 Gloss). There is a character $\phi$ of $F^\times$, of order $p$, that vanishes on the group of norms from $F'$. Choose $F'$ so as to minimize $\sw(\phi)$. The integer $c = \sw(\phi)$ is a jump of $\psi_{E/F}$ (1.3 Proposition), so $c\ge a$. We show that $c=a$. 
\par 
Suppose, for a contradiction, that $c>a$. Thus $a = \psi_{F'/F}(a)$ is a jump of $\psi_{E/F'}$ and indeed its least jump. By inductive hypothesis, $a$ is an integer and there is a character $\chi$ of ${F'}^\times$ such that $\chi\circ\N E{F'} = 1$. Since $c>a$, there is a unique character $\chi_1$ of $F^\times$ such that $\chi = \chi_1\circ \N{F'}F$. The character $\chi_1$ has order $p$, while $\sw(\chi_1) = a$ and $\chi_1\circ \N EF = 1$. The extension $F'_1/F$ that is class field to $\chi_1$ has the properties required of $F'/F$ but $\sw(\chi_1) < \sw(\phi)$. This contradicts our hypothesis, and proves (1). 
\par 
In (2), the group $D$ is an abelian $p$-group, since $[E{:}F]$ is a power of $p$. Let $\chi$ be a character of $F^\times$ and suppose that $\sw(\chi) = b$, $1\le b < a$. Since $b\notin J_{E/F}$, $\chi\circ \N EF$ is not trivial by 1.3 Proposition, so $\chi\notin D$. This proves the first assertion in (2). On the other hand, if $\chi\in D$, $\chi\neq 1$, then $\chi^p \in D$ and $\sw(\chi^p) < \sw(\chi)$. Therefore $\chi^p =1$ and it follows that $D$ is elementary abelian. 
\par 
To calculate the order of $D$, we first use part (1) to choose $\chi\in D$, $\chi\neq1$. Let $F'/F$ be class field to $\chi$. In particular, $F'\subset E$ and $F'/F$ is cyclic of degree $p$. The Herbrand function $\psi_{F'/F}$ has one jump, lying at $a$, and $w_a(F'|F) = p$. Composition with $\N {F'}F$ gives a homomorphism $D_{(1)}(E|F) \to D_{(1)}(E|F')$ with kernel of order $p$, generated by $\chi$. The function $\psi_{E/F'}$ has no jump strictly less than $a$, and $w_a(E|F') = p^{-1}w_a(E|F)$. If $w_a(E|F') = 1$, then $D_{(1)}(E|F')$ is trivial whence $D_{(1)}(E|F)$ has order $p =w_a(E|F)$. Assume therefore that $D_{(1)}(E|F)$ has order at least $p^2$, whence $D_{(1)}(E|F')$ has order at least $p$. 
\par 
Let $E'_1/F'$ be class field to the character group $D_{(1)}(E|F')$. Inductively, we can assume that $|D_{(1)}(E|F')| = w_a(E|F')$, so $\psi_{E/E'_1}$ has least jump strictly greater than $a$. If $\vD = \Gal{F'}F$, then $\vD = \vD^a = \vD_a$. Thus $\vD$ acts trivially on $U^1_{F'}/U^{1+a}_{F'}$. It follows that the extension $E_1'/F$ is Galois, of degree $pw_a(E|F') = w_a(E|F)$ and $\psi_{E'_1/F}$ has a unique jump, lying at $a$. Therefore  $\Gal{E'_1}F$ is elementary abelian and class field to a subgroup of $D_{(1)}(E|F)$. Comparing orders, this subgroup is the whole of $D_{(1)}(E|F)$, so $E'_1 = E_1$ and $D_{(1)}(E|F)$ has order $w_a(E|F)$. 
\par 
This completes the proof of (2). 
\par 
We now have 
$$ 
\psi_{E_1/F}(x) = \left\{\, \alignedat3 &x, &\quad &0\le x\le a, \\ &a+p^s(x{-}a), &\quad &a\le x, \endaligned \right. 
\tag 1.9.1 
$$ 
where $p^s = [E_1{:}F] = w_a(E|F)$. The function $\psi_{E/E_1}$ has no jump $j$ such that $j< a$. At $a = \psi_{E_1/F}(a)$, $w_a(E|E_1) = 1 = w_a(E|F)/w_a(E_1|F)$, so $a\notin J_{E/E_1}$. On the other hand, if $b>a$, then $b$ is not a jump of $\psi_{E_1/F}$ and therefore $w_{\psi_{E_1/F}(b)}(E|E_1) = w_b(E|F)$. In other words, $b$ is a jump of $E/F$ if and only if $\psi_{E_1/F}(b)$ is a jump of $E/E_1$. Part (3) follows straightaway. \qed 
\enddemo  
\proclaim{Corollary 1} 
Let $E/F$ be separable and absolutely wildly ramified. Let 
$$ 
j_1<j_2<\cdots\cdots < j_t 
$$ 
be the set of jumps of $\psi_{E/F}$. There is a unique tower of fields 
$$ 
F = E_0 \subset E_1\subset E_2 \subset \cdots\cdots \subset E_t = E 
\tag 1.9.2 
$$ 
with the following properties. 
\roster 
\item 
For $1\le k\le t$, the extension $E_k/E_{k-1}$ is elementary abelian of degree $w_{j_k}(E|F)$. 
\item 
For $1\le k\le t$, the function $\psi_{E_k/E_{k-1}}$ has a unique jump, namely $\psi_{E_{k-1}/F}(j_k)$. 
\endroster 
\endproclaim 
\demo{Proof} 
One applies the proposition to the absolutely wildly ramified extension $E/E_1$ and iterates. \qed 
\enddemo 
We refer to the tower (1.9.2) as the {\it elementary resolution\/} of the absolutely wild extension $E/F$. It gives a factorization 
$$ 
\psi_{E/F} = \psi_{E_t/E_{t-1}}\circ \psi_{E_{t-1}/E_{t-2}} \circ \cdots \circ \psi_{E_2/E_1}\circ \psi_{E_1/F} 
\tag 1.9.3 
$$ 
in which each factor $\psi_{E_k/E_{k-1}}$, $1\le k \le t$, has exactly one jump. 
\par 
We conclude with an application needed in section 10. 
\proclaim{Corollary 2} 
Let $E/F$ be a finite separable extension that is not tamely ramified. If $j_\infty = j_\infty(E|F)$ is the largest jump of $\psi_{E/F}$ then  
$$ 
j_\infty(E|F) = \roman{inf}\,\{x\in \Bbb R: \scr R_F(x)\subset \scr W_E\}. 
$$ 
In particular, $\scr W_E$ contains $\scr R_F^+(j_\infty)$ but not $\scr R_F(j_\infty)$. 
\endproclaim 
\demo{Proof} 
The assertion is unaffected by tamely ramified base field extension, so we may assume that $E/F$ is absolutely wild. We use the notation of Corollary 1 and proceed by induction on the number, $t$ say, of jumps. If $t=1$, then $E = E_1/F$ is elementary abelian with a single jump $j_1 = j_\infty(E|F)$. Every non-trivial character $\chi\in D_{(1)}(E|F)$ has Swan exponent $j_1$ and so is trivial on $\scr R^+_F(j_1)$, but not on $\scr R_F(j_1)$. Since $\scr W_E$ is the intersection of the kernels of all $\chi\in D_{(1)}(E|F)$, the assertion follows. 
\par 
So, we take $t>1$. Inductively we may assume that 
$$ 
\roman{inf}\,\{x: \scr R_{E_1}(x) \subset \scr W_E\} = j_\infty(E|E_1) = \psi_{E_1/F}(j_\infty(E|F)). 
$$ 
For $x > j_1 = \psi_{E_1/F}(j_1)$, we have $\scr R_F(x) = \scr R_{E_1}(\psi_{E_1/F}(x))$ by the first case and 1.4 Proposition. The assertion now follows. \qed 
\enddemo 
\head \Rm 
2. Certain simple characters 
\endhead 
The first part of this section provides a brief {\it aide m\'emoire\/} for those facts and methods from \cite{3}, \cite{13} and \cite{15} that will be used frequently. It relies on parts 2 and 3 of the Introduction for background but is focused on the detail of the special cases with which we are concerned. The later sub-sections 2.4--2.7 give partial results concerning Herbrand functions in those special cases. The notation we set out here remains standard throughout the paper. 
\subhead 
2.1 
\endsubhead 
Let $\Scr E(F)$ be the set of endo-classes of simple characters over $F$. When working with this set, we follow the scheme of \cite{13} 4.2 (apart from one minor adjustment of notation). 
\par
To each $\vT \in \Scr E(F)$, one attaches positive integer invariants $\deg\vT$, $e_\vT$ and a non-negative rational invariant $\vs_\vT$. (In \cite{13}, $\vs_\vT$ is $m_\vT$.) We will never be concerned with the case $\vs_\vT = 0$, so assume $\vs_\vT > 0$. Let $\mu_F$ be a character of $F$ {\it of level one.} By definition, $\mu_F$ is trivial on $\frak p_F$, but not trivial on $\frak o_F$. There exist a simple stratum $[\frak a,m,0,\beta]$ in a matrix ring $\M nF$ and a simple character $\theta\in \scr C(\frak a,0,\beta,\mu_F)$ of endo-class $\vT$. (Here, we have used the full notation of \cite{15} (3.2.1), (3.2.3), but we almost invariably abbreviate it to $\scr C(\frak a,\beta)$.) The algebra $E = F[\beta]$ is a field and 
$$
\deg\vT = [E{:}F], \quad e_\vT = e(E|F),\quad \vs_\vT = m/e_\frak a, 
$$ 
where $e_\frak a$ is the $\frak o_F$-period of the hereditary $\frak o_F$-order $\frak a$. We shall say that $\theta$ is a {\it realization\/} of $\vT$ on $[\frak a,m,0,\beta]$, and that $E/F$ is a {\it parameter field\/} for $\vT$. 
\par 
While $\deg\vT$, $e_\vT$ and $\vs_\vT$ are invariants of $\vT$, there will often be many choices for the field extension $E/F$, even up to isomorphism. The number $\vs_\vT$ has a useful interpretation. If $\pi \in\wG F$ contains a simple character of endo-class $\vT$, then, in the notation of the introduction, $\vs_\vT = \sw(\pi)/\roman{gr}(\pi)$. 
\par 
Let $\sigma \in \wwr F$. Thus $\sigma = \uL\pi$, for some $\pi\in \wG F$. If $\vT$ is the endo-class of a simple character contained in $\pi$, then $\sw(\sigma) = \sw(\pi)$ and  
$$ 
\sw(\sigma)/\dim\sigma = \sw(\pi)/\roman{gr}(\pi) = \vs_\vT. 
\tag 2.1.1 
$$ 
\subhead 
2.2 
\endsubhead 
Attached to $\vT \in \Scr E(F)$ is a {\it structure function\/} $\vF_\vT(x)$, $x\ge 0$, as defined in the Introduction. It is given by the explicit formula (4.4.1) of \cite{13} which we do not need to repeat: for the special cases considered here, see (2.4.1) below. If $\pi \in \wG F$ contains a simple character of endo-class $\vT$, the definition gives 
$$
\vF_\vT(0) = \sw(\check\pi \times \pi)/\roman{gr}(\pi)^2. 
\tag 2.2.1 
$$ 
\indent 
Let $\sigma \in \wW F$. The orbit $[\sigma]_0^+ \in \scr W_F\backslash \wP F$ and the canonical map $\Scr E(F) \to \scr W_F\backslash \wP F$, $\vT \mapsto \uL\vT$, are as in the Introduction. 
\par 
Attached to $\sigma$ is a {\it decomposition function\/} $\vS_\sigma(x)$, $x\ge 0$, defined as follows \cite{13} (3.1.2). Let $\sigma$ act on the vector space $V$, so that the semisimple representation $\check\sigma\otimes\sigma$ acts on $X = \check V\otimes V$. For $\delta > 0$, let $X(\delta)$ be the space of $\scr R_F^+(\delta)$-fixed points in $X$. This has a unique $\scr R^+_F(\delta)$-complement $X'(\delta)$ in $X$. The spaces $X(\delta)$, $X'(\delta)$ provide semisimple, smooth representations of $\scr W_F$. One sets 
$$ 
\vS_\sigma(\delta) = (\dim\sigma)^{-2} \big(\delta\dim X(\delta) + \sw\,X'(\delta)\big). 
\tag 2.2.2
$$ 
The function $\vS_\sigma$ depends only on the orbit $[\sigma]^+_0 \in \scr W_F\backslash \wP F$. 
\par 
Obviously, $\vS_\sigma(0) = \sw(\check\sigma\otimes\sigma)/(\dim\sigma)^2$. Let $\sigma  = \uL\pi$, $\pi \in \wG F$, and let $\vT$ be the endo-class of a simple character contained in $\pi$. Since the Langlands correspondence preserves Swan exponents of pairs, we have 
$$ 
\vS_\sigma(0) = \frac{\sw(\check\sigma\otimes\sigma)}{(\dim\sigma)^2} = \frac{\sw(\check\pi\times\pi)}{\roman{gr}(\pi)^2} = \vF_\vT(0). 
$$ 
\definition{Definition 1} 
Let $\vT \in \Scr E(F)$ and let $\sigma\in \wW F$ satisfy $[\sigma]^+_0 = \uL\vT$. Define the {\it Herbrand function\/} $\Psi_\vT$ of $\vT$ by $\Psi_\vT = \vF_\vT^{-1} \circ \vS_\sigma$. 
\enddefinition 
The function $\Psi_\vT$ is continuous, strictly increasing and piecewise linear. It does not depend on the choice of $\sigma$ in its definition. It satisfies $\Psi_\vT(0) = 0$ and $\Psi_\vT(x) = x$ for $x\ge \vs_\vT$. 
\definition{Definition 2} 
A {\it jump of\/} $\Psi_\vT$ is a point $x$, $0< x < \vs_\vT$, at which $\Psi'_\vT$ is not continuous. 
\enddefinition 
In many cases, the derivative $\Psi'_\vT$ has a discontinuity at $\vs_\vT$, but it holds no interest so we exclude it as a jump. The derivative $\Psi_\vT'$ takes only finitely many values, and the function $\Psi_\vT$ has only finitely many jumps. 
\par 
We often use the following property. Let $K/F$ be a finite, tamely ramified field extension and set $e = e(K|F)$. Let $\vT^K\in \Scr E(K)$ be a $K/F$-lift of $\vT$ \cite{3} 9.7. By 7.1 Proposition of \cite{13}, 
$$ 
\Psi_\vT(x) = \Psi_{\vT^K}(ex)/e, \quad x\ge 0. 
\tag 2.2.3 
$$ 
In Galois-theoretic terms, if $\sigma\in \wW F$ and $[\sigma]^+_0 = \uL\vT$, then $\uL(\vT^K) = [\tau]^+_0 \in \scr W_K\backslash \wP K$, for some irreducible component $\tau$ of $\sigma\Mid \scr W_K$: this follows from 6.2 Proposition of \cite{12}. 
\subhead 
2.3 
\endsubhead 
Let $\vT \in \Scr E(F)$. Say that $\vT$ is {\it totally wild\/} if $\deg\vT = e_\vT = p^r$, for an integer $r\ge0$. So, if $\vT$ is totally wild and if $E/F$ is a parameter field for $\vT$, then $E/F$ is totally ramified of degree $p^r$. If $\vT$ is totally wild and $K/F$ is a finite tame extension, then $\vT$ has a {\it unique\/} $K/F$-lift and that lift is totally wild. 
\par 
Suppose that $\vT\in \Scr E(F)$ is totally wild of degree $p^r$. Say that $\vT$ is {\it of Carayol type\/} if $r\ge1$ and the integer $p^r\vs_\vT$ is not divisible by $p$ ({\it cf\.} \cite{17}). 
\definition{Notation}  
Let $\ewc F$ denote the set of $\vT \in \Scr E(F)$ that are totally wild of Carayol type. 
\enddefinition 
Let $\vT \in \ewc F$ have degree $p^r$. There is a simple stratum $[\frak a,m,0,\alpha]$ in $M = \M{p^r}F$ carrying a realization of $\vT$. We describe this, following the definitions in Chapter 3 of \cite{15}. The integer $m$ is $p^r\vs_\vT$, the field extension $E = F[\alpha]/F$ is totally ramified of degree $p^r$ and $\frak a$ is the unique hereditary $\frak o_F$-order in $M$ that is stable under conjugation by $E^\times$. The integer $m = -\ups_E(\alpha)$ is not divisible by $p$, so the element $\alpha$ is {\it minimal over\/} $F$, in the sense of \cite{15} (1.4.14). We form the group 
$$ 
H^1(\alpha,\frak a) = U^1_E\,U^{1+[m/2]}_\frak a. 
$$ 
Set $\mu_M = \mu_F\circ \roman{tr}_M$, where $\roman{tr}_M:M \to F$ is the matrix trace. Define a function $\mu_M*\alpha$ on $M$ by 
$$ 
\mu_M*\alpha(x) = \mu_M(\alpha(x{-}1)), \quad x\in M. 
\tag 2.3.1 
$$ 
In particular, $\mu_M*\alpha$ represents a character of the group $U^{1+[m/2]}_\frak a$. It is trivial on $U^{1+m}_\frak a$ but non-trivial on $U^m_\frak a$. The set $\scr C(\frak a,\alpha) = \scr C(\frak a,0,\alpha,\mu_M)$ consists of all characters $\vt$ of $H^1(\alpha,\frak a)$ such that $\vt\,\big|\,U^{1+[m/2]}_\frak a = \mu_M*\alpha\,\big|\, U^{1+[m/2]}_\frak a$. By hypothesis, there exists $\theta\in \scr C(\frak a,\alpha)$ of endo-class $\vT$. 
\remark{Remarks} 
\roster 
\item 
The endo-class of any $\vt\in \scr C(\frak a,\alpha)$ is totally wild of Carayol type. 
\item 
Characters $\vt_1,\vt_2\in \scr C(\frak a, \alpha)$ are endo-equivalent if and only if they are equal: this follows from \cite{15} (3.3.2) and is peculiar to this situation. 
\item In the same vein, let $t$ be an integer, $0\le t \le [m/2]$. The restricted characters $\vt_i\Mid H^{1+t}(\alpha,\frak a)$ intertwine if and only if they are {\it equal.} 
\endroster 
\endremark 
In (3), $H^{1+t}(\alpha,\frak a)$ means $H^1(\alpha,\frak a)\cap U^{1+t}_\frak a$. 
\subhead 
2.4 
\endsubhead 
We specialize to the case of $\vT \in \ewc F$. 
\proclaim{Proposition} 
Let $\vT \in \ewc F$ have degree $p^r$. Choose $\sigma \in \wW F$ such that $[\sigma]_0^+ = \uL\vT$. 
\roster 
\item 
The function $\vF_\vT$ satisfies 
$$ 
\vF_\vT(x) = \left\{\, \alignedat3 &\vF_\vT(0) + p^{-r}x, &\quad &0\le x \le \vs_\vT, \\ &x, &\quad &x\le \vs_\vT. \endalignedat \right. 
\tag 2.4.1 
$$ 
\item $\Psi_\vT(0) = 0$ and $\Psi_\vT(x) = x$, for $x\ge \vs_\vT$. 
\item There exists $\ve>0$ such that 
$$ 
\Psi_\vT'(x) = \left\{\,\alignedat3 
&p^{-r},&\quad &0<x<\ve, \\ &p^r,&\quad &\vs_\vT{-}\ve < x < \vs_\vT. \endalignedat \right. 
$$ 
\item The function $\Psi_\vT$ is convex in the region $0<x<\vs_\vT$. 
\item 
If $0<x<\vs_\vT$, then  $0<\Psi_\vT(x) < x$. 
\item 
The jumps of $\Psi_\vT$ are the discontinuities of $\vS'_\sigma(x)$. 
\item 
If $\vs_\vT = m/p^r$ then 
$$ 
\vF_\vT(0) = \vS_\sigma(0) = m(p^r{-}1)/p^{2r}. 
\tag 2.4.2 
$$ 
\endroster 
\endproclaim 
\demo{Proof} 
Part (1) is the definition \cite{13} (4.4.1), and part (2) has already been noted. Part (3) is an instance of \cite{13} 7.6 Proposition. The function $\vS_\sigma$ is convex (2.2.2), and so (4) follows from (1). Part (5) now follows from (4) and (3). Part (6) follows from (1). Part (7) follows from (2.2.1) and \cite{13} 4.1 Proposition. \qed 
\enddemo 
\subhead 
2.5 
\endsubhead 
Key arguments will rely on the Interpolation Theorem of \cite{13} 7.5. We give an overview of that result, as it applies to $\vT \in \ewc F$. 
\definition{Definition}  
A {\it twisting datum\/} over $F$ is a triple $(k,c,\chi)$ in which 
\roster 
\item $k\ge 1$ is an integer; 
\item $c$ is an element of $F$ such that $\ups_F(c) = -k$; 
\item $\chi$ is a character of $F^\times$, of Swan exponent $k$, such that 
$$ 
\chi(x) = \mu_F*c(x), \quad x\in U^{1+[k/2]}_F. 
$$
\endroster 
\enddefinition 
Let $\vT\in \ewc F$ have degree $p^r$. Suppose that $\vT$ is the endo-class of $\theta\in \scr C(\frak a,\alpha)$, exactly as in 2.3. 
If $(k,c,\chi)$ is a twisting datum over $F$, the character $\chi\circ\det$ of $\GL{p^r}F$ satisfies 
$$ 
\chi(\det x) = \mu_M*c(x), \quad x\in U^{1+[p^rk/2]}_\frak a. 
$$ 
Following the discussion in \cite{13} 7.4, the quadruple $[\frak a,m,0,\alpha{+}c]$ is a simple stratum in $M$, such that $H^1(\alpha{+}c,\frak a) = H^1(\alpha,\frak a)$. The character $\chi\theta:x\mapsto \chi(\det x) \theta(x)$, $x\in H^1(\alpha,\frak a)$, lies in $\scr C(\frak a,\alpha{+}c)$. Denote by $\chi\vT$ the endo-class of $\chi\theta$. 
\par 
Let $\Bbb A$ be the ultrametric on $\Scr E(F)$ defined in \cite{13} 5.1 (see also the Notes below). We first give a preliminary version of the result, which follows from \cite{13} 7.3 Proposition. 
\proclaim{Proposition 1} 
Let $k\ge 1$ be an integer that is not a jump of $\Psi_\vT$. If $(k,c,\chi)$ is a twisting datum over $F$, then $\Psi_\vT(k) = \Bbb A(\chi\vT,\vT)$. In particular, $\Bbb A(\chi\vT,\vT)$ depends only on $k$, but not on $c$ or $\chi$. 
\endproclaim 
\remark{Notes} 
\roster 
\item 
In the context of the proposition, $\Bbb A(\chi\vT,\vT) = t/p^r$, where $t$ is the least integer such that the characters $\theta\Mid H^{1+t}(\alpha,\frak a)$, $\chi\theta\Mid H^{1+t}(\alpha,\frak a)$ intertwine in $\GL{p^r}F$: that is the definition of $\Bbb A$ in this case. 
\item 
The characters $\theta\Mid H^{1+t}(\alpha,\frak a)$, $\chi\theta\Mid H^{1+t}(\alpha,\frak a)$ intertwine in $\GL{p^r}F$ if and only if they are conjugate in $\GL{p^r}F$ \cite{15} (3.5.11). If this holds, the conjugation can be implemented by an element of $U^1_\frak a$. 
\item  
When $k$ is a jump of $\Psi_\vT$, $\Bbb A(\chi\vT,\vT)$ may depend on $\chi$, not only on $k$. 
\endroster 
\endremark  
We recall more about the notion of tame lifting, as it applies to $\vT \in \ewc F$. Let $K/F$ be a finite, tamely ramified field extension with $e(K|F) = e$. We form simple characters over $K$ relative to the character $\mu_K = \mu_F\circ \roman{Tr}_{K/F}$ of $K$. There is a unique simple stratum in $\M{p^r}K$ of the form $[\frak a^K,em,0,\alpha]$. Setting $EK = K[\alpha]\i \M{p^r}K$, there is a unique $\theta^K \in \scr C(\frak a^K,\alpha)$ such that $\theta^K(x) = \theta(\N {EK}E(x))$, $x\in U^1_{EK}$. The endo-class $\vT^K$ of $\theta^K$ lies in $\ewc K$ and is the unique $K/F$-lift of $\vT$. Combining Proposition 1 with (2.2.3), we get: 
\proclaim{Proposition 2} 
Let $K/F$ be a finite tame extension with $e = e(K|F)$, and let $\Bbb A_K$ be the canonical ultrametric on $\Scr E(K)$. Let $k\ge 1$ be an integer such that $k/e$ is not a jump of $\Psi_\vT$. If $(k,c,\chi)$ is a twisting datum over $K$, then 
$$ 
\Psi_\vT(k/e) = \Psi_{\vT^K}(k)/e = \Bbb A_K(\chi\vT^K,\vT^K)/e. 
$$ 
\endproclaim  
Proposition 2 summarizes the Interpolation Theorem. 
\subhead 
2.6 
\endsubhead 
Again let $\vT \in \ewc F$ be of degree $p^r$. Choose a simple stratum $[\frak a,m,0,\alpha]$ in $M = \M{p^r}F$ carrying a realization $\theta\in \scr C(\frak a,\alpha)$ of $\vT$ (as in 2.3). We use the Interpolation Theorem to determine $\Psi_\vT$ on half of (the interesting part of) its domain. 
\proclaim{Proposition} 
Writing $E = F[\alpha]/F$, the Herbrand function $\Psi_\vT$ satisfies  
$$ 
\Psi_\vT(x) = p^{-r}\,\psi_{E/F}(x), \quad 0\le x\le \vs_\vT/2,  
$$ 
where $\psi_{E/F}$ is the classical Herbrand function of \rom{1.1, 1.7.}
\endproclaim  
\demo{Proof} 
Let $k$ be an integer, $0<k<\vs_\vT/2$, which is not a jump of either function $\psi_{E/F}$, $\Psi_\vT$. Let $(k,\chi,c)$ be a twisting datum over $F$. The character $\chi\circ\det$ of $\GL{p^r}F$ is trivial on $U^{1+p^rk}_\frak a$. Since $p^rk \le [m/2]$, it is also trivial on the group $U^{1+[m/2]}_\frak a$. The character $\chi\theta:y\mapsto \chi(\det y)\theta(y)$, $y\in H^1(\alpha,\frak a)$, thus lies in $\scr C(\frak a,\alpha)$ (2.3). The characters $\chi\theta$, $\theta$ intertwine on a group $H^{1+t}(\alpha,\frak a) = H^1(\alpha,\frak a)\cap U^{1+t}_\frak a$, $t\ge 0$, if and only if they are equal there (2.3 Remark (3)). So, recalling 2.5 Note 1, $\Bbb A(\vT,\chi\vT) = t/p^r$ where $t$ is the least non-negative integer such that $\chi\circ\det$ is trivial on $U^{1+t}_E$. We have $\chi\circ\det(y) = \chi\circ\N EF(y)$, $y\in E^\times$. That $k$ is not a jump of $\psi_{E/F}$ implies $t = \psi_{E/F}(k)$ (1.3 Proposition) and the result follows from 2.5 Proposition  1 in this case.  
\par
In general, it is enough to prove the desired identity on a dense set of points $x$ satisfying $0 < x < \vs_\vT/2$. Take $x=a/b$, for positive integers $a$ and $b$ with $b$ not divisible by $p$. Assume that $x$ is not a jump of $\psi_{E/F}$ or $\Psi_\vT$. Let $K/F$ be a finite, tamely ramified field extension with $e(K|F) = b$. If $\vT^K$ is the unique $K/F$-lift of $\vT$, then $bx$ is not a jump of $\psi_{EK/K}$ or $\Psi_{\vT^K}$. The first case of the argument, 2.5 Proposition 2 and 1.1 Lemma together yield 
$$ 
\Psi_\vT(x) = \Psi_{\vT^K}(a)/b = p^{-r}\psi_{EK/K}(a)/b = p^{-r} \psi_{E/F}(x), 
$$ 
as required. \qed 
\enddemo 
\remark{Remark} 
In the context of the proposition, there is no reason to demand that $E/F$ be separable. This condition can be imposed, at the cost of a technical argument, but it is easier and more natural to extend the definition of the classical Herbrand function as in 1.7. 
\endremark 
\subhead 
2.7 
\endsubhead 
Remaining in the situation of 2.6, we refine the other part of 2.4 Proposition (3). We use the concept of {\it formal intertwining\/} of strata (as in \cite{15} 2.6). 
\proclaim{Proposition} 
Let $k$ be an integer, $0<k<\vs_\vT$, which is not a jump of $\Psi_\vT$. Let $t = p^r\Psi_\vT(k)$. If $2t > m$, then $t$ is the least integer such that the strata $[\frak a,m,t,\alpha]$, $[\frak a,m,t,\alpha{+}c]$ intertwine formally. 
\endproclaim  
\demo{Proof} 
Let $l$ be an integer such that $2l>m$. We have 
$$ 
\left. \aligned 
\theta(x) &= \mu_M*\alpha(x), \\ \chi\theta(x) &= \mu_M*(\alpha{+}c)(x), \endaligned \right\} \qquad x\in H^{1+l}(\alpha,\frak a) = U^{1+l}_\frak a. 
$$ 
In this situation, an element $g$ of $\GL{p^r}F$ intertwines $\theta\Mid U^{1+l}_\frak a$ with $\chi\theta\Mid U^{1+l}_\frak a$ if and only if $g^{-1}(\alpha+\frak p^{-l})g\cap (\alpha{+}c{+}\frak p^{-l}) \neq \emptyset$, that is, $g$ intertwines the strata $[\frak a,m,l,\alpha]$, $[\frak a,m,l,\alpha{+}c]$ formally. The result so follows from 2.5 Proposition 1. \qed 
\enddemo 
\head \Rm 
3. Functional equation 
\endhead 
Let $\vT \in \ewc F$ (2.3 Notation) be of degree $p^r$. In particular, $r\ge1$. In this section, we uncover a profound and surprising property of the function $\Psi_\vT$. 
\subhead 
3.1 
\endsubhead 
The main result is: 
\proclaim{Theorem} 
Let $\vT \in \ewc F$ be of degree $p^r$, $r\ge 1$. The Herbrand function $\Psi_\vT$ satisfies 
$$ 
\vs_\vT-x = \Psi_\vT\big(\vs_\vT - \Psi_\vT(x)\big), \quad 0\le x\le \vs_\vT. 
\tag 3.1.1 
$$ 
\endproclaim 
For many arguments, it is convenient to have an alternative formulation of (3.1.1). 
\proclaim{Symmetry} 
The function $\Psi_\vT$ satisfies $0\le \Psi_\vT(x) \le x$, for $0\le x\le \vs_\vT$. In that range, the graph $y = \Psi_\vT(x)$ is symmetric with respect to the line $x{+}y = \vs_\vT$. 
\endproclaim 
The first assertion here is 2.4 Proposition (5). Reflection in the line $x{+}y = \vs_\vT$ is the map 
$$ 
\eurm i_{\vs_\vT}: (x,y) \longmapsto (\vs_\vT{-}y,\vs_\vT{-}x),  
$$ 
so the two formulations are indeed equivalent. 
\par 
Before embarking on the proof of (3.1.1), we observe that it has a converse. As recalled in 2.5, the set $\Scr E(F)$ carries a canonical action $(\chi,\vT) \mapsto \chi\vT$ of the group of characters $\chi$ of $U^1_F$. It has the property $\Psi_{\chi\vT} = \Psi_\vT$ \cite{13} 7.4 Proposition. 
\proclaim{Corollary} 
Let $\vT \in \Scr E(F)$ be totally wild, with $\deg\vT = p^r$, $r \ge 1$. Suppose that $\vs_\vT \le \vs_{\chi\vT}$, for all characters $\chi$ of $U^1_F$. The function $\Psi_\vT$ satisfies \rom{(3.1.1)} if and only if $\vT \in \ewc F$. 
\endproclaim  
\demo{Proof} 
The hypothesis on $\vT$ is equivalent to $\vs_\vT = ap^{t-r}$, for integers $a\not\equiv 0 \pmod p$, $0\le t<r$ (7.6 Remark of \cite{13}). In particular, $\vT \in \ewc F$ if and only if $t=0$. By 7.6 Proposition of \cite{13}, there exist $\eps>0$, $\delta >0$, such that 
$$
\Psi_\vT'(x) = \left\{\,\alignedat3 &p^{-r}, &\quad &0<x<\eps, \\ &p^{r-t}, &\quad &\vs_\vT{-}\delta<x<\vs_\vT. \endalignedat \right. 
$$ 
If the functional equation holds for $\vT$, then $t=0$ and so $\vT \in \ewc F$. The converse is the theorem. \qed 
\enddemo 
The proof of (3.1.1) occupies the entire section. The first intermediate result, 3.4 Theorem, is entirely Galois-theoretic and applies to a relatively wide class of representations. The second, 3.5 Theorem, applies only to representations of Carayol type, and its proof depends on an intervention from the GL-side, in the form of a case of the conductor formula of \cite{14}. That result forms the first step in an inductive proof of the theorem above. 
\subhead 
3.2 
\endsubhead 
Let $\sigma \in \wW F$. Let $\vs_\sigma$ be the {\it slope\/} of $\sigma$. That is, 
$$ 
\aligned 
\vs_\sigma &= \roman{inf}\,\{\eps > 0: \scr R_F(\eps) \subset \roman{Ker}\,\sigma\} \\ 
&= \sw(\sigma)/\dim\sigma, 
\endaligned 
\tag 3.2.1 
$$ 
by \cite{22} Th\'eor\`eme 3.5. If $\vs_\sigma > 0$, then $\sigma\Mid \scr R_F(\vs_\sigma)$ does not contain the trivial character. 
\definition{Definition} 
Let $\sigma\in \wW F$. 
\roster 
\item  
Say that $\sigma$ is {\it totally wild\/} if the restriction $\sigma\Mid\scr P_F$ of $\sigma$ to $\scr P_F$ is irreducible. Let $\wwr F$ be the set of totally wild elements $\sigma$ of $\wW F$. Say that $\sigma \in \wwr F$ is {\it of Carayol type\/} if $p$ does not divide $\sw(\sigma)$ and $\dim\sigma \neq 1$. 
\item 
Let $\sigma \in \wwr F$ have dimension $p^r$. Say that $\sigma$ is {\it absolutely wild\/} if the associated projective representation $\bar\sigma:\scr W_F \to \roman{PGL}_{p^r}(\Bbb C)$ factors through a finite Galois group $\Gal EF$, with $E/F$ totally wildly ramified. Write $\awr F$ for the set of absolutely wild elements $\sigma$ of $\wwr F$. 
\endroster 
\enddefinition 
Remark that, if $\sigma \in \wwr F$, then $\dim\sigma = p^r$, for some $r\ge0$. 
\proclaim{Lemma} 
Let $\sigma \in \wwr F$. Let $K/F$ be a finite, tamely ramified field extension and set $e(K|F) = e$. The representation $\sigma^K = \sigma\Mid \scr W_K$ is irreducible. It lies in $\wwr K$ and 
$$ 
\vS_\sigma(x) = e^{-1}\vS_{\sigma^K}(ex), \quad x\ge 0. 
$$ 
One may choose $K/F$ so that $\sigma^K\in \awr K$. 
\endproclaim 
\demo{Proof} 
The relation between decomposition functions is \cite{13} 3.2 Proposition. The projective representation $\bar\sigma$ factors through a finite Galois group $\Gal EF$. The second assertion holds when $K/F$ is the maximal tame sub-extension of $E/F$. \qed 
\enddemo 
\subhead 
3.3 
\endsubhead 
Let $\sigma \in \wwr F$. Directly from the definition recalled in (2.2.2), $\vS_\sigma(x) = x$, for $x> \vs_\sigma{-}\eps$ and some $\eps >0$. Thus all discontinuities of $\vS_\sigma'(x)$ lie in the region $0<x<\vs_\sigma$. We call such points the {\it jumps\/} of $\vS_\sigma$. 
\par 
We assemble some properties of absolutely wild representations. 
\proclaim{Lemma 1} 
Let $\sigma\in \awr F$ have dimension $p^r$, $r\ge 1$. Let $a$ be the least jump of the function $\vS_\sigma$. 
\roster 
\item 
The jump $a$ is an integer and there exists a character $\chi$ of\/ $\scr W_F$, with $\sw(\chi) = a$, such that $\chi\otimes\sigma \cong \sigma$. 
\item 
If $\chi'$ is a non-trivial character of\/ $\scr W_F$ such that $\chi'\otimes \sigma\cong \sigma$, then $\sw(\chi') \ge a$. 
\item 
The character $\chi$ of \rom{(1)} has order $p$. If $K/F$ is the cyclic extension such that $\scr W_K = \roman{Ker}\,\chi$, there exists $\tau \in \awr K$ such that $\sigma \cong \Ind_{K/F}\,\tau$. The representation $\tau$ is uniquely determined up to conjugation by $\Gal KF$. 
\item 
Suppose, in \rom{(3)}, that $r\ge 2$. The representation $\tau$ is then of Carayol type if and only if $\sigma$ is of Carayol type. 
\endroster 
\endproclaim 
\demo{Proof} 
Parts (1)--(3) are \cite{13} 8.3 Theorem. Let $w_{K/F}$ be the wild exponent of the extension $K/F$ (1.6.1). The formula $\sw(\sigma) = \sw(\tau) + \dim(\tau)\,w_{K/F}$ (1.6 Lemma) gives $\sw(\sigma) \equiv \sw(\tau) \pmod p$ and part (4) follows. \qed 
\enddemo 
Continuing in the situation of Lemma 1 we gather some standard facts from section 1 and \cite{38}, for convenience of reference. 
\proclaim{Lemma 2} 
\roster 
\item 
The point $a$ is the unique ramification jump of the extension $K/F$, in either upper or lower numbering. 
\item  
The group $\scr W_K\cap \scr R_F(a)$ is of index $p$ in $\scr R_F(a)$ and $\scr R_F^+(a)\subset\scr W_K$, while $\scr W_F = \scr W_K\scr R_F(a)$. 
\item 
The following relations hold: 
$$ 
\align 
\scr R_K(\eps) &= \left\{\,\alignedat3 &\scr R_F(\eps) \cap \scr W_K, &\quad &0<\eps \le a, \\ 
&\scr R_F(\vf_{K/F}(\eps)), &\quad &a<\eps. \endalignedat \right. \\ 
\scr R_K^+(\eps) &= \scr R^+_F(\vf_{K/F}(\eps)), \quad a \le \eps. 
\endalign 
$$ 
\item 
The Herbrand function $\vf_{K/F}$ is given by
$$ 
\vf_{K/F}(x) = \left\{\,\alignedat3 &x, &\quad &0\le x\le a, \\ &a+(x{-}a)/p, &\quad &a\le x. \endalignedat \right. 
$$ 
\endroster 
\endproclaim 
\subhead 
3.4 
\endsubhead 
As the first part of the proof of 3.1 Theorem, we develop 3.3 Lemma 1 using the same notation. Thus $\sigma \in \awr F$ has dimension $p^r$, $r\ge 1$. The first jump of $\vS_\sigma$ is at $a$, $\chi$ is a character of $\scr W_F$ such that $\sw(\chi) = a$ and $\chi\otimes \sigma \cong \sigma$. Again, $\scr W_K = \roman{Ker}\,\chi$ and $\sigma = \Ind_{K/F}\,\tau$, $\tau\in \awr K$. 
\par 
For $\eps>0$, set 
$$ 
\align 
d_\eps(\sigma) &= \dim\End{\scr R_F(\eps)} \sigma, \\ 
d^+_\eps(\sigma) &= \dim\End{\scr R^+_F(\eps)} \sigma. 
\endalign 
$$ 
Since $\scr R_F(\eps)$, $\scr R_F^+(\eps)$ are normal subgroups of the pro-$p$ group $\scr P_F$, the integers $d_\eps(\sigma)$, $d^+_\eps(\sigma)$ are non-negative powers of $p$. Referring back to the definition (2.2.2) of $\vS_\sigma$, $p^{-2r}d_\eps(\sigma)$ is the {\it left\/} derivative of the piecewise linear function $\vS_\sigma$ at the point $\eps$. Likewise, $p^{-2r}d_\eps^+(\sigma)$ is the {\it right\/} derivative of $\vS_\sigma$ at $\eps$. It follows that $d_\eps(\sigma) = d_\eps^+(\sigma)$ unless $\eps$ is a jump of $\vS_\sigma$. If $\eps$ is a jump of $\vS_\sigma$, then $w_\eps(\sigma) = d^+_\eps(\sigma)/d_\eps(\sigma)$ is a positive power of $p$. Since 
$$ 
\vS_\sigma'(x) = \left\{\, \alignedat3 &p^{-2r}, &\quad &0<x<\delta, \\ &1, &\quad &\vs_\sigma{-}\delta<x, \endalignedat \right. 
$$
for some $\delta>0$, we have 
$$ 
\prod_{\eps>0} w_\eps(\sigma) = \prod_{\eps > 0} d^+_\eps(\sigma)/d_\eps(\sigma) = p^{2r}. 
$$ 
We make parallel definitions, 
$$ 
\align 
\pre Fd_\eps(\tau) &= \dim\End{\scr R_F(\eps)\cap \scr W_K} \tau, \\ 
\pre Fd^+_\eps(\tau) &= \dim\End{\scr R^+_F(\eps)\cap \scr W_K} \tau,  \\ 
\pre Fw_\eps(\tau) &= \pre Fd^+_\eps(\tau)/\pre Fd_\eps(\tau). 
\endalign 
$$ 
The quotient $w_\eps(\sigma)/\pre Fw_\eps(\tau)$ is a power of $p$, and 
$$ 
\prod_{\eps>0} w_\eps(\sigma)/\pre Fw_\eps(\tau) = p^2. 
\tag 3.4.1 
$$ 
\remark{Remark} 
One can define $d_\eps(\tau)$, etc., exactly as before, relative to the base field $K$. One then has $\pre Fd_\eps(\tau) = d_{\psi_{K/F}(\eps)}(\tau)$ ({\it cf\.} 3.3 Lemma 2) and similarly for the other functions. We use the notation $\pre Fd_\eps(\tau)$ to simplify comparison between the two base fields $F$ and $K$.  
\endremark 
We continue with the notation from the start of the sub-section: in particular $\sigma\in \awr F$. We prove: 
\proclaim{Theorem} 
Let $\gamma \in \Gal KF$, $\gamma \neq 1$. The quantity 
$$ 
c = c_{K/F}(\sigma) = \roman{inf}\,\big\{\eps>0: \Hom{\scr R_F(\eps)\cap \scr W_K}\tau{\tau^\gamma} \neq 0\big\} 
\tag 3.4.2 
$$ 
is independent of the choice of $\gamma$. The following properties hold. 
\roster 
\item 
$c \ge a$. 
\item 
If $c > a$, then $w_a(\sigma)/\pre Fw_a(\tau) = w_c(\sigma)/\pre Fw_c(\tau) = p$, while $w_\eps(\sigma)/\pre Fw_\eps(\tau) = 1$ for all other values of $\eps > 0$. 
\item 
If $c = a$, then $w_a(\sigma)/\pre Fw_a(\tau) = p^2$, while $w_\eps(\sigma)/\pre Fw_\eps(\tau) = 1$ for all other values of $\eps > 0$. 
\endroster 
\endproclaim 
\demo{Proof} 
Let $\eps > 0$. The irreducible components of the semisimple representation $\tau\Mid \scr W_K\cap\scr R_F(\eps)$ are all $\scr W_K$-conjugate and occur with the same multiplicity. Likewise for $\tau^\gamma\Mid \scr W_K\cap \scr R_F(\eps)$. Consequently, 
$$ \Hom{\scr R_F(\eps)\cap \scr W_K}\tau{\tau^\gamma}\neq 0 \quad \Longleftrightarrow \quad \tau^\gamma\Mid \scr R_F(\eps)\cap\scr W_K \cong \tau\Mid \scr R_F(\eps)\cap\scr W_K. 
$$ 
This condition is surely independent of $\gamma \neq 1$. If $0<\eps <a$, the function $\vS_\sigma$ is smooth at $\eps$, whence $\sigma\Mid \scr R_F(\eps)$ is irreducible. It is induced from $\tau\Mid \scr W_K\cap \scr R_F(\eps)$ whence follows part (1) of the theorem. 
\par 
To proceed, we need another litany of notation. Let $\eps>0$ and choose an irreducible component $\sigma_\eps$ of $\sigma\Mid \scr R_F(\eps)$. Let $l_\eps(\sigma)$ be the number of distinct $\scr W_F$-conjugates of $\sigma_\eps$, and $m_\eps(\sigma)$ the multiplicity of $\sigma_\eps$ in $\sigma\Mid \scr R_F(\eps)$. Thus $d_\eps(\sigma) = l_\eps(\sigma)m_\eps(\sigma)^2$ while $ l_\eps(\sigma)m_\eps(\sigma)$ is the Jordan-H\"older length of $\sigma\Mid \scr R_F(\eps)$. All of these numbers are non-negative powers of $p$. 
\par 
Similarly, choose an irreducible component $\sigma_\eps^+$ of $\sigma_\eps\Mid \scr R_F^+(\eps)$ and define $l_\eps^+(\sigma)$, $m^+_\eps(\sigma)$ in the same manner. Thus $d^+_\eps(\sigma) = l^+_\eps(\sigma)m^+_\eps(\sigma)^2$ and $ l^+_\eps(\sigma)m^+_\eps(\sigma)$ is the Jordan-H\"older length of $\sigma\Mid \scr R^+_F(\eps)$, all being non-negative powers of $p$. 
\par 
In exactly the same way, let $\tau_\eps$ be an irreducible component of $\tau\Mid \scr W_K\cap \scr R_F(\eps)$ and $\tau^+_\eps$ an irreducible component of $\tau_\eps\Mid \scr W_K\cap \scr R^+_F(\eps)$. We take $\sigma_\eps = \tau_\eps$ for $\eps >a$ and $\sigma_\eps^+ = \tau^+_\eps$ for $\eps \ge a$ ({\it cf\.} 3.3 Lemma 2). 
\proclaim{Lemma 1} 
If $\vS_\sigma$ is smooth at a point $\eps>0$ then $\vS_\tau$ is smooth at $\psi_{K/F}(\eps)$. 
\endproclaim 
\demo{Proof} 
Suppose first that $\eps<a$, so that $\psi_{K/F}(\eps) = \eps$. The definition of $a$ ensures that the function $\vS_\sigma$ is smooth at $\eps$. The representation $\tau$ is irreducible on $\scr R_K(a) = \scr R_F(a)\cap \scr W_K$, and so also on $\scr R_K(\eps)$. It follows that $\vS_\tau$ is smooth at $\eps$. 
\par 
The function $\vS_\sigma$ is not smooth at $a$, so take $\eps > a$. Since $\vS_\sigma$ is smooth at $\eps$, 8.1 Proposition of \cite{13} shows that the representation $\sigma_\eps$ is irreducible on $\scr R_F^+(\eps)$ and that $\sigma_\eps$ is not $\scr W_F$-conjugate to $\chi\otimes \sigma_\eps$, for any character $\chi\neq 1$ of $\scr R_F(\eps)/\scr R_F^+(\eps)$. We have taken $\tau_\eps = \sigma_\eps$, so $\tau_\eps$ is irreducible on $\scr R_F^+(\psi_{K/F}(\eps)) = \scr R_F^+(\eps)$, and it is not $\scr W_K$-conjugate to $\tau_\eps\otimes\phi$, for any non-trivial character $\phi$ of $\scr R_F(\psi_{K/F}(\eps))/\scr R_F^+(\psi_{K/F}(\eps)) = \scr R_F(\eps)/\scr R_F^+(\eps)$. Therefore $\vS_\tau$ is smooth at $\psi_{K/F}(\eps)$, as required. \qed 
\enddemo 
We assume henceforth that $\eps>a$ and use the notation introduced for Lemma 1. The $\scr W_F$-stabilizer of (the isomorphism class of) $\sigma_\eps$ is of the form $G_\eps = \scr W_{E_\eps}$, for a finite field extension $E_\eps/F$. Likewise, let $G_\eps^+ = \scr W_{E_\eps^+}$ denote the $\scr W_F$-stabilizer of $\sigma_\eps^+$. The $\scr W_K$-stabilizer of $\tau_\eps = \sigma_\eps$ is then $\scr W_K\cap G_\eps = \scr W_{KE_\eps}$, and similarly with $+$'s. 
\proclaim{Lemma 2} 
If $\eps > a$, then 
$$ 
\frac{d^+_\eps(\sigma)}{d_\eps(\sigma)} = \frac{\pre Fd^+_\eps(\tau)}{\pre Fd_\eps(\tau)}\,\frac{[K\cap E_\eps:F]}{[K\cap E^+_\eps:F]}\,. 
$$ 
The quotient of field degrees takes only the values $1$ and $p$. 
\endproclaim 
\demo{Proof} 
Since $\eps > a$, 
$$ 
\align 
m_\eps(\sigma) &= \sum_{\gamma\in \Gal KF} \dim\Hom{\scr R_F(\eps)} {\sigma_\eps}{\tau^\gamma} \\ 
&= \sum_{\gamma\in \Gal KF} \dim\Hom{\scr R_F(\eps)} {\sigma^\gamma_\eps}{\tau} . 
\endalign 
$$ 
If $\sigma_\eps^\gamma$ occurs in $\tau$, then $\sigma_\eps^\gamma = \sigma_\eps^\delta$, for some $\delta \in \scr W_K$, and conversely. The sum is therefore effectively taken over $\gamma\in \scr W_K\scr W_{E_\eps}/\scr W_K = \Gal K{K\cap E_\eps}$, so 
$$ 
m_\eps(\sigma) = \pre Fm_\eps(\tau)\,p\big/[K\cap E_\eps:F].  
$$ 
By definition, $l_\eps(\sigma) = [E_\eps:F]$ and $\pre Fl_\eps(\tau) = [KE_\eps:K] = [E_\eps:F]/[K\cap E_\eps:F]$. That is, 
$$ 
l_\eps(\sigma) = \pre Fl_\eps(\tau)\,[K\cap E_\eps:F]. 
$$ 
Consequently, 
$$ 
d_\eps(\sigma) = \pre Fd_\eps(\tau)\,p^2/[K\cap E_\eps:F],  
$$ 
and, likewise, 
$$ 
d^+_\eps(\sigma) = \pre Fd^+_\eps(\tau)\,p^2/[K\cap E^+_\eps:F]. 
$$ 
This proves the first assertion of the lemma. 
\par 
As $[K{:}F] = p$, the quotient $[K\cap E_\eps{:}F]/[K\cap E_\eps^+{:}F]$ may take only the values $1$, $p^{\pm1}$. It remains to show that the case $[K\cap E_\eps{:}F]/[K\cap E_\eps^+{:}F] = p^{-1}$ cannot arise. In other words, we have to show that $K\cap E_\eps = F$ implies $K\cap E_\eps^+ = F$. 
\par
Suppose, therefore, that $K\cap E_\eps = F$ or, as amounts to the same, $G_\eps\scr W_K = \scr W_F$. The restriction of $\tau$ to $\scr R_F(\eps)$ is a multiple of $\sum_\delta \sigma_\eps^\delta$, with $\delta$ ranging over $G_\eps\cap \scr W_K\backslash \scr W_K$, while $\sigma\Mid  \scr R_F(\eps)$ is a multiple of $\sum_\beta \sigma_\eps^\beta$, with $\beta\in G_\eps\backslash \scr W_F$. Our hypothesis $K\cap E_\eps = F$ implies that the natural map $G_\eps\cap \scr W_K\backslash \scr W_K \to G_\eps\backslash \scr W_F$ is bijective. We conclude that $\sigma\Mid  \scr R_F(\eps) = p\,\tau\Mid  \scr R_F(\eps)$, whence $\sigma\Mid  \scr R^+_F(\eps) = p\,\tau\Mid  \scr R^+_F(\eps)$. Put another way, 
$$ 
\frac{d^+_\eps(\sigma)}{d_\eps(\sigma)} = \frac{\pre Fd^+_\eps(\tau)}{\pre Fd_\eps(\tau)},  
$$ 
so $K\cap E_\eps^+ = F$, as required. \qed 
\enddemo 
For $c$ as in (3.4.2), observe that 
$$ 
\Hom{\scr R_F(c)\cap \scr W_K}\tau{\tau^\gamma} = 0. 
\tag 3.4.3 
$$ 
Otherwise, the representation $\check\tau\otimes\tau^\gamma$ would have an irreducible component $\lambda$ for which $\roman{Ker}\,\lambda$ contained $\scr R_F(c) \cap \scr W_K = \scr R_K(c')$, where $c' = \psi_{K/F}(c)$. In that case, $\roman{Ker}\,\lambda$ would contain $\scr R_K(c'')$, for some $c'' < c'$ (\cite{13} 2.1 Proposition 1). That is, $\Hom{\scr R_K(c'')}\tau{\tau^\gamma} \neq 0$, contrary to the definition of $c$. 
\proclaim{Lemma 3} 
If $\phi > c$, then $w_\phi(\sigma)/\pre Fw_\phi(\tau) = 1$. If $c > a$, then $w_c(\sigma)/\pre Fw_c(\tau) = p$. 
\endproclaim 
\demo{Proof} 
Let $\phi > c$, so that $\Hom{\scr R_F(\phi)}\tau{\tau^\gamma} \neq 0$. It follows that $\tau$ is $\scr R_F(\phi)$-isomorphic to $\tau^\gamma$, for all choices of $\gamma$. Therefore $\sigma\Mid \scr R_F(\phi)$ is a sum of $p$ copies of $\tau\Mid \scr R_F(\phi)$ and so $\sigma\Mid \scr R^+_F(\phi)$ is a sum of $p$ copies of $\tau\Mid \scr R^+_F(\phi)$. This implies $w_\phi(\sigma) = \pre Fw_\phi(\tau)$. 
\par
Suppose $c > a$. We have $\Hom{\scr R_F(c)\cap \scr W_K}\tau{\tau^\gamma} = 0$ while $\Hom{\scr R^+_F(c)}\tau{\tau^\gamma} \neq 0$. The second property implies that $G_c^+\scr W_K = \scr W_F$, whence $K\cap E_c^+ = F$ (notation as in the proof of Lemma 2). The first property implies $G_c\scr W_K \neq \scr W_F$, giving $K\subset E_c$. From Lemma 2, we deduce that $w_c(\sigma)/\pre Fw_c(\tau) = p$. \qed 
\enddemo 
Consider now the situation at the point $a$. 
\proclaim{Lemma 4} 
Let $\gamma$ generate $\Gal KF$. 
\roster 
\item 
If\/ $\Hom{\scr R_F^+(a)}\tau{\tau^\gamma} = 0$, then $c>a$ and $w_a(\sigma)/\pre Fw_a(\tau) = p$. 
\item  
If\/ $\Hom{\scr R_F^+(a)}\tau{\tau^\gamma} \neq 0$, then $c=a$ and $w_a(\sigma)/\pre Fw_a(\tau) = p^2$. 
\endroster 
\endproclaim 
\demo{Proof} 
The representation $\sigma\Mid  \scr R_F(a)$ is irreducible and 
$$ 
\align 
\sigma\Mid  \scr R_F(a) &= \sum_{x\in \scr W_K\backslash\scr W_F/\scr R_F(a)} \Ind_{\scr W_K\cap \scr R_F(a)}^{\scr R_F(a)} \,\tau^x\Mid  (\scr W_K\cap \scr R_F(a)) \\ 
 &= \Ind_{\scr W_K\cap \scr R_F(a)}^{\scr R_F(a)}\,\tau\Mid  (\scr W_K\cap \scr R_F(a)). 
\endalign 
$$ 
It follows that $\tau$ is irreducible on $\scr R_K(a) = \scr R_F(a)\cap \scr W_K$, and that the representations $\tau^\gamma\Mid \scr R_K(a)$, $\gamma \in \scr W_K\backslash \scr W_F$, are distinct. 
\par 
Next, 
$$ 
\align 
\sigma\Mid  \scr R^+_F(a) &= \sum_{x\in \scr W_K\backslash\scr W_F/\scr R^+_F(a)} \Ind_{\scr W_K\cap \scr R^+_F(a)}^{\scr R^+_F(a)}\, \tau^x\Mid  (\scr W_K\cap \scr R^+_F(a)) \\ 
&= \sum_{\gamma\in \scr W_K\backslash \scr W_F} \tau^\gamma\Mid  \scr R_F^+(a). 
\endalign 
$$ 
The restrictions $\tau^\gamma\Mid  \scr R_F^+(a)$ are either disjoint or identical. If they are disjoint, then 
$$ 
l^+_a(\sigma) = \pre Fl^+_a(\tau)p \quad \text{and} \quad m^+_a(\sigma) = \pre Fm^+_a(\tau). 
$$ 
In this case, 
$$ 
d^+_a(\sigma) = \pre Fd^+_a(\tau)p \quad \text{and} \quad\tau^\gamma\Mid  \scr R_F^+(a) \not\cong \tau\Mid  \scr R_F^+(a),\ \gamma\neq 1. 
$$ 
If the $\tau^\gamma\Mid  \scr R_F^+(a)$ are identical, then 
$$ 
l^+_a(\sigma) = \pre Fl^+_a(\tau), \quad m^+_a(\sigma) = \pre Fm^+_a(\tau)p, 
$$ 
yielding 
$$ 
d^+_a(\sigma) = \pre Fd^+_a(\tau)p^2 \quad \text{and} \quad \tau^\gamma\Mid  \scr R_F^+(a) \cong \tau\Mid  \scr R_F^+(a). 
$$ 
Since $d_a(\sigma) = \pre Fd_a(\tau) = 1$, the lemma follows. \qed 
\enddemo 
We prove the theorem. Part (1) has been done. Part (2) is given by Lemma 3, Lemma 4(1) and (3.4.1). Part (3) follows from (3.4.1) and Lemma 4(2). \qed 
\enddemo 
\subhead 
3.5 
\endsubhead 
We continue in the situation of 3.4, except that we now specialize to representations of Carayol type. Take $K/F$ and $c_{K/F}(\sigma)$ as in 3.4 Theorem. 
\proclaim{Theorem} 
Let $\sigma\in \awr F$ be of Carayol type and dimension $p^r$. Let $a_\sigma$ be the least jump of the function $\vS_\sigma$. The largest jump $z_\sigma$ of $\vS_\sigma$ is then 
$$
z_\sigma = c_{K/F}(\sigma) = \frac{\sw(\sigma)-a_\sigma}{p^r}. 
$$ 
\endproclaim 
\demo{Proof} 
We proceed by induction on $r$. Take $r=1$. We then have $\vS_\sigma(0) = (p{-}1)\sw(\sigma)/p^2$ (2.4.2) and $\vS_\sigma(x) = x$ for $x\ge \vs_\sigma = \sw(\sigma)/p$. In particular, $0<a_\sigma\le z_\sigma<\vs_\sigma$. In the region $0<x<\vs_\sigma$, the derivative $\vS_\sigma'(x)$ takes the values $p^{-2}$, $1$ and, possibly, $p^{-1}$ (as follows from (2.2.2)). If only the values $p^{-2}$, $1$ occur, then $a_\sigma$ is the only jump. It lies at the intersection of the lines $y= p^{-2}x+(p{-}1)\sw(\sigma)/p^2$ and $y= x$, that is $a_\sigma = \sw(\sigma)/(1{+}p^r) = (\sw(\sigma){-}a_\sigma)/p^r$, as required. If, on the other hand, $\vS_\sigma'$ takes the value $p^{-1}$ on some interval, then $z_\sigma$ is given by the intersection of the lines $y=x$ and $y{-}\vS_\sigma(a_\sigma) = (x{-}a_\sigma)/p$. Since $\vS_\sigma(a_\sigma) = p^{-2}a_\sigma+\vS_\sigma(0)$, the result follows from a quick calculation. 
\par 
Assume $r\ge 2$. From 3.3 Lemma 1 we recall: 
\proclaim{Lemma 1} 
The representation $\tau$ is absolutely wild of Carayol type and dimen\-sion $p^{r-1}$. 
\endproclaim 
We may therefore assume inductively that 
$$ 
z_\tau = (\sw(\tau){-}a_\tau)/p^{r-1}, 
$$ 
where $a_\tau \le z_\tau$ are the first and last jumps of $\vS_\tau$. We calculate a list of Swan exponents. 
\proclaim{Lemma 2} 
\roster 
\item 
$\sw(\check\sigma\otimes \sigma) = (p^r{-}1)\,\sw(\sigma)$. 
\item $\sw(\check\tau\otimes\tau) = (p^{r-1}{-}1)\,\sw(\tau)$. 
\item 
If $\gamma$ generates $\Gal KF$, then $\sw(\check\tau\otimes \tau^\gamma) = p^{r-1}(\sw(\tau){-}a_\sigma)$. 
\endroster 
\endproclaim 
\demo{Proof}
The representations $\sigma$, $\tau$ are of Carayol type, so (1) and (2) are given by (2.4.2) and (2.1.1). As in \cite{13} (2.5.3), set 
$$ 
\Delta_K(\rho_1,\rho_2) = \roman{inf}\,\{x>0: \Hom{\scr R_K(x)}{\rho_1}{\rho_2} \neq 0\}. 
$$ 
Thus (\cite{13} (3.1.4))
$$ 
\frac{\sw(\check\rho_1\otimes \rho_2)}{\dim(\rho_1)\,\dim(\rho_2)} = \vS_{\rho_1}(\Delta_K(\rho_1,\rho_2)), \quad \rho_i \in \wW K. 
$$ 
We started the proof of 3.4 Theorem by observing that, in effect, $\Delta_K(\tau,\tau^\gamma)$ is independent of $\gamma\in \Gal KF$, $\gamma \neq 1$. It follows that $\sw(\check\tau\otimes\tau^\gamma)$ does not depend on $\gamma$. With this in mind, we apply the 
induction formula for the Swan conductor (1.6 Lemma) to the relations 
$$ 
\align 
\check\tau \otimes \sigma\Mid \scr W_K &= \sum_{\gamma \in \Gal KF} \check\tau \otimes\tau^\gamma, \\ 
 \check\sigma \otimes\sigma &= \Ind_{K/F}\big(\check\tau \otimes \sigma\Mid \scr W_K\big).  
\endalign 
$$
By 1.6 Proposition, $w_{K/F} = (p{-}1)a_\sigma$. So, for any $\gamma \neq 1$, 
$$ 
(p{-}1)\,\sw(\check\tau\otimes \tau^\gamma) = \sw(\check\sigma\otimes \sigma) - \sw(\check\tau\otimes \tau) - p^{2r-1}(p{-}1)a_
\sigma,  
$$ 
whence (3) follows. \qed 
\enddemo 
\remark{Remark} 
The formul\ae\ in parts (1) and (2) of Lemma 2 rely ultimately on the conductor formula of \cite{14}. This is the only intervention from the 
GL-side in the proofs of the theorems of 3.4 and 3.5. It is, however, crucial. 
\endremark 
The definition of $c = c_{K/F}$ in (3.4.2) gives $\psi_{K/F}(c) = \Delta_K(\tau,\tau^\gamma)$. Since $c\ge a_\sigma$ (3.4 Theorem (1)), 
we have $\psi_{K/F}(c) = a_\sigma+p(c{-}a_\sigma)$. 
\proclaim{Lemma 3} 
If $\gamma \in \Gal KF$, $\gamma \neq 1$, then $\Delta_K(\tau,\tau^\gamma) \ge z_\tau$. Equality holds here if and only if $a_\sigma = a_\tau$. 
\endproclaim 
\demo{Proof} 
The relation $\vS_\tau(\Delta_K(\tau,\tau^\gamma)) = p^{2-2r}\sw(\check\tau\otimes\tau^\gamma)$ reduces us to proving 
$$ 
\sw(\check\tau\otimes\tau^\gamma) \ge p^{2r-2}\vS_\tau(z_\tau). 
$$ 
Since $z_\tau$ is the last jump of $\vS_\tau$, we have $\vS_\tau(y) = y$, for $y>z_\tau$. In particular, $\vS_\tau(z_\tau) = z_\tau$. The inductive hypothesis therefore yields 
$$ 
p^{2r-2}\vS_\tau(z_\tau) = p^{r-1}(\sw(\tau)-a_\tau). 
$$ 
On the other hand, $\sw(\check\tau\otimes \tau^\gamma) = p^{r-1}\sw(\tau){-}p^{r-1}a_\sigma$ by Lemma 2(3). By 3.4 Lemma 1, we have $a_\sigma \le a_\tau$ whence the result follows. \qed 
\enddemo 
\proclaim{Lemma 4} 
The element $c = c_{K/F}(\sigma)$ satisfies $c = z_\sigma \ge \vf_{K/F}(z_\tau)$. 
\endproclaim 
\demo{Proof} 
By definition, the number $\vf_{K/F}(z_\tau)$ is the infimum of $\eps > 0$ such that $\tau\Mid \scr W_K\cap \scr R_F(\eps)$ is a multiple of a character. Only numbers $\eps > a_\sigma$ enter and, by 3.3 Lemma 2, $\scr R_F(\eps) \subset \scr W_K$ for such $\eps$. That is, $\vf_{K/F}(z_\tau)$ is the infimum of $\eps > 0$ such that $\tau\Mid \scr R_F(\eps)$ is a multiple of a character. Lemma 3 gives 
$$
c =\vf_{K/F}(\Delta_K(\tau,\tau^\gamma)) \ge \vf_{K/F}(z_\tau) 
\tag 3.5.1 
$$ 
while, on the other hand, $c$ is the infimum of numbers $\eps$ such that $\tau\Mid \scr R_F(\eps) \cong \tau^\gamma \Mid \scr R_F(\eps)$. Thus (3.5.1) implies that $c$ is the infimum of numbers $\eps$ such that $\sigma\Mid  \scr R_F(\eps)$ is a multiple of a character. That is,  $c = z_\sigma \ge \vf_{K/F}(z_\tau)$, as required. \qed 
\enddemo 
Lemma 4 yields the first assertion of the theorem. We prove the second. To complete the induction, we have to show that 
$$ 
c = z_\sigma = p^{-r}(\sw(\sigma)-a_\sigma). 
$$ 
Abbreviating $\Delta = \Delta_K(\tau,\tau^\gamma)$, (3.5.1) asserts that  
$$ 
\psi_{K/F}(c) = a_\sigma+p(c{-}a_\sigma)  = \Delta. 
\tag 3.5.2 
$$ 
We have $\vS_\tau(y) = y$, for $y\ge z_\tau$, while Lemma 3 gives $\Delta \ge z_\tau$. So, 
$$ 
\Delta = \vS_\tau(\Delta) = \sw(\check\tau\otimes\tau^\gamma)/p^{2r-2} = p^{1-r}(\sw(\tau)-a_\sigma). 
$$ 
Combining with (3.5.2), we get 
$$ 
p^rc = \sw(\tau)+(p^r{-}p^{r-1}{-}1)a_\sigma. 
$$ 
However, $\sw(\tau) = \sw(\sigma)- p^{r-1}(p{-}1)a_\sigma$, whence 
$$ 
z_\sigma = c = p^{-r}(\sw(\sigma)-a_\sigma), 
\tag 3.5.3 
$$ 
as required. \qed 
\enddemo 
Keeping the notation of the theorem, we exhibit a consequence. 
\proclaim{Corollary 1} 
Let $\sigma \in \awr F$ be of Carayol type and degree $p^r$, $r\ge1$. Set $a = a_\sigma$. If $w_a(\sigma)/\pre Fw_a(\tau) = p^2$, then $a$ is the unique jump of the function $\vS_\sigma$. 
\endproclaim 
\demo{Proof} 
Lemma 4(2) of 3.4 implies $c = a_\sigma$. We have just shown that $c = z_\sigma$. The function $\vS_\sigma$ thus has a unique jump. \qed 
\enddemo 
\remark{Remark}  
The conclusion of the corollary has strong implications for the structure of the representation $\sigma$: see 8.4 Proposition below. 
\endremark 
To finish, we note that, because of (2.2.3), the theorem and its corollary apply equally to totally wild representations that are not absolutely wild. In particular, 
\proclaim{Corollary 2} 
Let $\sigma \in \wwr F$ be of Carayol type and dimension $p^r$. If $a_\sigma$ and $z_\sigma$ are the first and last jumps of the function $\vS_\sigma$ respectively, they are related by 
$$
z_\sigma = \frac{\sw(\sigma)-a_\sigma}{p^r}. 
$$ 
\endproclaim 
\subhead 
3.6 
\endsubhead 
We start the proof of the functional equation (3.1.1). The argument occupies the rest of the section. 
\par 
In 3.4, 3.5, we effectively worked with decomposition functions. We must now pass to Herbrand functions. To avoid the need for more notation, we work with endo-classes. Nontheless, the underlying technique is entirely Galois-theoretic and could be phrased in those terms. We start with the necessary translation. 
\proclaim{Proposition} 
Let $\vT \in \ewc F$ be of degree $p^r$, $r\ge1$. If $a_\vT\le z_\vT$ are the first and last jumps of $\Psi_\vT$, then 
$$ 
z_\vT = \vs_\vT-a_\vT/p^r = \vs_\vT-\Psi_\vT(a_\vT). 
\tag 3.6.1 
$$ 
\endproclaim 
\demo{Proof} 
There exists an irreducible cuspidal representation $\pi$ of $\GL{p^r}F$ that contains a simple character of endo-class $\vT$. The representation $\sigma = \uL\pi$ is therefore totally wild of dimension $p^r$. Moreover, $\sw(\sigma) = p^r\vs_\sigma = p^r\vs_\vT$ is not divisible by $p$, so $\sigma$ is of Carayol type. The formula in part (1) of 2.4 Proposition implies that the functions $\Psi_\vT$, $\vS_\sigma$ have the same jumps. In particular, $a_\vT = a_\sigma$ and $z_\vT = z_\sigma$. The first equality in (3.6.1) thus follows from 3.5 Corollary 2. In the range $0<x<a_\vT$, we have $\Psi_\vT'(x) = p^{-r}$ and so $a_\vT/p^r = \Psi_\vT(a_\vT)$, as required for the second equality. \qed 
\enddemo 
\subhead 
3.7 
\endsubhead 
Let $\vT \in \Scr E(F)$ be totally wild. Say that $\vT$ is {\it absolutely wild\/} if there exists $\sigma\in \awr F$ such that $\uL\vT = [\sigma]_0^+$. The relation $[\sigma]_0^+ = \uL\vT$ determines $\sigma$ up to tensoring with a tame character of $\scr W_F$ \cite{12} 1.3 Proposition. So, if one choice of $\sigma$ is absolutely wild, then all are. 
\par 
For given $\vT$, there surely exists a finite tame extension $T/F$ so that the unique $T/F$-lift $\vT^T$ of $\vT$ is absolutely wild. We have $\vs_{\vT^T} = e(T|F)\vs_\vT$. From (2.2.3) we deduce that if (3.1.1) holds for $\vT^T$ it also holds for $\vT$. We therefore proceed on the basis that the given endo-class $\vT$ is {\it absolutely wild.} 
\par 
For the next result, take $\vT\in \ewc F$, absolutely wild of degree $p^r$. Choose $\sigma \in \awr F$ so that $[\sigma]^+_0 = \uL\vT$. Define $a = a_\sigma$, $K/F$ and $\tau$, relative to $\sigma$, as in 3.3 Lemma 1. Let $c = c_{K/F}(\sigma)$ as in (3.4.2), and note that $a = a_\vT$. 
\proclaim{Proposition} 
There exists a unique $\vU \in \Scr E(K)$ such that $[\tau]^+_0 = \uL\vU$. If $r\ge 2$, the endo-class $\vU$ is absolutely wild of degree $p^{r-1}$ while, otherwise, $\deg\vU  =1$. In either case, it satisfies 
$$ 
\Psi_\vT(x) = p^{-1}\Psi_\vU(\psi_{K/F}(x)), \quad 0\le x\le c. 
$$ 
\endproclaim 
\demo{Proof} 
The existence and uniqueness of $\vU$ are clear. If $r\ge 2$, then $\tau$ is absolutely wild, whence so is $\vU$. In the region $0\le x\le a$, we have $\Psi_\vT(x) = p^{-r}x$ while $\Psi_\vU(\psi_{K/F}(x)) = \Psi_\vU(x) = p^{1-r}x$. The required relation therefore holds in this range. In the case $a=c$, there is nothing left to do so we assume $a<c$. 
\par 
If $a<x<c$, 3.4 Theorem gives $w_x(\sigma) = \pre Fw_x(\tau)$. In other words, the ratio of the derivatives of $\Psi_\vT$ and $\Psi_\vU\circ\psi_{K/F}$ is constant on the interval $a<x<c$. For $a<x<a{+}\delta$, with $\delta$ small and positive, this ratio is equal to $p$: this follows from the relation $w_a(\sigma)/\pre Fw_a(\tau) = p$. Integrating the derivative relation, the result follows. \qed 
\enddemo 
\subhead 
3.8 
\endsubhead 
We prove (3.1.1). Let $\vT \in \ewc F$ be absolutely wild, of degree $p^r$. We first dispose of a singular case. 
\proclaim{Proposition} 
Suppose that $\Psi_\vT$ has a unique jump $a$. The functional equation \rom{(3.1.1)} then holds for $\vT$ and $a = p^r\vs_\vT/(1{+}p^r)$. 
\endproclaim 
\demo{Proof} 
Appealing to 2.4 Proposition part (3), the graph of $\Psi_\vT$, in the range $0\le x\le \vs_\vT$, comprises only segments of the two lines $y=p^{-r}x$, $y=p^rx{-}(p^r{-}1)\vs_\vT$. The latter has slope $p^r$ and passes through the point $(\vs_\vT,\vs_\vT)$. These two lines intersect at the point $(a,p^{-r}a)$, where $a = p^r\vs_\vT/(1{+}p^r)$. Using the Symmetry formulation of 3.1,  the result is clear in this case. \qed 
\enddemo 
We assume henceforth that $\Psi_\vT$ has at least two jumps and proceed by induction on $r$. Suppose $r=1$. In this case, $\Psi_\vT$ has exactly two jumps, and they are related as in 3.6 Proposition. The graph consists of segments of the two lines $y=p^{-1}x$, $y= px{-}(p{-}1)\vs_\vT$ and a non-empty segment of a third line of slope $1$. Using the symmetry formulation, the result is clear in this case.
\par 
Suppose therefore that $r\ge 2$ and that $\Psi_\vT$ has at least two distinct jumps. Let $a = a_\vT$ be the least jump. There exists a character $\chi$ of $F^\times$, of Swan exponent $a$ and order $p$, such that $\chi\vT = \vT$ (as follows from 3.3 Lemma 1). View $\chi$ as a character of $\scr W_F$ and let $\scr W_K = \roman{Ker}\,\chi$. Take $\vU \in \ewc K$ as in 3.7 Proposition. By the inductive hypothesis, 
$$ 
\vs_\vU-y = \Psi_\vU\big(\vs_\vU-\Psi_\vU(y)\big), \quad 0\le y \le \vs_\vU. 
$$ 
Let $z = z_\vT$ be the largest jump of $\Psi_\vT$ and $z_K$ that of $\Psi_\vU\circ \psi_{K/F}$. It follows from 3.5 Lemma 4 that $z_K\le z$. In the range $z<x<\vs_\vT$, we have 
$$ 
\Psi_\vT(x) = \vs_\vT - p^r(\vs_\vT{-}x). 
$$ 
Also, $\vs_\vT{-}x<\vs_\vT{-}z = a/p^r$, by 3.5 Theorem. Therefore 
$$ 
\Psi_\vT(\vs_\vT - \Psi_\vT(x)) = \Psi_\vT(p^r(\vs_\vT{-}x)) = \vs_\vT{-}x, 
$$
as desired. If, on the other hand, $0<x<a$, then $\Psi_\vT(x) = x/p^r$, whence 
$$
\vs_\vT-\Psi_\vT(x) = \vs_\vT-x/p^r > \vs_\vT-a/p^r = z. 
$$ 
Therefore $\Psi_\vT(\vs_\vT{-}\Psi_\vT(x)) = \vs_\vT{-}x$. 
\par 
It remains to treat the range $a<x<z$. Here, $\vs_\vT{-}\Psi_\vT(x) < \vs_\vT{-}\Psi_\vT(a) = \vs_\vT{-}a/p^r = z$. We may therefore apply 3.7 Proposition and (3.5.3) to get 
$$ 
\Psi_\vT(\vs_\vT - \Psi_\vT(x)) = p^{-1}\,\Psi_\vU(\psi_{K/F}(\vs_\vT{-}\Psi_\vT(x))). 
$$ 
We have 
$$ 
\Psi_\vT(x) < \Psi_\vT(z) = \vs_\vT - p^r(\vs_\vT{-}z) = \vs_\vT - a. 
$$
That is,  $\vs_\vT{-}\Psi_\vT(x) > a$. It follows that 
$$ 
\align 
\psi_{K/F}(\vs_\vT{-}\Psi_\vT(x)) &= \psi_{K/F}(\vs_\vT)-p\Psi_\vT(x) \\ 
&= \vs_\vU-p\Psi_\vT(x). 
\endalign  
$$ 
Therefore 
$$ 
\align 
\Psi_\vT(\vs_\vT - \Psi_\vT(x)) &= p^{-1}\,\Psi_\vU(\vs_\vU-p\Psi_\vT(x)) \\ 
&= p^{-1}\,\Psi_\vU(\vs_\vU-\Psi_\vU(\psi_{K/F}(x))) \\ 
&= p^{-1}(\vs_\vU-\psi_{K/F}(x)), 
\endalign 
$$ 
applying the inductive hypothesis at the last step. Finally, 
$$ 
p^{-1}(\vs_\vU-\psi_{K/F}(x)) = p^{-1}(\psi_{K/F}(\vs_\vT)-\psi_{K/F}(x)) = \vs_\vT-x, 
$$ 
and the proof is complete. \qed 
\head\Rm 
4. Symmetry and the bi-Herbrand function 
\endhead 
We turn attention to the GL-side. Let $\vT \in \ewc F$ be of degree $p^r$ (notation of 2.3). In particular, $r\ge 1$. We observed in 3.1 that the functional equation (3.1.1) can be interpreted as a symmetry property of the graph of $\Psi_\vT$. This leads us to define a family of more transparent ``bi-Herbrand functions'' with the same properties of symmetry and convexity. Our objective, realized in section 7, is to calculate $\Psi_\vT$ explicitly as a bi-Herbrand function. However, 4.6 Example at the end of the section does exactly that in a substantial family of cases. 
\subhead 
4.1 
\endsubhead 
We draw out some useful features of the graph $y = \Psi_\vT(x)$. For $\lambda > 0$, let $\eurm i_\lambda$ be the reflection in the line $x{+}y=\lambda$. That is, 
$$ 
\eurm i_\lambda:(x,y) \longmapsto (\lambda{-}y,\lambda{-}x). 
$$ 
\proclaim{Proposition} 
Let $\vT \in \ewc F$ be of degree $p^r$ and abbreviate $\vs = \vs_\vT$. 
\roster 
\item 
The graph $y=\Psi_\vT(x)$, $0\le x\le \vs$, is stable under the reflection $\eurm i_\vs$. 
\item 
There is a unique point $c_\vT$ such that $c_\vT{+}\Psi_\vT(c_\vT) = \vs$. The following conditions are equivalent. 
\itemitem{\rm (a)} The point $c_\vT$ is not a jump of $\Psi_\vT$. 
\itemitem{\rm (b)} The function $\Psi_\vT$ has an even number of jumps. 
\itemitem{\rm (c)} The function $\Psi'_\vT$ takes the value $1$ on a non-empty open subset of the region $0<x<\vs$. 
\itemitem{\rm (d)} The set $I$ of $x$ for which $\Psi'_\vT(x) = 1$ is an open interval containing $c_\vT$. 
\item 
If the conditions \rom{(2)(a)--(d)} hold, then 
$$ 
\Psi_\vT(x) = x-2c_\vT+\vs, \quad x\in I. 
$$ 
\item 
Let $0\le x\le \vs$. In all cases, $\Psi'_\vT(x) \le 1$ if $x{+}\Psi_\vT(x) \le \vs$, while $\Psi'_\vT(x) \ge 1$ if $x{+}\Psi_\vT(x) \ge \vs$. 
\endroster 
\endproclaim 
\demo{Proof} 
Part (1) has been proved in 3.1, as a consequence of (3.1.1). The function $\Psi_\vT$ is strictly increasing, giving the first assertion in (2). The equivalence of (a), (b) and (d) follows from the symmetry of part (1). Suppose (c) holds, and let $I$ be the set of $x$, $0<x<\vs$, for which $\Psi'_\vT(x) = 1$. The convexity of $\Psi_\vT$ implies that $I$ is an interval and symmetry implies $c_\vT \in I$. Thus (c) implies (d), and surely (d) implies (c). 
\par 
In part (3), there is a neighbourhood of $c_\vT$ on which $\Psi_\vT(x) = x{-}b$, for some constant $b$. Thus $\vs = c_\vT{+}\Psi_\vT(c_\vT) = 2c_\vT{-}b$, whence $b = 2c_\vT{-}\vs$, as required. Part (4) follows from the convexity of $\Psi_\vT$ and the symmetry property of (1). \qed 
\enddemo 
\remark{Remark} 
The function $\Psi_\vT$ is continuous and strictly increasing. The condition $x{+}\Psi_\vT(x) \le \vs$ of part (4) is therefore equivalent to $x\le c_\vT$. 
\endremark 
We frequently use the following simple observation, so we exhibit it as a corollary. 
\proclaim{Corollary} 
The function $\Psi_\vT$ has an odd number of jumps if and only if $c_\vT$ is a jump. In that case, $c_\vT$ is the middle one. 
\endproclaim 
\demo{Proof} 
The reflection $\eurm i_\vs$ stabilizes the set of jumps of $\Psi_\vT$ but fixes the point $(c_\vT,\Psi_\vT(c_\vT))$. \qed 
\enddemo  
\subhead 
4.2 
\endsubhead 
We construct a family of $\frak i_\vs$-symmetric functions using more transparent data. They have properties analogous to those in 4.1 Proposition. To specify them, we need two families of auxiliary functions defined using the classical Herbrand functions $\psi_{E/F}$, $\vf_{E/F}$ of section 1. 
\definition{Definition} 
Let $E/F$ be a totally ramified field extension of degree $p^r$, $r\ge1$. Let $\vs = m/p^r$, where $m$ is a positive integer not divisible by 
$p$. Define 
$$ 
\left. \aligned 
\Psi^\times_{(E/F,\vs)}(x) &= p^{-r}\psi_{E/F}(x), \\ 
\Psi^+_{(E/F,\vs)}(x) &= \vs - \vf_{E/F}(p^r(\vs{-}x)), 
\endaligned\,\right\} \qquad 0\le  x  \le \vs.  
\tag 4.2.1 
$$ 
\enddefinition 
The functions $\Psi^\times_{(E/F,\vs)}$, $\Psi^+_{(E/F,\vs)}$ are continuous, strictly increasing, convex and piecewise linear in the region $0\le x\le \vs$. They have only finitely many jumps there. 
\proclaim{Lemma} 
\roster 
\item 
The functions $\Psi_{(E/F,\vs)}^\times$, $\Psi_{(E/F,\vs)}^+$ satisfy 
$$ 
\align 
\vs-x &= \Psi_{(E/F,\vs)}^+\big(\vs-\Psi_{(E/F,\vs)}^\times(x)\big) \\  
&= \Psi_{(E/F,\vs)}^\times\big(\vs-\Psi_{(E/F,\vs)}^+(x)\big). 
\endalign 
$$ 
\item 
There is a unique point $c = c_{(E/F,\vs)}$ such that $c{+}\Psi_{(E/F,\vs)}^\times(c) = \vs$. It further satisfies $c{+}\Psi_{(E/F,\vs)}^
+(c) = \vs$. 
\item 
Let $j_\infty = j_\infty(E|F)$ be the largest jump of $\psi_{E/F}$. If $j_\infty < \vs$ then $j_\infty$ is the largest jump of $\Psi^\times_{(E/F,\vs)}$ and 
$$ 
\bar\jmath_\infty = \vs-\Psi^\times_{(E/F,\vs)}(j_\infty) 
\tag 4.2.2 
$$ 
is the least jump of $\Psi^+_{(E/F,\vs)}$. If $j_\infty < c$, then $c<\bar\jmath_\infty < \vs$. 
\endroster 
\endproclaim 
\demo{Proof} 
Part (1) follows from a simple manipulation of the definition (4.2.1). In (2), the function $\Psi_{(E/F,\vs)}^\times$ is strictly increasing and $
\Psi_{(E/F,\vs)}^\times(0) = 0$, giving the first assertion. For the second, we abbreviate the notation in the obvious way. From (1), $\vs{-}c = \Psi^+(\vs{-}\Psi^\times(c)) = \Psi^+(c)$, as required. The graphs $y = \Psi^\times_{(E/F,\vs)}(x)$, $y = \Psi^+_{(E/F,\vs)}(x)$ are interchanged by the involution $\eurm i_\vs$, whence (3) follows. \qed 
\enddemo 
We define the {\it bi-Herbrand function\/} $\biP EF\vs$ by 
$$ 
\biP EF\vs(x) = \roman{max}\big\{\Psi_{(E/F,\vs)}^\times(x), \Psi_{(E/F,\vs)}^+(x)\big\}, \quad 0\le x\le \vs. 
\tag 4.2.3 
$$ 
When speaking of the jumps of $\biP EF\vs$, we mean the discontinuities of its derivative in the region $0<x<\vs$. 
\proclaim{Proposition} 
Let $j_\infty = j_\infty(E|F)$ and write $c = c_{(E/F,\vs)}$, as in the lemma. 
\roster 
\item 
The function $\biP EF\vs$ is continuous, strictly increasing, piecewise linear and convex, with only finitely many jumps. The graph $y = \biP EF\vs(x)$ is symmetric with respect to the line $x{+}y=\vs$. 
\item 
Suppose $j_\infty \ge c$. The function $\biP EF\vs$ has an odd number of jumps, of which $c$ is the middle one. The derivative $\biP EF\vs'$ does not take the value $1$. Moreover, 
$$ 
\biP EF\vs(x) = \left\{\,\alignedat3 &\Psi_{(E/F,\vs)}^\times(x) > \Psi_{(E/F,\vs)}^+(x), &\quad &0 < x <c, \\ 
&\Psi_{(E/F,\vs)}^+(x) > \Psi_{(E/F,\vs)}^\times(x), &\quad &c<x<\vs. 
\endalignedat \right.  
$$ 
\item 
Suppose $j_\infty < c$. Defining $\bar\jmath_\infty$ as in \rom{(4.2.2),} we have $j_\infty < c <\bar\jmath_\infty$. 
\itemitem{\rm (a)} 
If $j_\infty < x < \bar \jmath_\infty$, then 
$$ 
\gather 
\biP EF\vs'(x) = \Psi_{(E/F,\vs)}^{\times\prime}(x) = \Psi_{(E/F,\vs)}^{+\prime}(x) = 1, \\
\biP EF\vs(x) = \Psi_{(E/F,\vs)}^\times(x) = \Psi_{(E/F,\vs)}^+(x) = x{-}p^{-r}w_{E/F}. 
\endgather 
$$
\itemitem{\rm (b)} 
If $0<x<j_\infty$, then $\Psi_{(E/F,\vs)}^{+\prime}(x) = 1 > \Psi_{(E/F,\vs)}^{\times\prime}(x)$ and 
$$ 
\biP EF\vs(x) = \Psi_{(E/F,\vs)}^\times(x) > \Psi_{(E/F,\vs)}^+(x). 
$$ 
\itemitem{\rm (c)} 
If $\bar \jmath_\infty < x<\vs$, then $\Psi_{(E/F,\vs)}^{\times\prime}(x) = 1 < \Psi_{(E/F,\vs)}^{+\prime}(x)$ and 
$$ 
\biP EF\vs(x) = \Psi_{(E/F,\vs)}^+(x) > \Psi_{(E/F,\vs)}^\times(x). 
$$ 
\item"" In particular, $\biP EF\vs$ has an even number of jumps. 
\endroster 
\endproclaim 
\demo{Proof} 
In (1),  only convexity requires comment, and that is obvious from parts (2) and (3). 
\par 
The index $(E/F,\vs)$ will be constant throughout, so we omit it for the rest of this argument. We have $\Psi^\times(c) = \Psi^+(c) = \upr{2\,}\Psi(c)$. We examine the functions in a small neighbourhood of $x=c$. The values of $\Psi^{\times\prime}(x)$ are of the form $p^{-s}$, and those of $\Psi^{+\prime}(x)$ are $p^s$, for various integers $s$ such that $0\le s\le r$. In part (2), the left derivative of $\Psi^\times$ at $c$ is, at most, $p^{-1}$, while the right derivative of $\Psi^+$ at $c$ is at least $p$. So, $c$ is a jump of $\upr{2\,}\Psi$. The other assertions in (2) follow from the convexity of the functions $\Psi^\times$ and $\Psi^+$. 
\par 
In part (3), the functions $\Psi^\times$, $\Psi^+$ agree, and have derivative $1$, on the interval $j_\infty < x<\bar\jmath_\infty$ (which contains $c$). The derivative relations are clear from the definitions, and imply the main points readily. \qed 
\enddemo 
\remark{Remark} 
By 1.6 Proposition, the condition $j_\infty \ge c$ amounts to 
$$ 
j_\infty{+}\Psi^\times_{(E/F,\vs)}(j_\infty) = 2j_\infty - p^{-r}w_{E/F} \ge \vs. 
$$ 
By 1.6 Corollary, this will hold if $w_{E/F} \ge m(p^r{-}1)/(p^r{+}1)$. 
\endremark 
\subhead 
4.3 
\endsubhead 
We restate 2.6 Proposition in terms of the bi-Herbrand function. 
\proclaim{Proposition} 
Let $\vT\in \ewc F$ be of degree $p^r$. Let $\theta\in \scr C(\frak a,\alpha)$ be a realization of $
\vT$ on a simple stratum $[\frak a,m,0,\alpha]$ in $\M{p^r}F$. If $\vs = \vs_\vT = m/p^r$ and $E = F[\alpha]$ then 
$$ 
\alignedat3
\Psi_\vT(x) = \biP EF\vs (x) &= \Psi^\times_{(E/F,\vs)}(x), &\quad &0\le x\le \vs/2, \\ 
\Psi_\vT(x) = \biP EF\vs (x) &= \Psi^+_{(E/F,\vs)}(x), &\quad &\vs/2 \le \Psi^+_{(E/F,\vs)}(x) \le \vs. 
\endalignedat 
$$ 
\endproclaim 
\demo{Proof} 
The first assertion combines 2.6 Proposition with 4.2 Proposition. The second follows from the symmetry properties of $\Psi_\vT$ and $\biP EF\vs$. \qed 
\enddemo 
\subhead 
4.4 
\endsubhead 
We record the effect of tame lifting on these functions. 
\proclaim{Proposition} 
Let $E/F$ be totally ramified of degree $p^r$ and let $\vs = m/p^r$, for a positive integer $m$ not divisible by $p$. If $K/F$ is a finite tame extension and $e = e(K|F)$, then 
$$ 
\left. \aligned  
\Psi^\times_{(E/F,\vs)}(x) &= \Psi^\times_{(EK/K,e\vs)}(ex)/e, \\ 
\Psi^+_{(E/F,\vs)}(x) &= \Psi^+_{(EK/K,e\vs)}(ex)/e, \\
\upr{2\,}\Psi_{(E/F,\vs)}(x) &= \upr{2\,}\Psi_{(EK/K,e\vs)}(ex)/e,  
\endaligned\ \right\} \qquad 0\le x\le \vs. 
$$ 
\endproclaim 
\demo{Proof} 
This combines the definitions (4.2.1), (4.2.3) with 1.1 Lemma. \qed 
\enddemo 
\subhead 
4.5 
\endsubhead 
The second assertion of 4.3 Proposition determines $\Psi_\vT$ where $\Psi_\vT(x) > \vs/2$. That has already been done in 2.7 Proposition, but in a rather different way. Reconciliation of the two approaches reveals a fundamental property of $\Psi^+_{(E/F,\vs)}$. See 2.5 Definition for the notion of ``twisting datum'' . 
\proclaim{Proposition} 
Let $[\frak a,m,0,\alpha]$ be a simple stratum in $\M{p^r}F$, in which $E = F[\alpha]/F$ is totally ramified of degree $p^r$ and $m$ is not divisible by $p$. Set $\vs = m/p^r$. If $(k,c,\chi)$ is a twisting datum over $F$ such that $k < m/p^r$ is not a jump of $\Psi^+_{(E/F,\vs)}$ then 
$$ 
\Psi^+_{(E/F,\vs)}(k) = t/p^r, 
$$ 
where $t$ is the least integer for which the congruence 
$$ 
u^{-1}\alpha u \equiv \alpha{+}c \pmod{\frak p^{-t}} 
\tag 4.5.1 
$$ 
admits a solution $u\in U^1_\frak a$. 
\endproclaim 
\demo{Proof} 
Assume initially that $2t > m$. For comparison purposes, choose $\theta\in \scr C(\frak a,\alpha)$ and let $\vT$ be the endo-class of $\theta$. Thus $\vT$ is totally wild and of Carayol type. By 4.3 Proposition, $k$ is not a jump of $\Psi_\vT$ and so, by 2.7 Proposition, $t/p^r = \Psi_\vT(k) = \Psi^+_{(E/F,\vs)}(k)$. Because of the jump condition, $t$ depends on $k$ but not on the element $c\in \frak p_F^{-k} \smallsetminus \frak p_F^{1-k}$. 
\par 
We now admit the possibility $2t\le m$. The integer $t$ depends on $\alpha$ and $c$, so we define a function $T(\alpha,c) = p^{-r}t$ where, as before, $t$ is the least integer for which (4.5.1) admits a solution. Let $n$ be a positive integer and take $\nu\in F$ with $\ups_F(\nu) = -n$. 
Thus $[\frak a,m{+}p^rn,0,\nu\alpha]$ is a simple stratum in $
\M{p^r}F$. The congruences 
$$ 
\align 
u^{-1}\alpha u &\equiv \alpha{+}c \pmod{\frak p^{-t}}, \\ 
u^{-1}\nu\alpha u &\equiv \nu(\alpha{+}c) \pmod{\frak p^{-(t+p^rn)}} 
\endalign 
$$ 
have the same sets of solutions $u\in U^1_\frak a$. Consequently, 
$$ 
T(\nu\alpha, \nu c) = T(\alpha,c)+n. 
$$ 
Provided $2T(\nu\alpha, \nu c) > \vs{+}n$, we therefore have 
$$ 
T(\nu\alpha, \nu c) = \Psi^+_{(E/F,\vs{+}n)}(k{+}n).  
$$ 
The definition of the functions $\Psi^+_{(E/F,\vs)}$ implies 
$$ 
\Psi^+_{(E/F,\vs{+}n)}(x{+}n) = \Psi^+_{(E/F,\vs)}(x)+n, 
$$ 
so $k{+}n$ is not a jump of $\Psi^+_{(E/F,\vs{+}n)}(x{+}n)$. The condition $2T(\nu\alpha, \nu c)> \vs{+}n$ thus reduces to $2T(\alpha,c) > \vs-n$. So, for integers $k = -\ups_F(c)$ satisfying $2T(\nu\alpha, \nu c)> k > \vs{+}n$, we have $\Psi^+_{(E/F,\vs)}(k) = T(\alpha,c)$. Allowing $n$ to increase 
without bound, we see that $\Psi^+_{(E/F,\vs)}(k) = T(\alpha,c)$, for all integers $k$ that are not jumps of $\Psi^+_{(E/F,\vs)}$. \qed 
\enddemo 
\remark{Remark} 
The relation between the function $\Psi^+_{(E/F,\vs)}$ and intertwining properties of simple strata was observed in more general work of 
E.-W\. Zink \cite{40}, \cite{41} on a corresponding problem in $F$-division algebras. 
\endremark 
\subhead 
4.6 
\endsubhead 
To finish the section with an example, we calculate $\Psi_\vT$ in a large family of cases. 
\proclaim{Example} 
Let $\vT\in \ewc F$ be of degree $p^r$, $r\ge1$. Let $\theta\in \scr C(\frak a,\alpha)$ be a realization of $\vT$ on a simple stratum $[\frak a,m,0,\alpha]$ in $\M{p^r}F$. Write $\vs = \vs_\vT = m/p^r$ and $E = F[\alpha]$. If $j_\infty(E|F) < \vs/2$, then 
$$ 
\Psi_\vT(x) = \biP EF\vs(x), \quad 0\le x\le \vs. 
$$ 
\endproclaim 
\demo{Proof} 
By 4.3 Proposition, $\Psi_\vT(x) = \upr{2\,}\Psi_{(E/F,\vs)}(x)$ for $0\le x\le \vs/2$. Likewise for $\vs{-}\Psi_\vT(\vs/2) \le x \le \vs$ by symmetry. In particular, $\Psi'_\vT(x) = \upr{2\,}\Psi'_{(E/F,\vs)}(x) = 1$ for $j_\infty<x<\vs/2$. Thus 4.1 Proposition (2) applies. It shows that $\Psi'_\vT(x) = 1$ on the set $j_\infty < x < \vs{-}\Psi_\vT(j_\infty)$. The same argument, using 4.2 Proposition, applies to $\biP EF\vs$, whence $\Psi_\vT(x) = \biP EF\vs(x)$ on this range. Overall,  $\Psi_\vT(x) = \biP EF\vs(x)$ for $0\le x\le \vs$. \qed 
\enddemo 
\proclaim{Gloss} 
The hypothesis $j_\infty < \vs/2$ holds if $w_{E/F} < (p{-}1)m/2p$. 
\endproclaim 
\demo{Proof} 
By 1.6 Corollary, $j_\infty \le p^{1-r}w_{E/F}/(p{-}1)$. \qed 
\enddemo 
\head \Rm 
5. Characters of restricted level 
\endhead 
Let $[\frak a,m,0,\alpha]$ be a simple stratum in $M = \M{p^r}F$, $r\ge 1$, satisfying the usual conditions: 
\roster 
\item 
$E = F[\alpha]/F$ is totally ramified of degree $p^r$, 
\item $m$ is not divisible by $p$ and $\vs = m/p^r$. 
\endroster 
Let $\|\scr C(\frak a,\alpha)\|$ be the set of endo-classes of simple characters $\theta\in \scr C(\frak a,\alpha)$. Thus any $\vT \in \|\scr C(\frak a,\alpha)\|$ lies in $\ewc F$ and has degree $p^r$. In this section, we fix $\alpha$ and identify a set of $\vT \in \|\scr C(\frak a,\alpha)\|$ for which $\Psi_\vT = \biP{F[\alpha]}F\vs$. This will be the set called $\scr L_\alpha$ in the Introduction. In substance, the section is a sequence of increasingly delicate conjugacy calculations. These are progressively interpreted in terms of intertwining properties of simple characters, using the elementary properties of the graphs of the various functions ``$\Psi$'' laid out in section 4.  
\subhead 
5.1 
\endsubhead 
We recall, in the special case to hand, some of the machinery of \cite{15} Chapter 1. Let $\frak p$ be the Jacobson radical of $\frak a$. Define 
$$ 
\align 
A_\alpha:M&\longrightarrow M, \\ x &\longmapsto \alpha x\alpha^{-1}{-}x. 
\endalign 
$$ 
Let $s_{E/F}:M \to E$ be a tame corestriction on $M$, relative to $E/F$. By definition, $s_{E/F}$ is an $(E,E)$-bimodule homomorphism $M \to E$ such that $s_{E/F}(\frak a) = \frak o_E$. For integers $i<j$, we have exact sequences 
$$ 
\gathered 
0\to \frak p^i_E \longrightarrow \frak p^i @>{\ A_\alpha\ }>> \frak p^i @>{\ s_{E/F}\ }>> \frak p^i_E \to 0, \\
0\to \frak p^i_E/\frak p^j_E \longrightarrow \frak p^i/\frak p^j @>{\ A_\alpha\ }>> \frak p^i/\frak p^j @>{\ s_{E/F}\ }>> \frak p^i_E/\frak p^j_E \to 0 . 
\endgathered 
\tag 5.1.1 
$$ 
\indent  
As in 2.1, let $\mu_F$ be a character of $F$ of level one and set $\mu_M = \mu_F\circ \roman{tr}_M$. Let $w_{E/F}$ denote the wild exponent of the field extension $E/F$ (1.6.1). 
\proclaim{Lemma} 
\roster 
\item 
There is a unique character $\mu_E$ of $E$, of level one, such that 
$$ 
\mu_M(x) = \mu_E(s_{E/F}(x)), \quad x\in M. 
\tag 5.1.2 
$$ 
\item 
There is a unique $d \in E$, of valuation $w_{E/F}$, such that $s_{E/F}(y) = yd$, $y\in E$. 
\endroster 
\endproclaim 
\demo{Proof} 
Part (1) is \cite{15} (1.3.7), part (2) follows from \cite{15} (1.3.8). \qed 
\enddemo 
\subhead 
5.2 
\endsubhead 
We introduce a new parameter. 
\definition{Definition} 
Let $\theta \in \scr C(\frak a,\alpha)$. Define $l_E(\theta)$ as the least integer $l\ge0$ for which the character $\theta\Mid U^{l+1}_E$ is trivial. 
\enddefinition 
\proclaim{Proposition} 
Abbreviate $w = w_{E/F}$ and let $\theta\in \scr C(\frak a,\alpha)$.  
\roster 
\item 
If $m > 2w$, then $l_E(\theta) = m{-}w$. 
\item 
If $m\le 2w$, then $0\le l_E(\theta) \le m/2$. If $l$ is an integer, $0\le l\le m/2$, there exists $\vt\in \scr C(\frak a,\alpha)$ such that $l_E(\vt) = l$. 
\endroster 
\endproclaim 
\demo{Proof}
Let $y\in E$, $\ups_E(y) \ge [m/2]{+}1$. The description (2.3.1) of $\theta$ gives 
$$ 
\theta(1{+}y) = \psi_M*\alpha(1{+}y) = \mu_E(\alpha s_{E/F}(y)), 
$$ 
for a tame corestriction $s_{E/F}$ and a character $\mu_E$ of $E$, as in 5.1 Lemma. Also, $\ups_E(s_{E/F}(y)) = \ups_E(y){+}w$. Consequently, if $2w<m$, the character $\theta$ is non-trivial on $U^{1+[m/2]}_E$ and $l_E(\theta) = m{-}w$. Otherwise, $\theta$ is trivial on $U^{1+[m/2]}_E$ and assertion (2) follows from the description in 2.3. \qed 
\enddemo 
\remark{\bf Warning} 
The variation of $l_E(\theta)$ with $E$ is unstable and quite subtle. We explore and exploit this in section 6. 
 \endremark 
\subhead 
5.3 
\endsubhead 
We use the notation $j_\infty$, $\bar\jmath_\infty$ of (4.2.2). We spend the rest of this section proving the following.  
\proclaim{Theorem} 
Let $[\frak a,m,0,\alpha]$ be a simple stratum in $M = \M{p^r}F$, $r\ge1$, in which $E = F[\alpha]/F$ is totally ramified of degree $p^r$ and $p$ does not divide $m$. Set $\vs = m/p^r$ and let $w = w_{E/F}$. Let $\theta\in \scr C(\frak a,\alpha)$ have endo-class $\vT$ and suppose that 
$$ 
l_E(\theta) \le \roman{max}\,\{0,m{-}w\}. 
\tag 5.3.1 
$$ 
\roster 
\item 
If\/ $\upr{2\,}\Psi_{(E/F,\vs)}(x)$ has an odd number of jumps, then $\Psi_\vT(x) = \upr{2\,}\Psi_{(E/F,\vs)}(x)$, $0\le x\le \vs$. 
\item 
If $m>2w$, then $l_E(\theta) = m{-}w$ and $\Psi_\vT(x) = \upr{2\,}\Psi_{(E/F,\vs)} (x)$, $0\le x\le \vs$. 
\item 
If $w$ is divisible by $p$, then $\Psi_\vT(x) = \upr{2\,}\Psi_{(E/F,\vs)}(x)$, $0\le x\le \vs$. 
\item 
Suppose that $m > w\ge m/2$, that $w$ is not divisible by $p$, and that $\upr{2\,}\Psi_{(E/F,\vs)}$ has an even number of jumps. There is a unique character $\phi$ of $U^{m-w}_E$, trivial on $U^{1+m-w}_E$, with the following property. 
\itemitem{\rm (a)} The relation $\Psi_\vT(x) = \upr{2\,}\Psi_{(E/F,\vs)}(x)$ holds for all $x$, $0\le x\le \vs$, if and only if $\theta \Mid U^{m-w}_E \neq \phi$. 
\itemitem{\rm (b)} If $\theta \Mid U^{m-w}_E = \phi$, then 
$$ 
\alignedat3 \Psi_\vT(x) &= \upr{2\,}\Psi_{(E/F,\vs)}(x), &\quad &0\le x \le j_\infty,\ \bar\jmath_\infty\le x\le \vs, \\ 
\Psi_\vT(x) &< \upr{2\,}\Psi_{(E/F,\vs)}(x), &\quad &j_\infty<x<\bar\jmath_\infty.  
\endalignedat 
$$ 
\endroster 
\endproclaim 
\remark{Remarks} 
\roster 
\item 
The hypothesis of part (1) holds if and only if $\biP EF\vs'(x) \neq 1$ for $0<x<\vs$ (4.2 Proposition). It is valid if $w \ge m(p^r{-}1)/(p^r{+}1)$ (4.2 Remark). In particular, if $w\ge m$ then part (1) applies. 
\item 
In part (2), the hypothesis (5.3.1) holds for all $\theta\in \scr C(\frak a,\alpha)$ (5.2 Proposition). This case therefore subsumes 4.6 Example (but we will use this example in the proof of the theorem).
\item 
Regarding (3), the case $w\equiv 0\pmod p$ can only occur when $F$ has characteristic $0$ (1.8). 
\item 
A form of the character $\phi$ in part (4) is written down in (5.12.3) below. A different version is given in 7.3 Remark below, showing that it may or may not be trivial. 
\endroster 
In the theorem, the division into cases (1)--(4) is not exclusive. Certainly (3) can overlap either (1) or (2). When $p=2$, (1) and (2) can overlap (6.2 Example below). Case (4) overlaps no other. 
\endremark 
After preparatory work, part (1) of the theorem is proved in 5.6. Following more preparation in 5.8 and 5.9, parts (2), (3) and (4) are proved in 5.10, 5.11 and 5.12 respectively. 
\subhead 
5.4 
\endsubhead 
Let $\frak p$ be the Jacobson radical of $\frak a$. Let $c\in F$, with $\ups_F(c) = -k$ and $k<\vs = m/p^r$. Let $t < p^rk$ be an integer. As a first step, we consider formal intertwining between the simple strata $[\frak a,m,t,\alpha]$ and $[\frak a,m,t,\alpha{+}c]$. That is, we analyze the congruence 
$$ 
u^{-1}\alpha u \equiv \alpha{+}c \pmod{\frak p^{-t}}, \quad u \in U^1_\frak a. 
\tag 5.4.1 
$$ 
\proclaim{Lemma} 
The set of solutions $u\in U^1_\frak a$ of \rom{(5.4.1)} is either empty or constitutes one coset $u U^1_EU^{m-t}_\frak a\in U^1_EU^{m-p^rk}_\frak a/U^1_EU^{m-t}_\frak a$. 
\endproclaim 
\demo{Proof} 
Let $u\in U^1_\frak a$ satisfy (5.4.1). Thus $u$ conjugates the equivalence class of the simple stratum $[\frak a,m,t,\alpha]$ to that of $[\frak a,m,t,\alpha{+}c]$. If $v\in U^1_\frak a$ and $uv$ satisfies (5.4.1), then $v$ conjugates the equivalence class of the stratum $[\frak a,m,t,\alpha{+}c]$ to itself. Equivalently, $v \in U^1_EU^{m-t}_\frak a$ \cite{15} (1.5.8), so the coset $uU^1_EU^{m-t}_\frak a$ is uniquely determined by (5.4.1). On the other hand, $u$ conjugates the equivalence class of $[\frak a,m,p^rk,\alpha]$ to itself, so $u\in U^1_EU^{m-p^rk}_\frak a$ {\it loc\. cit\.} \qed 
\enddemo 
\remark{Remark} 
Since $U^1_E$ commutes with $\alpha$, we need only ever consider solutions $u$ of (5.4.1) that satisfy $u \in U^{m-p^rk}_\frak a$. 
\endremark 
\subhead 
5.5 
\endsubhead 
We continue with the same notation. In (5.4.1), write $u=1{+}a$, $a\in \frak p^{m-p^rk}$. In this form, (5.4.1) amounts to 
$$ 
(1{+}a)^{-1}\alpha(1{+}a) \equiv \alpha{+}c \pmod{\frak p^{-t}} 
\tag 5.5.1 
$$ 
or, equivalently, 
$$ 
\alpha a{-}a\alpha \equiv c(1{+}a) \pmod{\frak p^{-t}}. 
\tag 5.5.2 
$$ 
We use the standard notation $[x,y] = xy{-}yx$, for $x,y\in M$. 
\proclaim{Proposition} 
Let $a\in \frak p^{m-p^rk}$ satisfy \rom{(5.5.1).} If $y\in E$, $\ups_E(y) = b\ge 1$, then 
$$ 
(1{+}a)(1{+}y)(1{+}a)^{-1} \equiv 1{+}\bar y \pmod{\frak p^{b+m-t}}, 
$$ 
for an element $\bar y\in E$ such that $\bar y \equiv y \pmod{\frak p_E^{b+m-p^rk}}$. 
\endproclaim  
\demo{Proof} 
We re-arrange the conjugation as 
$$
(1{+}a)(1{+}y)(1{+}a)^{-1} = 1+y+[a,y](1{+}a)^{-1}. 
$$ 
Applying the defining relations (5.5.1), (5.5.2), we get 
$$ \multline 
\left[\alpha, [a,y](1{+}a)^{-1}\right] \\ 
\alignedat3 
&= \alpha[a,y](1{+}a)^{-1} - [a,y](1{+}a)^{-1}\alpha & & \\ 
&\equiv \alpha[a,y](1{+}a)^{-1} - [a,y](\alpha{+}c)(1{+}a)^{-1} &\quad &\pmod{\frak p^{b+m-p^rk-t}} \\ 
&\equiv \big([\alpha,a]y - y[\alpha,a] - [a,y]c\big)(1{+}a)^{-1} &\quad &\pmod{\frak p^{b+m-p^rk-t}} \\ 
&\equiv \big(c(1{+}a)y-yc(1{+}a)-[a,y]c\big)(1{+}a)^{-1} &\quad &\pmod{\frak p^{b-t}} \\ 
&\equiv 0 &\quad &\pmod{\frak p^{b-t}} .  
\endalignedat \endmultline 
$$
The exact sequences (5.1.1) imply $[a,y](1{+}a)^{-1} = v{+}h$, for $v\in \frak p_E^{b+m-p^rk}$ and $h\in \frak p^{b+m-t}$, as required. \qed 
\enddemo 
\subhead 
5.6 
\endsubhead 
We continue with the same notation, especially $\vs = m/p^r$ and $w = w_{E/F}$. 
\proclaim{Proposition 1} 
Let $I$ be an open sub-interval of $(0,\vs)$  on which $\Psi^\times_{(E/F,\vs)}$ and $\Psi^+_{(E/F,\vs)}$ are both smooth and satisfy 
$$ 
\Psi^\times_{(E/F,\vs)}(x) > \Psi^+_{(E/F,\vs)}(x), \quad x\in I. 
\tag 5.6.1 
$$ 
Let $\theta\in \scr C(\frak a,\alpha)$, and suppose   
$$ 
l = l_E(\theta) \le \roman{max}\,\{0,m{-}w\}. 
\tag 5.6.2 
$$ 
If $\theta$ has endo-class $\vT$, then 
$$ 
\Psi_\vT(x) = \Psi^\times_{(E/F,\vs)}(x) = \upr{2\,}\Psi_{(E/F,\vs)}(x), \quad x\in I. 
$$ 
\endproclaim 
\demo{Proof} 
By (4.2.1), (4.2.3), we have $\Psi^{\times\prime}_{(E/F,\vs)}(x) \le 1\le \Psi^{+\prime}_{(E/F,\vs)}(x)$, $0 < x <\vs$. By 4.2 Proposition,
the hypothesis (5.6.1) implies that $\Psi^{\times\prime}_{(E/F,\vs)}(x) < 1$, $x\in I$. The convexity of $\psi_{E/F} = p^r\Psi^\times_{(E/F,\vs)}$ and 1.6 Proposition now imply $\Psi^\times_{(E/F,\vs)}(x) > x{-}p^{-r}w$, $x\in I$. 
\par 
As in the proof of 2.6 Proposition, the tame lifting properties of 4.4 Proposition and 2.5 Proposition 2 show it is enough to prove the result when $x$ is an integer. So, let $k$ be an integer, $k\in I$, that is not a jump of $\Psi_\vT$. Let $(k,c,\chi)$ be a twisting datum (2.5). We apply 5.5 Proposition with $t = p^r\Psi^+_{(E/F,\vs)}(k)$. By 4.5 Proposition, $t$ is the least integer for which the congruence (5.5.1) admits a solution $a$. By 5.4 Lemma, we may take $a\in \frak p^{m-p^rk}$. The definition of $I$ implies that $k$ is not a jump of $\psi_{E/F}$, so $\psi_{E/F}(k) = \sw(\chi\circ\N EF)$ is an integer (1.3 Proposition). 
\par 
Write $v = \psi_{E/F}(k)$ and let $y\in E$ have valuation $1{+}v$. In particular, $\chi\circ\det(1{+}y) = 1$ ({\it cf\.} 1.3 Proposition). Our hypothesis (5.6.1) amounts to 
$$ 
\psi_{E/F}(k) = p^r\Psi^\times_{(E/F,\vs)}(k)  > p^r\Psi^+_{(E/F,\vs)}(k) = t, 
$$ 
so $v>t$. Thus 5.5 Proposition gives 
$$ 
(1{+}a)(1{+}y)(1{+}a)^{-1} \equiv 1{+}\bar y \pmod{\frak p^{2+m}}, 
$$ 
whence $1{+}a$ normalizes the group $H^{1+v}(\alpha,\frak a)$ and $\theta^{1+a}(1{+}y) = \theta(1{+}\bar y)$. Taking first the case $l = 0$, we get $\theta^{1+a}(1{+}y) = \theta(1{+}\bar y) = 1 = \chi\theta(1{+}y)$. In the other case $0 < l \le m{-}w$, 
$$ 
\ups_E(\bar y{-}y) \ge 1{+}v{+}m{-}p^rk = 1+ \psi_{E/F}(k) +m -p^rk \ge 1+ m-w, 
$$
since $\psi_{E/F}(k) \ge p^rk{-}w$. It follows that $\theta^{1+a}(1{+}y) = \theta(1{+}\bar y) = \theta(1{+}y) = \chi\theta(1{+}y)$. By hypothesis (5.6.1), $t< v$ so the definition of $a$ ensures that $1{+}a$ conjugates $\theta$ to $\chi\theta$ on $H^{1+v}(\alpha,\frak a)$. Therefore $\Psi_\vT(k) \le v/p^r = \biP EF\vs(k)$. 
\par 
We go through the same process with $\ups_E(y) = v = \psi_{E/F}(k)$. We choose $y$ so that $\chi\circ\det(1{+}y) = \chi\circ\N EF(1{+}y) \neq 1$. If $m>w$, then 
$$ 
\ups_E(\bar y{-}y) \ge v{+}m{-}p^rk > m{-}w \ge l, 
$$ 
whence $\theta^{1+a}(1{+}y) = \theta(1{+}y) \neq \chi\theta(1{+}y)$. The element $1{+}a$ therefore normalizes $H^v(\alpha,\frak a)$ but does not conjugate $\theta$ to $\chi\theta$ on that group. If $m\le w$ then $l = 0$ and the same conclusion holds. 
\par 
Suppose there exists $1{+}b\in U^1_\frak a$ that intertwines $\theta$ with $\chi\theta$ on $H^v(\alpha,\frak a)$: that is, it conjugates $\theta$ to $\chi\theta$ on that group. It therefore conjugates $\theta$ to $\chi\theta$ on $H^{1+v}(\alpha,\frak a)$ and so is of the form $1{+}b = u(1{+}a)$, where $u\in U^1_\frak a$ conjugates $\theta\Mid H^{1+v}(\alpha,\frak a)$ to itself. 
\proclaim{Lemma} 
We have $v\le [m/2]$. 
\endproclaim 
\demo{Proof} 
The hypothesis (5.6.1) implies that $k$ is strictly less that the largest jump of $\psi_{E/F}$. Therefore $v = \psi_{E/F}(k)\le p^{r-1}k$. On the other hand, $k< \vs = m/p^r$. Suppose that $v>[m/2]$. Since $v$ is an integer, this implies $v>m/2$ and so 
$$ 
m/2 < v \le p^{r-1}k < m/p, 
$$ 
which is ridiculous. \qed 
\enddemo 
Following the lemma, the element $u$ conjugates $\theta$ to itself on $H^1(\alpha,\frak a)$, as follows from \cite{15} (3.3.2). Therefore 
$$ 
\theta^{1+a}\Mid H^v(\alpha,\frak a) = \theta^{1+b}\Mid H^v(\alpha,\frak a) = \chi\theta\Mid H^v(\alpha,\frak a), 
$$ 
which is false. We conclude that $\theta$ does not intertwine with $\chi\theta$ on $H^v(\alpha,\frak a)$, and so $\Psi_\vT(k) = \Bbb A(\vT,\chi\vT) = v/p^r = \biP EF\vs(k)$, as required. \qed 
\enddemo 
Proposition 1 has a ``mirror image'' as follows. 
\proclaim{Proposition 2} 
Let $I$ be an open sub-interval of $(0,\vs)$  on which $\Psi^\times_{(E/F,\vs)}$, $\Psi^+_{(E/F,\vs)}$ are smooth and satisfy 
$$ 
\Psi^\times_{(E/F,\vs)}(x) < \Psi^+_{(E/F,\vs)}(x), \quad x\in I. 
\tag 5.6.3 
$$ 
Let $\theta\in \scr C(\frak a,\alpha)$, and suppose   
$$ 
l = l_E(\theta) \le \roman{max}\,\{0,m{-}w\}. 
\tag 5.6.4 
$$ 
If $\theta$ has endo-class $\vT$, then 
$$ 
\Psi_\vT(x) = \Psi^+_{(E/F,\vs)}(x) = \upr{2\,}\Psi_{(E/F,\vs)}(x), \quad x\in I. 
$$ 
\endproclaim 
\demo{Proof} 
The symmetry property of $\Psi_\vT$ (3.1.1) and the corresponding properties (4.2 Lemma) connecting $\Psi^\times$ with $\Psi^+$ together show that this proposition is equivalent to Proposition 1. \qed 
\enddemo 
\demo{Proof of\/ {\rm 5.3 Theorem (1)}}
Here, $\biP EF\vs$ has an odd number of jumps. The interval $0<x<\vs$, with the jumps of $\biP EF\vs$ removed, is covered by a finite union of open intervals $I_j$ on which either (5.6.1) or (5.6.3) holds. The propositions imply that $\Psi_\vT(x) = \biP EF\vs(x)$ for $x\in \bigcup_j I_j$. By continuity, the functions are equal for $0\le x\le \vs$. \qed 
\enddemo 
The argument used to prove part (1) of 5.3 Theorem has broader applicability. As before, $\vT$ is the endo-class of a simple character $\theta\in \scr C(\frak a,\alpha)$ satisfying (5.3.1). 
\proclaim{Corollary 1} 
If\/ $\biP EF\vs$ has an even number of jumps, then 
$$ 
\Psi_\vT(x) = \biP EF\vs(x) 
$$ 
for all $x$ such that $0\le x\le j_\infty$ or $\bar\jmath_\infty \le x\le \vs$. 
\endproclaim 
\demo{Proof} 
In the region $0<x<j_\infty$, we have $\Psi^\times_{(E/F,\vs)}(x) > \Psi^+_{(E/F,\vs)}(x)$ by 4.2 Proposition. Proposition 1 then implies $\Psi_\vT(x) = \Psi^\times_{(E/F,\vs)}(x) = \biP EF\vs(x)$ for $0\le x\le j_\infty$. Proposition 2 implies $\Psi_\vT(x) = \Psi^+_{(E/F,\vs)}(x) = \biP EF\vs(x)$ for $\bar\jmath_\infty\le x\le \vs$. \qed 
\enddemo 
We can push this train of thought a little further. 
\proclaim{Corollary 2} 
If\/ $\biP EF\vs$ has an even number of jumps, then 
$$ 
\Psi_\vT(x) \le \biP EF\vs(x),\quad j_\infty< x < \bar\jmath_\infty.
$$ 
The following conditions are equivalent. 
\roster 
\item $\Psi_\vT(x_0) = \biP EF\vs(x_0)$, for some $x_0$ such that  $j_\infty < x_0 < \bar\jmath_\infty$.  
\item $\Psi_\vT(x) = \biP EF\vs(x)$ for all $x$ such that $j_\infty < x < \bar\jmath_\infty$. 
\endroster 
\endproclaim 
\demo{Proof} 
For $j_\infty \le x\le \bar\jmath_\infty$, we have $\biP EF\vs(x) = x{-}p^{-r}w$. The functions $\Psi_\vT$, $\biP EF\vs$ agree at the end-points $j_\infty$, $\bar\jmath_\infty$. As $\Psi_\vT$ is convex in this region, so (1) implies (2). The converse is trivial. \qed 
\enddemo 
Corollary 2 provides the basis of a strategy for proving the remaining assertions of 5.3 Theorem. 
\subhead 
5.7 
\endsubhead 
Before we can develop this strategy, we need a minor result derived from elementary linear algebra. 
\par 
Let $\bk k$ be a field and $V$ a $\bk k$-vector space of finite dimension $n$. Let $\frak n$ be a {\it regular\/} nilpotent endomorphism  of $V$. The $\frak n$-stable subspaces of $V$ are then $V_j = \frak n^j(V)$, $0\le j\le n$. 
\proclaim{Lemma 1} 
Let $\frak n'$ be a nilpotent endomorphism of $V$ that commutes with $\frak n$. There exists $a = a(V,\frak n,\frak n')\in \bk k$ such that 
$$ 
\frak n'(v) \equiv a\frak n(v) \pmod{V_{j+2}}, \quad v\in V_j, 
$$ 
for $0\le j\le n{-}2$. The element $a$ is non-zero if and only if $\frak n'$ is regular. 
\endproclaim 
\demo{Proof} 
Let $\frak m \in \End{\bk k}V$ commute with $\frak n$. There is a unique polynomial $\phi(X) \in \bk k[X]$, of degree at most $n{-}1$, such that $\frak m = \phi(\frak n)$. The endomorphism $\frak m$ is nilpotent if and only if $\phi(0) = 0$. If this holds, the linear coefficient $a = \phi'(0)$ satisfies $\frak m(v) \equiv a\frak n(v) \pmod{V_{j+2}}$, $v\in V_j$, as required. \qed 
\enddemo 
We apply Lemma 1 in the following context. Let $[\frak a,m,0,\alpha]$ be a simple stratum in $M = \M{p^r}F$, $E = F[\alpha]$, as in the theorem. Let $\Bbbk_F = \frak o_F/\frak p_F$ be the residue field of $F$. If $\frak p$ is the Jacobson radical of $\frak a$, the $\Bbbk_F$-algebra $\frak a/\frak p$ is isomorphic to $\Bbbk_F^{p^r}$ and $\alpha$ acts on it by conjugation. 
\proclaim{Lemma 2} 
The endomorphism of $\frak a/\frak p$, induced by $A_\alpha$, is regular nilpotent. 
\endproclaim 
\demo{Proof} 
As an endomorphism of the $\Bbbk_F$-space $\frak a/\frak p$, $A_\alpha = \roman{Ad}\,\alpha-1$ satisfies $(A_\alpha)^p = A_{\alpha^p}$, and so $(A_\alpha)^{p^r} = A_{\alpha^{p^r}}$. However, $\alpha^{p^r} \in F^\times U_\frak a$, whence $\roman{Ad}\,\alpha^{p^r}$ induces the identity map on $\frak a/\frak p$. That is, $(A_\alpha)^{p^r} = 0$ and so $A_\alpha$ is nilpotent. By (5.1.1), $\roman{Ker}\,A_\alpha$ is the $1$-dimensional subspace $\frak o_E/\frak p_E$ of $\frak a/\frak p$, so $A_\alpha$ is regular. 
\qed 
\enddemo 
\proclaim{Proposition} 
Let $V_j = A_\alpha^j(\frak a/\frak p)$. Let $s$ be an integer and write $\zeta = s/m \in \frak o_F$. If $\beta \in E$ has valuation $\ups_E(\beta) = -s$, then 
$$ 
A_\beta(v) \equiv \zeta A_\alpha(v) \pmod{V_{j+2}}, 
$$ 
for $v\in V_j$, $0\le j\le p^r{-}2$. 
\endproclaim 
\demo{Proof} 
The set of indecomposable idempotents of the $\Bbbk_F$-algebra $\frak a/\frak p$ provides a $\Bbbk_F$-basis that is permuted cyclically by $\roman{Ad}\,\alpha$, with period $p^r$. We have $A_\alpha = \roman{Ad}\,\alpha - 1$. Similarly for $A_\beta$, and $A_\beta{+}1 = (A_\alpha{+}1)^t$, for an integer $t$, $0\le t\le p^s{-}1$, such that $s\equiv mt \pmod{p^s}$. The linear term in $(A_\alpha{+}1)^t$ is $t A_\alpha$, whence the result follows. \qed 
\enddemo 
\subhead 
5.8 
\endsubhead 
We return to the proof of 5.3 Theorem, as it was left at the end of 5.6. We may now assume that $\biP EF\vs$ has an even number of  jumps. Let $I$ be the non-empty open interval $j_\infty < x< \bar\jmath_\infty$. So, for $x\in I$, 
$$ 
\biP EF\vs(x) = \Psi^\times_{(E/F,\vs)}(x) = \Psi^+_{(E/F,\vs)}(x) = x{-}p^{-r}w, 
$$ 
where $w = w_{E/F}$. Let $(k,c,\chi)$ be a twisting datum with $k\in I$: in particular, $w< p^rk$. Our aim, in this sub-section and the next, is to refine 5.5 Proposition in this more restricted context. 
\par 
By 4.5 Proposition, the congruence
$$ 
(1{+}a)^{-1}\alpha(1{+}a) \equiv \alpha{+}c \pmod{\frak p^{-t}} 
\tag 5.8.1 
$$ 
admits a solution $a$ if and only if $t \ge p^r\Psi^+_{(E/F,\vs)}(k) =  p^rk{-}w$. We examine these solutions $a$ more closely when $t = p^rk{-}w$. As in 5.4 Remark, we need only consider elements $a \in \frak p^{m-p^rk}$. 
\par 
Re-write (5.8.1) in the form 
$$ 
A_\alpha(a) \equiv (1{+}a)c\alpha^{-1} \pmod{\frak p^{m-p^rk+w}}, 
\tag 5.8.2 
$$ 
and set 
$$ 
\eps = A_\alpha(a) - (1{+}a)c\alpha^{-1} \in \frak p^{m-p^rk{+}w}. 
\tag 5.8.3 
$$ 
By 4.5 Proposition, the congruence 
$$ 
A_\alpha(a') \equiv (1{+}a')c\alpha^{-1} \pmod{\frak p^{1+m-p^rk+w}} 
\tag 5.8.4   
$$ 
has no solution $a'$. 
\proclaim{Lemma} 
The element $\eps$ of \rom{(5.8.3)} satisfies $\ups_E(s_{E/F}(\eps)) = m-p^rk+w$ and so $\ups_E(s_{E/F}(a) ) \ge w$. 
\endproclaim 
\demo{Proof} 
Write $t = p^rk{-}w$. Suppose, for a contradiction, that $\ups_E(s_{E/F}(\eps)) > m{-}t$. Take $a\in \frak p^{m-p^rk}$ satisfying (5.8.2): the element $a$ is then determined modulo $\frak p_E^{m-p^rk}{+}\frak p^{m-t}$ (5.4 Lemma). Let $y\in \frak p^{m-t}$ and consider the congruence 
$$ 
A_\alpha(a{+}y) \equiv (1{+}a{+}y)c\alpha^{-1} \pmod{\frak p^{1+m-t}}.  
$$ 
Since $m>p^rk$, we can neglect the term $yc\alpha^{-1}$, so this congruence amounts to 
$$ 
A_\alpha(a{+}y) \equiv (1{+}a)c\alpha^{-1}\pmod{\frak p^{1+m-t}},  
$$ 
that is, 
$$ 
A_\alpha(y) \equiv -\eps \pmod{\frak p^{1+m-t}}.  
$$ 
We have assumed that $\ups_E(s_{E/F}(\eps)) > m{-}t$ so, by (5.1.1), this last congruence admits a solution $y\in \frak p^{m-t}$. The element $a' = a{+}y$ then satisfies (5.8.4), which is impossible. This proves the first assertion. Now apply $s_{E/F}$ to the definition (5.8.3). Since $s_{E/F}(1)$ has valuation $w$, the second assertion follows directly. \qed 
\enddemo 
\subhead 
5.9 
\endsubhead 
We continue in the situation of 5.8. In particular, $(k,c,\chi)$ is a twisting datum such that $j_\infty < k < \bar\jmath_\infty$. In particular, $p^rk > w$, by 4.2 Proposition (3)(a). Going forward, we impose the following simplification: 
\proclaim{Assumption} 
We henceforward assume that $\biP EF\vs(\bar\jmath_\infty) \le \vs/2$. 
\endproclaim 
\demo{Justification} 
For, if $\biP EF\vs(\bar\jmath_\infty) > \vs/2$, the functional equation implies $j_\infty < \vs/2$ and we are in the situation of 4.6 Example. In that case, we know that $\Psi_\vT(x) = \biP EF\vs(x)$ for $0\le x\le \vs$, as demanded by part (2) of the theorem. \qed 
\enddemo 
\proclaim{Proposition 1} 
Let $\theta\in \scr C(\frak a,\alpha)$ satisfy $l_E(\theta) \le m{-}w$. Let $a\in \frak p^{m-p^rk}$ be a solution of 
\rom{(5.8.2).} Define $\eps$ by \rom{(5.8.3)} and set $\zeta = w/m \in \frak o_F$. If $y\in E$ and $\ups_E(y) \ge p^rk{-}w$, then 
$$ 
\theta^{1+a}(1{+}y)/\theta(1{+}y) = \theta(1{-}\zeta c\alpha^{-1} y)\,\mu_M(-\alpha\zeta\eps y). 
\tag 5.9.1 
$$ 
\endproclaim 
\demo{Proof} 
Suppose first that $\ups_E(y) > p^rk{-}w$. This implies $\theta(1{-}\zeta c\alpha^{-1} y) = 1$, while $\mu_M(-\alpha\zeta\eps y) = 1$ by 5.8 Lemma. The right hand side of (5.9.1) thus equals $1$. Application of 5.5 Proposition gives the same for the left hand side. Assume now that $\ups_E(y) = p^rk{-}w$. 
\proclaim{Lemma 1} 
If $y\in E$ and $\ups_E(y) = p^rk{-}w$, then 
$$ 
(1{+}a)(1{+}y)(1{+}a)^{-1} = 1{+}\bar y+h, 
$$ 
for elements $\bar y$ of $E$ and $h$ of $\frak p^m$ such that 
$$ 
\align  \bar y &\equiv y \pmod{\frak p_E^{m-w}}, \\ 
h&\equiv -\zeta\eps y \pmod{A_\alpha(\frak p^m)+\frak p^{m+1}}. 
\endalign 
$$ 
\endproclaim 
\demo{Proof} 
We re-write the defining relation (5.8.1), with $t=p^rk{-}w$, as 
$$ 
(1{+}a)^{-1}\alpha(1{+}a) = \alpha{+}c{+}\delta. 
\tag 5.9.2 
$$ 
Thus $\delta \in \frak p^{w-p^rk}$ and 
$$ 
\align 
[\alpha,a] &= (1{+}a)(c{+}\delta), \\ 
A_\alpha(a) &= (1{+}a)(c{+}\delta)\alpha^{-1}. 
\endalign 
$$ 
Therefore $\eps = (1{+}a)\delta\alpha^{-1}$. We start from  the identity 
$$ 
(1{+}a)(1{+}y)(1{+}a)^{-1} = 1+y+[a,y](1{+}a)^{-1} 
\tag 5.9.3 
$$ 
and evaluate, using (5.9.2). We find 
$$ 
\align 
\big[\alpha,[a,y](1{+}a)^{-1}\big] &= \alpha[a,y](1{+}a)^{-1}-[a,y](1{+}a)^{-1}\alpha \\ 
&= \alpha[a,y](1{+}a)^{-1} - [a,y](\alpha{+}c{+}\delta)(1{+}a)^{-1} \\ 
&= \big([\alpha,a]y-y[\alpha,a]-[a,y](c{+}\delta)\big)(1{+}a)^{-1} \\ 
&= \big((1{+}a)(c{+}\delta)y-y(1{+}a)(c{+}\delta) - [a,y](c{+}\delta)\big)(1{+}a)^{-1} \\ 
&= \big((1{+}a)\delta y-y(1{+}a)\delta -[a,y]\delta\big)(1{+}a)^{-1} \\ 
&= (1{+}a)[\delta,y](1{+}a)^{-1}. 
\endalign 
$$ 
Substituting for $\delta$, we get 
$$ 
\align 
(1{+}a)[\delta,y](1{+}a)^{-1} &= (1{+}a)\big[(1{+}a)^{-1}\eps\alpha,y\big](1{+}a)^{-1} \\ 
&=(1{+}a)\big((1{+}a)^{-1}\eps\alpha y - y(1{+}a)^{-1}\eps\alpha\big)(1{+}a)^{-1} \\ 
&\equiv [\eps\alpha,y]  \pmod{\frak p}, 
\endalign 
$$ 
since $[\eps\alpha,y] \in \frak a$ and $a\in \frak p$. Thus 
$$ 
\align 
A_\alpha\big((1{+}a)(1{+}y)(1{+}a)^{-1}\big) &= A_\alpha([a,y](1{+}a)^{-1}) \\ 
&\equiv [\eps\alpha,y]\alpha^{-1} \equiv [\eps,y] \pmod{\frak p^{m+1}}. 
\endalign 
$$ 
We have $[a,y](1{+}a)^{-1} \in \frak p^{m-w}$ and $[\eps,y]\in \frak p^m$. It follows that 
$$ 
(1{+}a)(1{+}y)(1{+}a)^{-1} = 1+\bar y+h, 
$$ 
where $\bar y\in \frak p_E^{p^rk-w}$ satisfies $\bar y\equiv y \pmod{\frak p_E^{m-w}}$ and $h\in \frak p^m$ satisfies 
$$ 
A_\alpha(h) \equiv [\eps,y] \pmod{\frak p^{1+m}}. 
$$ 
By 5.7 Proposition, 
$$ 
[\eps,y] = -A_y(\eps)y \equiv -\zeta A_\alpha(\eps)y \pmod{A_\alpha^2(\frak p^m)+\frak p^{m+1}}. 
$$ 
Adjusting $\bar y$ by an element of $\frak p_E^m$, which changes nothing, we may choose $h$ to satisfy 
$$ 
h\equiv -\zeta\eps y \pmod{A_\alpha(\frak p^m)+\frak p^{1+m}},   
$$ 
as required. \qed 
\enddemo 
The elementary identity (5.9.3) implies 
$$ 
(1{+}a)(1{+}y)(1{+}a)^{-1} \equiv 1+y+[a,y] \pmod{\frak p^{1+m-w}}. 
\tag 5.9.4 
$$ 
\proclaim{Lemma 2} 
Let $\ups_E(y) = p^rk{-}w$. If $\zeta = w/m\in \frak o_F$, then $[a,y] \equiv -\zeta A_\alpha(a)y \pmod{\frak p^{1+m-w}}$. 
\endproclaim 
\demo{Proof} 
The defining relation $A_\alpha(a)\equiv (1{+}a)c\alpha^{-1} \pmod{\frak p^{m-p^rk+w}}$ implies that $A_\alpha^2(a) \in \frak p^{1+m-p^rk}$. That is, 
$$ 
A_\alpha(a) \in \frak p_E^{m-p^rk}{+}\frak p^{1+m-p^rk} = A_\alpha^{p^r-1}(\frak p^{m-p^rk}){+}\frak p^{1+m-p^rk}. 
$$ 
Therefore $a\in  A_\alpha^{p^r-2}(\frak p^{m-p^rk}){+}\frak p^{1+m-p^rk}$. We apply 5.7 Proposition to get 
$$ 
[a,y] = -A_y(a)y \equiv -\zeta A_\alpha(a)y \pmod{A_\alpha^{p^r}(\frak p^{m-w})+\frak p^{1+m-w}}. 
$$ 
Since $A_\alpha^{p^r}(\frak p^{m-w})\i\frak p^{1+m-w}$, we have the result. \qed 
\enddemo 
Lemmas 1 and 2 imply $[a,y] \equiv -\zeta A_\alpha(a)y \equiv -\zeta c\alpha^{-1} y \pmod{\frak p^{1+m-w}}$, and the proposition follows from (5.9.4). \qed 
\enddemo 
\remark{Remark} 
Consider the right hand side of the equation (5.9.1). The dependence on $a$ enters only via the element $\eps$, and the expression depends only on $s_{E/F}(\eps)$ modulo $\frak p_E^{1+m-p^rk+w}$. The element $a\in \frak p^{m-p^rk}$ is only determined, as a solution of (5.8.2), modulo $\frak p_E^{m-p^rk}{+}\frak p^{m-p^rk+w}$ (5.4 Lemma). The definition (5.8.3) of $\eps$ implies that $s_{E/F}(\eps)+\frak p_E^{1+m-p^rk+w}$, {\it does not depend on the choice of the solution $a$.} It follows that (5.9.1) holds equally for all solutions $a$ of (5.8.2). 
\endremark 
\proclaim{Corollary} 
In the notation of \rom{Proposition 1,} the following conditions are equivalent. 
\roster 
\item $\Psi_\vT(k) < k{-}p^{-r}w = \biP EF\vs(k)$. 
\item 
$\theta(1{-}\zeta c\alpha^{-1}y)\,\mu_M(-\alpha\zeta\eps y) = \mu_M(cy)$, for all $y\in \frak p_E^{p^rk-w}$. 
\endroster 
\endproclaim 
\demo{Proof} 
If $y\in \frak p_E^{1+p^rk-w}$, the proposition gives 
$$ 
\theta^{1+a}(1{+}y)/\theta(1{+}y) = 1 = \mu_M(cy). 
$$ 
Our Assumption implies $\theta^{1+a}(1{+}x) = \chi\theta(1{+}x)$, for $x\in \frak p^{1+[m/2]}$. Therefore $1{+}a$ conjugates $\theta$ to $\chi\theta$ on $H^{1+p^rk-w}(\alpha,\frak a)$. For the same reason, if (2) holds, then $1{+}a$ conjugates $\theta$ to $\chi\theta$ on $H^{p^rk-w}(\alpha,\frak a)$, which implies (1). 
\par 
Conversely, suppose that (1) holds: there exists $1{+}b\in U^1_\frak a$ that conjugates $\theta$ to $\chi\theta$ on $H^{p^rk-w}(\alpha,\frak a)$. Thus $(1{+}b) = u(1{+}a)$, for some $u\in U^1_\frak a$ that conjugates $\theta\Mid H^{1+p^rk-w}(\alpha,\frak a)$ to itself. By the Assumption again, any such $u$ conjugates $\theta$ to itself, whence $\theta^{1+a}\Mid H^{p^rk-w}(\alpha,\frak a) = \chi\theta\Mid H^{p^rk-w}(\alpha,\frak a)$ and this implies (2). \qed 
\enddemo 
We shall apply Proposition 1 in combination with the following result. 
\proclaim{Proposition 2} 
Let $k \in I$ be an integer and suppose that $\Psi_\vT$ is smooth at $k$. The following conditions are equivalent. 
\roster 
\item There is a twisting datum $(k,c,\chi)$ relative to which 
$$ 
\theta^{1+a_c}(1{+}y)/\theta(1{+}y) = \mu_M(cy), 
$$ 
for all $y\in E$ such that $\ups_E(y) \ge p^rk{-}w$ and all $a_c\in \frak p^{m-p^rk}$ such that 
$$ 
(1{+}a_c)^{-1}\alpha(1{+}a_c) \equiv \alpha{+}c \pmod{\frak p^{w-p^rk}}. 
$$ 
\item $\Psi_\vT(k) < \biP EF\vs(k)$. 
\item $\Psi_\vT(x) < \biP EF\vs(x)$ for all $x$, $j_\infty <x < \bar\jmath_\infty$. 
\item 
For any twisting datum $(h,d,\phi)$, where $h\in I$ is an integer at which $\Psi_\vT$ is smooth, we have 
$$ 
\theta^{1+a_d}(1{+}y)/\theta(1{+}y) = \mu_M(dy), 
$$ 
for all $y\in E$ such that $\ups_E(y) \ge p^rh{-}w$ and all $a_d\in \frak p^{m-p^rh}$ such that 
$$ 
(1{+}a_d)^{-1}\alpha(1{+}a_d) \equiv \alpha{+}d \pmod{\frak p^{w-p^rh}}. 
$$ 
\endroster 
\endproclaim 
\demo{Proof} 
The equivalence of (1) and (2) is the preceding Corollary. The equivalence of (2) and (3) is 5.6 Corollary 2. 
\par 
Certainly (4) implies (1), so suppose that (4) fails: there is a twisting datum $(h,d,\phi)$ such that $\theta^{1+a_d}(1{+}y)/\theta(1{+}y) \neq  \mu_M(dy)$, for some $y\in \frak p_E^{p^rh-w}$. Thus $\Bbb A(\phi\vT,\vT) = \Psi_\vT(h) = h{-}p^{-r}w = \biP EF\vs(h)$. Corollary 2 of 5.6 now implies that (3) fails. \qed 
\enddemo 
\subhead 
5.10 
\endsubhead 
We start the proofs of the parts of 5.3 Theorem that allow $\biP EF\vs$ to have an even number of jumps, the case of an odd number of jumps having been dispatched in 5.6. In this sub-section, we prove part (2) of the theorem. 
\proclaim{Proposition} 
Suppose $m>2w_{E/F}$. If $\theta\in \scr C(\frak a,\alpha)$ has endo-class $\vT$, then $\Psi_\vT(x) = \upr{2\,}\Psi_{(E/F,\vs)}(x)$, $0\le x\le \vs_\vT$. 
\endproclaim 
\demo{Proof} 
If $\biP EF\vs$ has an odd number of jumps, the result follows from part (1) of the theorem, proved in 5.6. We therefore assume that $\biP EF\vs$ has an even number of jumps and continue with the notation of 5.8, 5.9. As argued at the beginning of 5.9, we may assume that $\Psi_\vT(\bar \jmath_\infty) = \biP EF\vs(\bar\jmath_\infty) \le \vs/2$. 
\par 
 The line segment $y=x{-}w/p^r = \biP EF\vs(x)$, $x\in I$, crosses the axis of symmetry $x{+}y=\vs$ where $p^rx = (m{+}w)/2$. So, we choose an integer $k$, at which $\Psi_\vT$ is smooth, to satisfy $j_\infty < k<(m{+}w)/2p^r$. That is, 
$$ 
m{-}p^rk > (m{-}w)/2 > w/2. 
\tag 5.10.1 
$$ 
Let $(k,c,\chi)$ be a twisting datum over $F$. Define $a_c$ as in 5.9 Proposition 2. We apply the definition (5.8.3), with $\ups_E(y) \ge p^rk{-}w$, to get  
$$ 
\align 
\theta(1{-}\zeta c\alpha^{-1} y)\,\mu_M(-\alpha\zeta\eps y) &= \theta(1{-}\zeta c\alpha^{-1} y)\, \mu_M(\alpha\zeta c\alpha^{-1} y)\, \mu_M(\alpha\zeta a_cc\alpha^{-1} y) \\ 
&= \mu_M(-\alpha\zeta c\alpha^{-1} y)\, \mu_M(\alpha\zeta c\alpha^{-1} y)\, \mu_M(\alpha\zeta a_c c\alpha^{-1} y) \\ 
&= \mu_M(a_cc\zeta y). 
\endalign 
$$ 
So, by 5.9 Proposition 1, 
$$ 
\theta^{1+a_c}(1{+}y)/\theta(1{+}y) = \mu_M(a_c c\zeta y). 
$$ 
We show that the character 
$$ 
1{+}y \longmapsto \mu_M((1{-}\zeta a_c)cy),\quad y\in \frak p_E^{p^rk-w}, 
\tag 5.10.2 
$$ 
is not trivial, for some choice of $c\in \frak p_F^{-k}/\frak p_F^{1-k}\smallsetminus \{0\}$. The proposition will then follow from 5.9 Proposition 2. 
\par 
The defining relation $A_\alpha(a_c) \equiv (1{+}a_c)c\alpha^{-1} \pmod{\frak p^{m-p^rk+w}}$ (5.8.2) implies $\ups_E(s_{E/F}(a_c)) \ge w$. If $\ups_E(s_{E/F}(a_c)) > w$, (5.10.2) reduces to $1{+}y \mapsto \mu_M(cy)$, which is surely not trivial. We therefore assume that $s_{E/F}(a_c)$ has valuation $w$ for all $c \in \frak p_F^{-k}/\frak p_F^{1-k}$, $c\neq 0$. We show that this hypothesis is untenable.
\par 
We put $a_0 = 0$ and let $c,c'\in \frak p_F^{-k}/\frak p_F^{1-k}$. Conjugating the defining relation 
$$ 
(1{+}a_c)^{-1}\alpha(1{+}a_c) \equiv \alpha{+}c \pmod{\frak p^{w-p^rk}} 
$$ 
by $(1{+}a_{c'})$, condition (5.10.1) yields 
$$ 
\left. \aligned 
a_{c+c'} &\equiv a_c+a_{c'} \pmod{\frak p^{w+1}} , \\ 
s_{E/F}(a_{c+c'}) &\equiv s_{E/F}(a_c+a_{c'}) \pmod{\frak p_E^{w+1}},  
\endaligned \right\}
 \qquad c,c' \in \frak p_F^{-k}/\frak p_F^{1-k}. 
$$ 
Thus $c\mapsto s_{E/F}(a_c)$ is a homomorphism $\frak p_F^{-k}/\frak p_F^{1-k} \to \frak p_E^w/\frak p_E^{1+w}$. By 5.8 Lemma, $s_{E/F}(1{+}a_c) \notin \frak p_E^{1+w}$. That is, the non-zero element $-s_{E/F}(1)$ of $\frak p_E^w/\frak p_E^{1+w}$ is not of the form $s_{E/F}(a_c)$. So, the homomorphism $\frak p_F^{-k}/\frak p_F^{1-k} \to \frak p_E^w/\frak p_E^{1+w}$, $c\mapsto s_{E/F}(a_c)$, cannot be surjective. It therefore has a non-trivial kernel, contradicting our hypothesis, and the proposition follows. \qed 
\enddemo 
\subhead 
5.11 
\endsubhead 
We prove part (3) of 5.3 Theorem. 
\proclaim{Proposition} 
Suppose $m/2 < w< m$ and that $w\equiv 0 \pmod p$. Let $\theta\in \scr C(\frak a,\alpha)$ satisfy $l_E(\theta) \le m{-}w$. If $\vT$ is the endo-class of $\theta$, then $\Psi_\vT(x) = \upr{2\,}\Psi_{(E/F,\vs)}(x)$, $0\le x\le \vs$. 
\endproclaim 
\demo{Proof} 
We may again assume that $\biP EF\vs$ has an even number of jumps and proceed as before. In the formula (5.9.1), 
$$ 
\theta^{1+a}(1{+}y)/\theta(1{+}y) = \theta(1{-}\zeta c\alpha^{-1} y)\,\mu_M(-\alpha\zeta\eps y), 
$$ 
we have $\zeta \equiv 0 \pmod{\frak p_F}$, so it reduces to 
$$ 
\theta^{1+a}(1{+}y)/\theta(1{+}y) = 1 \neq \chi\circ\N EF(1{+}y) = \mu_M(cy), 
$$ 
for some choice of $y \in \frak p_E^{p^rk-w}$. The result now follows from 5.9 Proposition 2. \qed 
\enddemo 
\subhead 
5.12 
\endsubhead 
We prove part (4) of 5.3 Theorem. Thus $\biP EF\vs$ has an even number of jumps, and we may continue with the notation of 5.8, 5.9. In particular, $I$ is the interval $j_\infty < x < \bar\jmath_\infty$. Here, the element $\zeta$ of 5.9 Proposition 1 is a unit in $F$. We have to prove: 
\proclaim{Proposition} 
Suppose that $m>w\ge m/2$ and $w\not\equiv 0 \pmod p$. Assume that $\biP EF\vs$ has an even number of jumps. There is a unique character $\phi$ of $U^{m-w}_E/U^{1+m-w}_E$ with the following property: a character $\theta\in \scr C(\frak a,\alpha)$, with $l_E(\theta) \le m{-}w$ and endo-class $\vT$, satisfies $\Psi_\vT = \biP EF\vs$ if and only if $\theta\Mid U^{m-w}_E \neq \phi$. 
\par 
If $\theta\Mid U^{m-w}_E = \phi$, then $\Psi_\vT(x) < \biP EF\vs(x)$ for all $x\in I$. 
\endproclaim 
\demo{Proof} 
We write out again the formula (5.9.1) 
$$ 
\theta^{1+a_c}(1{+}y)/\theta(1{+}y) = \theta(1{-}\zeta c\alpha^{-1} y)\,\mu_M(-\alpha\zeta\eps_c y), 
\tag 5.12.1 
$$ 
where $a_c\in \frak p^{m-p^rk}$ is a solution of the congruence (5.8.1): 
$$ 
(1{+}a_c)^{-1}\alpha(1{+}a_c) \equiv \alpha{+}c \pmod{\frak p^{w-p^rk}} 
\tag 5.12.2 
$$ 
and $\eps_c$ is given by (5.8.3), relative to the element $a_c$. 
\par 
We use (5.8.3) to re-write the last factor in (5.12.1) as 
$$ 
\mu_M(-\alpha\zeta\eps_c y) = \mu_M(\zeta cy)\, \mu_M(\zeta a_ccy). 
$$ 
For $1{+}y \in U^{p^rk-w}_E/U^{1+p^rk-w}_E$, write 
$$ 
\Xi_{\theta,c}(1{+}y) = \theta(1{-}\alpha^{-1}\zeta cy)\, \mu_M(\zeta cy)\, \mu_M(\zeta a_ccy)\, \mu_M(-cy).
$$ 
That is, $\Xi_{\theta,c}(1{+}y)$ is the product of the right hand side of (5.12.1) and $\mu_M(-cy)$. Therefore, invoking 5.9 Proposition 2, we have: 
\proclaim{Lemma} 
The character $\Xi_{\theta,c}$ is trivial if and only if $\Psi_\vT(x) < \biP EF\vs(x)$ for all $x\in I$. This condition holds for one element $c\in \frak p_F^{-k}/\frak p_F^{1-k} \smallsetminus \{0\}$ if and only if it holds for all.  
\endproclaim  
Write $z = \zeta\alpha^{-1}cy$. Thus $y\mapsto z$ induces an isomorphism of $U^{p^rk-w}_E/U^{1+p^rk-w}_E$ with $U^{m-w}_E/U^{1+m-w}_E$ and 
$$ 
\align 
\Xi_{\theta,c}(1{+}y) &= \theta(1{-}z)\,\mu_M(\alpha z) \mu_M(\alpha a_c   z) \mu_M(-\zeta^{-1}\alpha z) \\ 
&= \theta(1{+}z)^{-1}\,\mu_M\big((1 +  a_c - \zeta^{-1})\alpha z\big). 
\endalign 
$$ 
Since $\ups_E(s_{E/F}(a_c)) \ge w$ (5.8 Lemma), the formula 
$$ 
\xi_c(1{+}z) = \mu_M((1 +  a_c - \zeta^{-1})\alpha z) 
$$ 
defines a character of $U^{m-w}_E/U^{1+m-w}_E$ which is independent of the character $\theta \in \scr C(\frak a,\alpha)$ such that $l_E(\theta) \le m{-}w$. For fixed $\theta$, the character $\theta^{-1}\xi_c\Mid U^{m-w}_E$ is either trivial for all $c\in \frak p_F^{-k}/\frak p_F^{1-k} \smallsetminus \{0\}$, or else it is non-trivial for all such $c$, by the lemma. Given any character $\phi$ of $U^{m-w}_E/U^{1+m-w}_E$, there exists $\theta\in \scr C(\frak a,\alpha)$ agreeing with $\phi$ on $U^{m-w}_E$. We conclude that, if $c,c'\in \frak p_F^{-k}/\frak p_F^{1-k} \smallsetminus \{0\}$, then $\xi_c = \xi_{c'}$. 
\par 
The proposition therefore holds for the character 
$$ 
\phi(1{+}z) = \xi_c(1{+}z) = \mu_M((1 +  a_c - \zeta^{-1})\alpha z) , \quad z\in \frak p_E^{m-w}, 
\tag 5.12.3 
$$ 
for any non-trivial element $c$ of $\frak p_F^{-k}/\frak p_F^{1-k}$. \qed 
\enddemo 
We have completed the proof of 5.3 Theorem. \qed 
\remark{Remark} 
We have noted that the character $\phi$ of (5.12.3) does not depend on the parameter $c\in \frak p_F^{-k}/\frak p_F^{1-k} \smallsetminus \{0\}$. Indeed, any twisting datum $(h,b,\xi)$ with $j_\infty < h < \bar\jmath_\infty$ will, by 5.9 Proposition 2, give rise to the same character. 
\endremark 
\head\Rm 
6. Variation of parameters 
\endhead 
In section 5, we fixed the stratum $[\frak a,m,0,\alpha]$ and calculated $\Psi_\vT$, in many cases, under the restriction (5.3.1). Here, we investigate the scope for changing the stratum {\it without\/} changing the set $\scr C(\frak a,\alpha)$, in order to avoid the condition (5.3.1) and to clarify the dichotomy in part (4) of the theorem. 
\subhead 
6.1 
\endsubhead 
Let $[\frak a,m,0,\alpha]$ be a simple stratum in $M = \M{p^r}F$, $r\ge 1$, satisfying the usual conditions: $F[\alpha]/F$ is totally ramified of degree $p^r$ and $m$ is not divisible by $p$. Set $\vs = m/p^r$. Define $\roman P(\frak a,\alpha)$ as the set of $\beta\in \GL{p^r}F$ for which $[\frak a,m,0,\beta]$ is a simple stratum such that $\scr C(\frak a,\beta) = \scr C(\frak a,\alpha)$. We summarize the main properties of such elements $\beta$. As usual, $\frak p$ is the Jacobson radical of $\frak a$. 
\proclaim{Proposition} 
Write $E = F[\alpha]$. 
\roster 
\item 
If $\beta\in \roman P(\frak a,\alpha)$, the field extension $F[\beta]/F$ is totally ramified of degree $p^r$. Moreover, 
$$ 
\beta\equiv \alpha \pmod{\frak p^{-[m/2]}}. 
$$ 
\item 
Let $k\le [m/2]$ be an integer and let $b\in \frak p^{-k}$. The element $\beta = \alpha{+}b$ then lies in $\roman P(\frak a,\alpha)$. 
\endroster 
\endproclaim 
\demo{Proof} 
In (1), the first assertion is an instance of \cite{15} (3.5.1), while the second follows from the definition in 2.3. In (2), the stratum $[\frak a,m,0,\beta]$ is simple with the required properties, as follows from \cite{15} (2.2.3). \qed 
\enddemo 
\remark{Remarks} 
\roster 
\item 
Any element of $\roman P(\frak a,\alpha)$ arises as in part (2) of the proposition \cite{15} (2.4.1).  
\item 
If $\beta\in \roman P(\frak a,\alpha)$, then $\roman P(\frak a,\beta) = \roman P(\frak a,\alpha)$. 
\item 
Let $[\frak a,q,0,\gamma]$ be a simple stratum in $M$. If the set $\scr C(\frak a,\alpha)\cap \scr C(\frak a,\gamma)$ is non-empty, then $q=m$ and  $\scr C(\frak a,\alpha) = \scr C(\frak a,\gamma)$ \cite{15} (3.5.8), (3.5.11), whence $\|\scr C(\frak a,\alpha)\| = \|\scr C(\frak a,\gamma)\|$. 
\item 
Let $\scr K_\frak a$ be the group of $x\in \GL{p^r}F$ such that $x\frak ax^{-1} = \frak a$. For $\beta_1,\beta_2 \in \roman P(\frak a,\alpha)$, say that $\beta_1\sim \beta_2$ if $\beta_1U^m_\frak a$ is $\scr K_\frak a$-conjugate to $\beta_2U^m_\frak a$. It is shown in \cite{16} that the sets $\roman P(\frak a,\alpha)/\!\sim\,$ and $\|\scr C(\frak a,\alpha)\|$ are in (non-canonical) bijection. 
\endroster 
\endremark 
\subhead 
6.2 
\endsubhead 
We give a first application of this concept. 
\proclaim{Proposition} 
Suppose that $m>2w_{F[\alpha]/F}$. If $\beta\in \roman P(\frak a,\alpha)$ then $w_{F[\beta]/F} = w_{F[\alpha]/F}$. 
\endproclaim 
\demo{Proof} 
Let $\frak p = \roman{rad}\,\frak a$. Abbreviate $w_\alpha = w_{F[\alpha]/F}$, $w_\beta = w_{F[\beta]/F}$. By hypothesis, $m{-}w_\alpha \ge \left[\frac m2\right]{+}1$. A character $\theta\in \scr C(\frak a,\alpha)$, by definition (2.3.1), agrees with $\mu_M*\alpha$ on the group $H^{1+[m/2]}(\alpha,\frak a) = U^{1+[m/2]}_\frak a$. The integer $l_E(\theta)$ is the least integer $k\ge0$ such that $\theta$ is trivial on $U^{1+k}_EU^{1+m}_\frak a = 1{+}\frak p_E^{1+k}{+}\frak p^{1+m}$. In this case, $l_E(\theta) = m{-}w_\alpha$, as in 5.2 Proposition.  However, $\frak p_E^{1+k}{+}\frak p^{1+m}$ is the kernel of the adjoint map $A_\alpha$ on $\frak p^{1+k}/\frak p^{1+m}$, so $l_E(\theta)$ is the least integer $k\ge0$ such that $\mu_M*\alpha$ is trivial on $1+\roman{Ker}\, A_\alpha\Mid \frak p^{1+k}/\frak p^{1+m}$. The same analysis applies relative to $\beta$ in place of $\alpha$. 
\par 
By hypothesis, $\theta$ also agrees with $\mu_M*\beta$ on $U^{1+[m/2]}_\frak a$, so $\beta\equiv \alpha \pmod{\frak p^{-[m/2]}}$. The maps $A_\alpha$, $A_\beta$ therefore agree on the group $\frak p^{1+[m/2]}/\frak p^{1+m}$, and the result follows. \qed 
\enddemo 
The following corollary does not form part of the main development, but is included to illuminate the division into cases in 5.3 Theorem: see Example below. 
\proclaim{Corollary} 
In the context of the proposition, we have 
$$ 
\align 
\biP {F[\beta]}F\vs(x) &= \biP {F[\alpha]}F\vs(x),  \quad 0\le x\le \vs, \\ 
\psi_{F[\beta]/F} &= \psi_{F[\alpha]/F}. 
\endalign 
$$ 
If $p\ge3$, the function $\Psi_\vT$ has an even number of jumps. 
\endproclaim 
\demo{Proof} 
Let $\vT \in \|\scr C(\frak a,\alpha)\|$. Part (2) of 5.3 Theorem gives 
$$ 
\Psi_\vT(x) = \upr{2\,}\Psi_{(F[\alpha]/F,\vs)}(x) = \upr{2\,}\Psi_{(F[\beta]/F,\vs)}(x), \quad 0\le x\le \vs, 
\tag 6.2.1 
$$ 
whence follows the first assertion. 
\par 
If $\Psi$ denotes any of the functions appearing in (6.2.1), define $c$ by $c{+}\Psi(c) = \vs$, so that $\psi_{F[\alpha]/F}(x) = \psi_{F[\beta]/F}(x)$ for $0\le x\le c$. Let $j_\infty$ be the last jump of $\psi_{F[\alpha]/F}$. 
\proclaim{Lemma} 
If $j_\infty < c$, then $\psi_{F[\beta]/F} = \psi_{F[\alpha]/F}$ and $\Psi_\vT$ has an even number of jumps. 
\endproclaim 
\demo{Proof} 
The second assertion follows from 4.2 Proposition and 5.3 Theorem (2). For $0< x < c$, we have $p^r\Psi_\vT(x) = \psi_{F[\alpha]/F}(x) = \psi_{F[\beta] /F}(x)$. Thus $\psi_{F[\beta]/F}$ has a jump at $j_\infty$ and $\psi'_{F[\beta]/F}(x) = p^r$, for $j_\infty < x <c$. Therefore $j_\infty$ is the last jump of $\psi_{F[\beta]/F}$ ({\it cf\.} 1.6 Proposition) and the lemma follows. \qed 
\enddemo 
We use 1.6 Corollary:  
$$ 
j_\infty\le \frac{w_\alpha}{p^{r-1}(p{-}1)},   
$$ 
where $w_\alpha = w_{F[\alpha]/F}$. Invoking also 1.6 Proposition, we get 
$$ 
\align 
j_\infty+p^{-r} \psi_{F[\alpha]/F}(j_\infty) &= 2j_\infty - p^{-r}w_\alpha  
\le \frac{2w_\alpha}{p^{r-1}(p{-}1)} -\frac {w_\alpha}{p^r} \\ 
&=\frac{w_\alpha}{p^r}\, \frac{p{+}1}{p{-}1}  
< \frac m{2p^r}\,\frac{p{+}1}{p{-}1} = \vs\,\frac{p{+}1}{2(p{-}1)}, 
\endalign 
$$ 
since $w_\alpha < m/2$. So, if $p\ge3$, the point $(j_\infty,p^{-r} \psi_{F[\alpha]/F}(j_\infty))$ lies strictly below the line $x{+}y=\vs$. That is, $j_\infty < c$ and, in this case, the corollary follows from the lemma. 
\par 
Suppose therefore that $p=2$. If $r=1$, the graph of $\psi_{F[\alpha]/F}$ consists of segments of the two lines $y=x$ and $y = 2x{-}w_{F[\alpha]/F}$. Likewise for $\beta$. The proposition gives $w_\alpha = w_\beta$, whence the result in this case. 
\par 
Consider next the case where $r\ge2$ and $j_\infty$ is the only jump of $\psi_{F[\alpha]/F}$. Here, $j_\infty = w_\alpha/(2^r{-}1)$ (1.6 Corollary) and so $2^{-r} \psi_{F[\alpha]/F}(j_\infty) = 2^{-r}j_\infty$. Therefore  
$$ 
j_\infty+2^{-r} \psi_{F[\alpha]/F}(j_\infty) = \frac {w_\alpha}{2^r}\,\frac{2^r{+}1}{2^r{-}1} < \frac\vs2\,\frac{2^r{+}1}{2^r{-}1} \le \vs. 
$$ 
Thus $j_\infty{+}2^{-r} \psi_{F[\alpha]/F}(j_\infty) < \vs$, so $j_\infty < c$, and the result in this case also follows from the lemma. 
\par 
We are left with the case where $r\ge2$ and $\psi_{F[\alpha]/F}$ has at least two jumps. If $j_\infty < c$, there is nothing more to do, so we assume $j_\infty \ge c$. Let $j'$ be the penultimate jump of $\psi_{F[\alpha]/F}$. In particular,  
$$ 
2^{-r}\psi_{F[\alpha]/F}(x) \le x/4, \quad 0\le x\le j'. 
$$ 
We show that $j'<c$. 
\par 
Abbreviate $a = w_\alpha/2^r$, so that $2^{-r} \psi_{F[\alpha]/F}(x) = x{-}a$ for $x\ge j_\infty$, while $2^{-r}\psi_{F[\alpha]/F}(x) > x{-}a$ when $0\le x<j_\infty$. Thus 
$$ 
x{-}a <2^{-r}\psi_{F[\alpha]/F}(x) \le x/4, \quad 0\le x\le j'. 
$$ 
The lines $y=x{-}a$, $y=x/4$ meet at the point $(4a/3,a/3)$, so $j' < 4a/3$. Since $\frac{4a}3+\frac a3 = \frac{5a}3 < 2a<\vs$, this point of intersection lies below the line $x{+}y = \vs$. Therefore $j'< 4a/3 < c$.  
\par 
We have $\psi_{F[\alpha]/F}(x) = \psi_{F[\beta]/F}(x)$ in the region $0\le x < c$. The same analysis applies with $\beta$ replacing $\alpha$, so $j'$ is also the penultimate jump of $\psi_{F[\beta]/F}$. Let $\tilde\psi(x)$ be the piecewise linear function agreeing with $\psi_{F[\alpha]/F}(x) = \psi_{F[\beta]/F}(x)$ for $x<c$ and smooth for $x> j'$. In the region $x\ge 0$, we then have 
$$ 
\align 
\psi_{F[\alpha]/F}(x) &= \roman{max}\{\tilde\psi(x),x{-}2^{-r}w_\alpha\} \\ 
&= \roman{max}\{\tilde\psi(x),x{-}2^{-r}w_\beta\} = \psi_{F[\beta]/F}(x), 
\endalign 
$$ 
as required. \qed 
\enddemo 
\example{Example} 
Suppose $p = 2$ and let $\vT \in \ewc F$ have degree $2$. Thus $\Psi_\vT$ has an odd number of jumps (in fact one jump) if and only if $j_\infty \ge c$, using the notation of the Corollary. By 3.8 Proposition, this is equivalent to $m\le 3w$ ({\it cf\.} Kutzko \cite{33}). So, for $p^r =2$, there are examples of endo-classes $\vT$ for which $m>2w$ while $\Psi_\vT$ has an odd number of jumps. 
\endexample 
\subhead 
6.3 
\endsubhead 
We use the notation from the start of 6.1, except that we write $E = F[\alpha]$ and assume $m\le 2w_{E/F}$. This case is more complex and interesting. We first investigate the possibility of changing $\alpha$ to {\it raise\/} the exponent $w_{E/F}$. 
\par 
Let $\frak p$ be the Jacobson radical $\roman{rad}\,\frak a$ of $\frak a$. Let $ s_{E/F}:M\to E$ be a tame corestriction. 
\proclaim{Proposition}
Suppose that $m/2\le w_{E/F} < m$ and that $w  = w_{E/F} \not\equiv 0 \pmod p$. Let $\zeta = (w{-}m)/m \in \frak o_F$. There exists $b\in \frak p^{w-m}$ such that 
$$ 
(\zeta{+}1) s_{E/F}(b) \equiv s_{E/F}(\alpha) \pmod{\frak p_E^{1+w-m}}. 
\tag 6.3.1 
$$ 
For any such $b$, the element $\beta =  \alpha{-}b$ lies in $\roman P(\frak a,\alpha)$ and $w_{F[\beta]/F} > w$. 
\endproclaim 
\demo{Proof} 
The hypothesis $w\not\equiv 0\pmod p$ implies that $\zeta\not\equiv -1 \pmod{\frak p_F}$. The exact sequences (5.1.1) then give an element $b$ with the necessary properties. 
\par 
The hypothesis $m\le 2w$ implies $\beta \equiv \alpha \pmod{\frak p^{-[m/2]}}$ and, following 6.1 Proposition, $\beta\in \roman P(\frak a,\alpha)$. Write $E' = F[\beta]$. 
\proclaim{Lemma} 
\roster 
\item 
Let $y\in \frak p_E^{m-w}$. There exist $y'\in \frak p_{E'}^{m-w}$ and $h\in \frak p^m$ so that 
$$ 
y = y'+h, 
\tag 6.3.2 
$$ 
The map $y\mapsto y'$ induces an isomorphism $\frak p_E^{m-w}/\frak p_E^m \to \frak p_{E'}^{m-w}/\frak p_{E'}^m$. 
\item 
The decomposition \rom{(6.3.2)} may be chosen so that, additionally, 
$$ 
h \equiv \zeta by\beta^{-1} \pmod{\frak p^{m+1}{+}A_\beta(\frak p^m)}. 
\tag 6.3.3 
$$ 
\endroster 
\endproclaim  
\demo{Proof} 
In (1), the relation $H^1(\beta,\frak a) = H^1(\alpha,\frak a)$ implies that any $y\in \frak p^{m-w}_E$ takes the form $y = y'{+}h$, with $y'\in \frak p^{m-w}_{E'}$ and $h\in \frak p^{1+[m/2]}$. The element $[\beta,y] = [\beta,h] = -[b,y]$ lies in $\frak a$. By (5.1.1), we may choose the decomposition so that $h\in \frak p^m$. The second assertion is immediate. 
\par 
In (2), 5.7 Proposition gives 
$$ 
[\beta,h] = -[b,y] = A_y(b)y \equiv \zeta A_\beta(b)y \pmod{\frak p{+}A^2_\beta(\frak a)}. 
$$ 
So, we may further refine (6.3.2) to get (6.3.3). \qed 
\enddemo 
In multiplicative terms, the definition of $\beta$ gives $\beta\equiv \alpha \pmod{U^w_\frak a}$ and therefore $\beta^{-1}\equiv \alpha^{-1} \pmod{\frak p^{m+w}}$. It follows that  
$$ 
\aligned 
\zeta by\alpha^{-1} &\equiv \zeta by\beta^{-1} \pmod{\frak p^{m+1}}, \quad \text{and}\\ 
A_\alpha(\frak p^m){+}\frak p^{m+1} &= A_\beta(\frak p^m){+}\frak p^{m+1}. 
\endaligned 
\tag 6.3.4 
$$  
The relation (6.3.2) gives $\mu_M(\beta y) = \mu_M(\beta y')\,\mu_M(\beta h)$, while (6.3.3), (6.3.4) yield $\mu_M(\beta h) = \mu_M(\zeta by)$. On the other hand, $\mu_M(\beta y) = \mu_M((\alpha{-}b)y)$ by definition, so 
$$
\mu_M(\beta y') = \mu_M\big((\alpha {-} (\zeta{+}1)b)y\big) = \mu_E(s_{E/F}(\alpha-(\zeta{+}1)b)y) = 1, 
$$ 
for all $y\in \frak p_E^{m-w}$, by (6.3.1). Part (1) of the lemma now shows that $\mu_M(\beta y') = 1$ for all $y'\in \frak p_{E'}^{m-w}$. Therefore $w_{E'/F}>w$, as required. \qed 
\enddemo 
\proclaim{Corollary} 
Suppose that $m\le 2w_{E/F}$. There exists $\beta = \alpha{-}b \in \roman P(\frak a,\alpha)$, where $b\in \frak p^{w_{E/F}-m}$ satisfies \rom{(6.3.1),} with the following property. If $E' = F[\beta]$, then either 
\roster 
\item $w_{E'/F}\ge w_{E/F}$ and $w_{E'/F} \equiv 0 \pmod p$, or 
\item $w_{E'/F} \ge m$. 
\endroster 
\endproclaim 
\demo{Proof} 
If $w_{E/F}$ is divisible by $p$, there is nothing to do. Otherwise, we construct $E_1 = F[\beta]$ following the proposition. If either $w_{E_1/F} \ge m$ or $w_{E_1/F} \equiv 0 \pmod p$, we are finished. So, assume that  $w_{E_1/F} < m$ and $w_{E_1/F} \not\equiv 0 \pmod p$. Set $w_1 = w_{E_1/F}$. Following the procedure as before, we construct an element $\gamma \equiv \beta \pmod{\frak p^{w_1-m}}$ such that $w_{F[\gamma]/F} > w_1$. The congruence condition on $\gamma$ ensures that $b_1 = \alpha{-}\gamma$ satisfies (6.3.1). We iterate this procedure as necessary until we achieve either (1) or (2). \qed 
\enddemo 
\subhead 
6.4 
\endsubhead 
We retain the notation of 6.3, in particular $E = F[\alpha]$ and $w = w_{E/F}$. The elements $\beta$ of 6.3 Proposition have useful properties relative to certain simple characters. 
\proclaim{Proposition} 
Let $\beta = \alpha{-}b \in \roman P(\frak a,\alpha)$, where $b\in \frak p^{w-m}$ satisfies \rom{(6.3.1).} 
\roster 
\item 
If $\xi\in \scr C(\frak a,\alpha)$ satisfies $\xi(1{+}y) = \mu_M(\alpha y)$, $y\in \frak p_E^{m-w}$, then 
$$ 
l_{F[\beta]}(\xi) = l_E(\xi) = m{-}w. 
$$ 
\item 
If $\xi\in \scr C(\frak a,\alpha)$ satisfies $l_E(\xi) > m{-}w$, then $l_{F[\beta]}(\xi) = l_E(\xi)$. 
\endroster 
\endproclaim 
\demo{Proof} 
For part (1), we use (6.3.4) to evaluate  
$$ 
\xi(1{+}y') = \mu_M(\alpha y)\mu_M(\alpha h) = \mu_M((\alpha{-}\zeta b)y),  
$$ 
where $y\in \frak p_E^{m-w}$, $y'\in \frak p_{E'}^{m-w}$ and $h\in \frak p^m$ are related as in 6.3 Lemma. As $1{+}\zeta$ is a unit of $\frak o_F$, we have $\zeta s_{E/F}(b) \equiv \zeta s_{E/F}(\alpha)(\zeta{+}1)^{-1} \pmod{\frak p_E^{1+w-m}}$ and so, by (6.3.1), 
$$ 
\xi(1{+}y') = \mu_E\big((\zeta{+}1)^{-1}s_{E/F}(\alpha)y\big), \quad y'\in \frak p_{F[\beta]}^{m-w}. 
$$ 
We may choose $y$ so that $\xi(1{+}y') \neq 1$ and part (1) of the proposition follows. 
\par 
In part (2), let $l = l_E(\xi)$. Let $y\in \frak p_E^l$. Since $l> m{-}w$, we use 6.3 Lemma to write $y = y'{+}h$, where $y'\in \frak p_{F[\beta]}^l$ and $h\in \frak p^{m+1}$. Thus $\xi(1{+}y') = \xi(1{+}y)$ and we may choose $y$ so that $\xi(1{+}y) \neq 1$. If, however, $y\in \frak p_E^{1+l}$, then $\xi(1{+}y') = \xi(1{+}y) = 1$, so $l_{F[\beta]}(\xi) = l$, as required. \qed 
\enddemo 
Note the very restrictive hypothesis on $\xi$ in this corollary. 
\subhead 
6.5 
\endsubhead 
We turn to the question of {\it lowering\/} of the exponent $w_{E/F}$. Following 6.2 Proposition, we are restricted to the case where $2w_{E/F} > m$. The consequences for simple characters are complementary to those of 6.4, but we get much more detail. 
\proclaim{Theorem} 
Let $[\frak a,m,0,\alpha]$ be a simple stratum in $M = \M{p^r}F$ in which $E = F[\alpha]/F$ is totally ramified of degree $p^r$ and $p$ does not divide $m$. Suppose $m<2w_{E/F}$. Let $d$ be an integer such that 
$$ \gathered
1\le d \le m/2, \quad d> \roman{max}\,\{0,m{-}w_{E/F}\}, \\ d\not\equiv m \pmod p. \endgathered 
\tag 6.5.1  
$$ 
Let $b\in \frak p^{-d}$ satisfy $\ups_E(s_{E/F}(b)) = -d$. The element $\beta = \alpha{+}b$ lies in $\roman P(\frak a,\alpha)$ and 
$$ 
w_{E'/F} = m{-}d < w_{E/F}, \quad E' = F[\beta]. 
\tag 6.5.2 
$$ 
Let $\theta \in \scr C(\frak a,\alpha)$ and write $l = l_E(\theta)$. For any such $\beta$, the following hold. 
\roster 
\item Suppose $l<d$. 
\itemitem{\rm (a)} If $d\not\equiv 0 \pmod p$, then $l_{E'}(\theta) = d$. 
\itemitem{\rm (b)} If $d\equiv 0 \pmod p$, then $l_{E'}(\theta) < d$. 
\item If $l>d$, then $l_{E'}(\theta) = l$. 
\item Suppose $l = d$. 
\itemitem{\rm (a)} If $d\not\equiv 0 \pmod p$, then $l_{E'}(\theta) \le d$, with both equality and inequality occurring. 
\itemitem{\rm (b)} If $d\equiv 0 \pmod p$, then $l_{E'}(\theta) = d$. 
\endroster 
\endproclaim 
\demo{Proof} 
Writing $\frak p = \roman{rad}\,\frak a$, let $b\in \frak p^{-d}$ satisfy $\ups_E(s_{E/F}(b)) = -d$. As 6.1 Proposition, the element $\beta = \alpha{+}b$ lies in $\roman P(\frak a,\alpha)$. Put $E' = F[\beta]$. 
\proclaim{Lemma} 
\roster 
\item 
Let $y\in \frak p_E^d$. There exist $y'\in \frak p_{E'}^d$ and $h\in \frak p^m$ such that $y=y'{+}h$. The map $y\mapsto y'$ induces an isomorphism $\frak p_E^d/\frak p_E^m \to \frak p_{E'}^d/\frak p_{E'}^m$. 
\item 
If $y\in \frak p_E^{d+1}$, then $y'\in \frak p_{E'}^{d+1}$ and one may take $h\in \frak p^{m+1}$. 
\endroster 
\endproclaim 
\demo{Proof} 
This is identical to the proof of part (1) of 6.3 Lemma, so we omit the details. \qed 
\enddemo 
Set $w' = w_{E'/F}$. We first show that $\mu_M(\beta z) = 1$, for $z\in \frak p_{E'}^{1+d}$. By the lemma, there exist $y\in \frak p_E^{1+d}$ and $h\in \frak p^{m+1}$ such that $y = z{+}h$. The condition $d>m{-}w_{E/F}$ implies $\mu_M(\alpha y) = 1$. As $by\in \frak p$, so $\mu_M(by) = 1$. Altogether, $\mu_M(\beta y) = \mu_M(\alpha y) \mu_M(by) = 1$. Therefore $1= \mu_M(\beta z) \mu_M(\beta h)  = \mu_M(\beta z)$, as asserted. It follows that $d\ge m{-}w'$. 
\par 
Now take $z\in E'$ with $\ups_{E'}(z) = d$. Thus $z=y{-}h$, where $y\in E$ satisfies $\ups_E(y) = d$ and $h\in \frak p^m$. Consequently, $[\beta,h] = [\beta,y] = [b,y]$. Setting $\zeta = -d/m$, 5.7 Proposition gives  
$$ 
[b,y] = -A_y(b)y \equiv -\zeta A_\alpha(b)y \pmod{A_\alpha^2(\frak a)+\frak p}. 
$$ 
Since $\alpha \equiv \beta \pmod{U_\frak a^{m-d}}$, we have 
$$ 
A_\alpha(a) \equiv A_\beta(a) \pmod{\frak p^{k+m-d}}, \quad a\in \frak p^k, 
$$ 
for any integer $k$. So, 
$$ 
[\beta,h] = [b,y] \equiv -\zeta A_\beta(b)y \pmod{A_\beta^2(\frak a)+\frak p}. 
$$ 
We may therefore choose the decomposition $y=z{+}h$ so that 
$$ 
h\equiv -\zeta by\beta^{-1} \equiv -\zeta by\alpha^{-1} \pmod{A_\beta(\frak p^m)+\frak p^{m+1}}. 
\tag 6.5.3 
$$ 
We apply the character $\mu_M*\beta$ to the relation $y=z{+}h$. Since $\mu_M(\alpha y) = 1$ (because $d>m{-}w$), we  get 
$$ 
\align 
\mu_M(b y) = \mu_M(\beta y) &= \mu_M(\beta z)\mu_M(\beta h) \\ 
&= \mu_M(\beta z) \mu_M(\alpha h) \\ 
&= \mu_M(\beta z)\mu_M(-\zeta by),  
\endalign 
$$ 
whence $\mu_M((1{+}\zeta)by) = \mu_M(\beta z)$. Our hypothesis $d\not\equiv m \pmod p$ implies that $\zeta \not\equiv -1 \pmod p$ so, for some choice of $z$, we get $\mu_M(\beta z) \neq 1$. In combination with the previous argument, this proves that $w'=m{-}d$ and the first assertion (6.5.2) of the theorem. 
\par  
Let $\theta \in \scr C(\frak a,\alpha) = \scr C(\frak a,\beta)$ and suppose $l = l_E(\theta) <d = m{-}w'$. We calculate the $E'$-level $l_{E'}(\theta)$. If $y\in E$, $\ups_E(y) = 1{+}d$, we write $y=z{+}h$ as above, with $z\in E'$ of valuation $1{+}d$ and $h\in A_\alpha(\frak p^m){+}\frak p^{m+1}$. This gives $1=\theta(1{+}y) = \theta(1{+}z)\theta(1{+}h) = \theta(1{+}z)$. Thus $l_{E'}(\theta)\le d$. Now take $y\in E$ of valuation $d$ and write $y=z{+}h$, where $\ups_{E'}(z) = d$ and $h\in \frak p^m$. Indeed, we may take $h \equiv -\zeta by\alpha^{-1} \pmod{A_\alpha(\frak p^m){+}\frak p^{m+1}}$ as before. This gives 
$$ 
1 = \theta(1{+}y) = \theta(1{+}z)\mu_M(\alpha h), 
$$ 
and 
$$ 
\mu_M(\alpha h) = \mu_M(-\zeta by) = \mu_E(-\zeta ys_\alpha(b)). 
$$ 
Suppose $d\not\equiv 0\pmod p$. Thus $\zeta \not\equiv 0 \pmod{\frak p_F}$ and we may choose $y\in \frak p_E^d$ so that $\mu_E(-\zeta ys_\alpha(b)) \neq 1$.  Thus $\theta(1{+}z) \neq 1$, whence $l_{E'}(\theta) = d$ as required for (1)(a). If $d\equiv 0 \pmod p$, then $\zeta \equiv 0 \pmod{\frak p_F}$ and $\theta(1{+}z) = 1$. Thus $l_{E'}(\theta) < d$, as required for (1)(b). 
\par 
Part (2) follows from a similar, but easier, argument. 
\par 
Part (3) is given by a counting argument as follows. Let $q$ be the cardinality of the residue field $\frak o_F/\frak p_F$. For an integer $k\le [m/2]$, let $\scr C(\alpha;\,{\le} k)$ be the set of $\theta\in \scr C(\frak a,\alpha)$ such that $l_{F[\alpha]}(\theta) \le k$. We use the obvious variations. Note that $\scr C(\alpha;\,{\le} k)$ has exactly $q^k$ elements while $\scr C(\alpha;\,{>} k)$ has $q^{[m/2]-k}$ elements. 
\par 
Part (2) gives $\scr C(\alpha;\,{>}d) \i \scr C(\beta;\,{>}d)$, hence $\scr C(\alpha;\,{>}d) = \scr C(\beta;\,{>}d)$ and also $\scr C(\alpha;\,{\le} d) = \scr C(\beta;\,{\le} d)$. Assertions (3)(a) and (3)(b) now follow from (1)(a) and (1)(b) respectively. \qed 
\enddemo 
We refine the final step of the argument, retaining the notation of the theorem.  
\proclaim{Corollary 1} 
\roster 
\item 
There is a unique character $\xi$ of\/ $U^d_{E'}/U^{1+d}_{E'}$ with the following property: if $\theta\in \scr C(\frak a,\alpha)$ has $l_{E'}(\theta) = d$, then $l_E(\theta) < d$ if and only if $\theta \Mid U^d_{E'} = \xi$. 
\item 
Let $\theta_0$ be the unique element of $\scr C(\frak a,\alpha)$ such that $l_E(\theta_0) = 0$. It satisfies $l_{E'}(\theta_0) \le d$ and the character $\xi$ of \rom{(1)} is given as $\xi = \theta_0 \Mid U^d_{E'}$. 
\item 
The character $\xi$ is trivial if and only if $d\equiv 0 \pmod p$. 
\endroster 
\endproclaim 
\demo{Proof}
Let $\theta_0\in \scr C(\frak a,\alpha)$ have $l_E(\theta_0) = 0$ and endo-class $\vT_0$. Let $\xi$ be the restriction of $\theta_0$ to $U_{E'}^d$. By assertion (1) of the theorem, this character $\xi$ is trivial if and only if $d\equiv 0 \pmod p$. Let $\theta'\in \scr C(\frak a,\alpha)$ have endo-class $\vT'$. If $\Bbb A$ is the canonical ultrametric on $\Scr E(F)$, then $l_E(\theta') < d$ if and only if $\Bbb A(\vT_0,\vT') < p^{-r}d$. This condition is also equivalent to $\theta'$ agreeing with $\theta_0$ on $U^d_{E'}$. \qed 
\enddemo 
\proclaim{Corollary 2} 
Let $\theta\in \scr C(\frak a,\alpha)$ satisfy $l_E(\theta) = d$. In the theorem, one may choose $\beta$ so that $l_{E'}(\theta) = d$. 
\endproclaim 
\demo{Proof} 
If $d \equiv 0 \pmod p$, there is nothing more to do, so we assume the contrary. Let $y\in \frak p_E^d$, and write $y = z{+}h$, for $z\in \frak p_{E'}^d$ and $h\in \frak p^m$, satisfying (6.5.3). Thus $\theta(1{+}y) = \theta(1{+}z)\,\mu_M(\alpha h) = \theta(1{+}z)\,\mu_M(-\zeta by)$. The function $1{+}y \mapsto \mu_M(-\zeta by)$ represents a non-trivial character of $U^d_E/U^{1+d}_E$. We may choose $b$ at the beginning so that $\mu_M(-\zeta by) \neq \theta(1{+}y)$, for some $y\in \frak p_E^d$. This gives $\theta(1{+}z) \neq 1$ and $l_{E'}(\theta) = d$, as required. \qed 
\enddemo 
\head\Rm 
7. The Herbrand function 
\endhead 
We continue with a simple stratum $[\frak a,m,0,\alpha]$ as in the previous sections. We recall that $\|\scr C(\frak a,\alpha)\|$ is the set of endo-classes of simple characters $\theta\in \scr C(\frak a,\alpha)$ and that the canonical map $\scr C(\frak a,\alpha) \to \|\scr C(\frak a,\alpha)\|$ is a bijection (2.3 Remark (2)). 
\par 
In this section, we state and prove the main results concerning the Herbrand function $\Psi_\vT$, $\vT \in \|\scr C(\frak a,\alpha)\|$. In 7.2 Theorem 1 and the supplementary 7.5 Proposition, we describe these functions in coherent families, rather along the lines of 5.3 Theorem but exploiting the flexibility gained in section 6. In 7.2 Theorem 2, we take a rather different approach. We fix $\alpha$ and specify, via an explicit formula, a non-empty subset $\scr C^\star(\frak a,\alpha)$ of $\scr C(\frak a,\alpha)$, the elements of which are the simple characters that {\it conform to\/} $\alpha$. If $\vT$ is the endo-class of $\theta\in \scr C^\star(\frak a,\alpha)$, we show $\Psi_\vT = \biP{F[\alpha]}F\vs$. All Herbrand functions $\Psi_\vT$, $\vT\in \ewc F$, are captured this way. The description given by Theorem 2 has particularly good properties with respect to the Langlands correspondence (section 10 below), but its proof relies on Theorem 1. 
\subhead 
7.1 
\endsubhead 
We introduce a new concept. 
\definition{Definition} 
Let $\vT \in \ewc F$ have degree $p^r$. Let $\theta\in \scr C(\frak a,\alpha)$ be a realization of $\vT$, on a simple stratum $[\frak a,m,0,\alpha]$ in $M = \M{p^r}F$. Let $E = F[\alpha]$ and $l = l_E(\theta)$ (5.2).  Say that $\alpha$ is {\it $\theta$-conformal\/} if 
$$ 
\theta(1{+}y) = \mu_M(\alpha y), \quad y\in \frak p_E^{1+[l/2]}. 
$$ 
Say $\alpha$ is {\it weakly $\theta$-conformal\/} if 
$$ 
\theta(1{+}y) = \mu_M(\alpha y), \quad y\in \frak p_E^l. 
$$ 
\enddefinition 
In this situation, we might equally say that $\theta$ is $\alpha$-conformal. Let $\scr C^\star(\frak a,\alpha)$ be the set of $\alpha$-conformal $\theta\in \scr C(\frak a,\alpha)$. Surely $\scr C^\star(\frak a,\alpha)$ is not empty. 
\par 
Let $\|\scr C^\star(\frak a,\alpha)\|$ be the set of endo-classes of elements of $\scr C^\star(\frak a,\alpha)$. The canonical map $\scr C^\star(\frak a,\alpha) \to \|\scr C^\star(\frak a,\alpha)\|$ is a bijection. 
\proclaim{Proposition} 
Let $\vT \in \ewc F$ be of degree $p^r$. The endo-class $\vT$ has a realization $\theta\in \scr C(\frak a,\alpha)$, on a simple stratum $[\frak a,m,0,\alpha]$ in $M = \M{p^r}F$, such that $\alpha$ is $\theta$-conformal. That is, $\vT \in \|\scr C^\star(\frak a,\alpha)\|$. 
\endproclaim 
\demo{Proof} 
Let $\theta\in \scr C(\frak a,\alpha)$ be a realization of $\vT$ and let $\frak p = \roman{rad}\,\frak a$. Let $\nu_\theta(\alpha)$ be the least integer $\nu$ for $\theta(1{+}y) = \mu_M(\alpha y)$, $y\in \frak p_E^{1+\nu}$. Certainly $\nu \le [m/2]$ (2.3.1). Write $E = F[\alpha]$ and $d_\alpha = m{-}w_{E/F}$. We have $\nu\ge [d_\alpha/2]$ since, otherwise, the function $\mu_M*\alpha$ does not represent a character of $U^{1+\nu}_E$. If $\nu = [d_\alpha/2]$, there is nothing more to do. 
\proclaim{Lemma} 
Set $\nu = \nu_\theta(\alpha)$, and assume that $\nu > [d_\alpha/2]$. There exists $\beta \in \roman P(\frak a,\alpha)$ such that $\beta\equiv \alpha\pmod{\frak p^{-\nu}}$ and $\nu_\theta(\beta) \le \nu{-}1$. This condition determines the stratum $[\frak a,m,\nu{-}1,\beta]$ uniquely, up to formal intertwining. 
\endproclaim 
\demo{Proof} 
Recall that $\nu\le [m/2]$. By hypothesis, the function 
$$ 
\xi(1{+}x) = \theta(1{+}x)\,\mu_M(-\alpha x), \quad x\in \frak p_E^\nu, 
$$ 
represents a non-trivial character of $U^\nu_E$, trivial on $U^{1+\nu}_E$. Consequently, there exists $z\in \frak p^{-\nu}$ so that 
$$ 
\xi(1{+}x) = \mu_M(zx),\quad x\in \frak p_E^\nu. 
\tag 7.1.1 
$$ 
Choose a tame corestriction $s_{E/F}:M \to E$ and let $\mu_E$ be the character of $E$ for which $\mu_E\circ s_{E/F} = \mu_M$. Thus (7.1.1) reads $\xi(1{+}x) = \mu_E(s_{E/F}(z)x)$, for all $x$ as before. As $\xi$ defines a non-trivial character of $U^\nu_E/U^{1+\nu}_E$, so $\ups_E(s_{E/F}(z)) = -\nu$ and $s_{E/F}(z)$ is uniquely determined, by $\theta$, modulo $\frak p_E^{1-\nu}$. We invoke \cite{15} (2.2.3): the stratum $[\frak a,m,\nu{-}1,\alpha{+}z]$ is simple (whence $\beta\in \roman P(\frak a,\alpha)$) and uniquely determined up to formal intertwining \cite{15} (2.2.1). 
\par 
Set $L = F[\beta]$. We show that 
$$ 
\theta(1{+}x) = \mu_M(\beta x), \quad x\in \frak p_L^\nu. 
\tag 7.1.2   
$$ 
This will imply $\nu_\vT(\beta) \le \nu{-}1$, as required to complete the proof of the lemma. 
\par 
Since $x\in L = F[\beta]$, there is a polynomial $f(T) \in F[T]$, of degree at most $p^r{-}1$, such that $x = f(\beta)$. Write 
$$ 
f(T) = a_0+a_1T+\dots\dots + a_{p^r-1}T^{p^r-1}. 
$$  
The $L$-valuations of the terms $a_i\beta^i$, $0\le i\le p^r{-}1$, are distinct. The condition $\ups_L(x) = \nu$ translates as $\nu\le p^r\ups_F(a_i){-}mi$ for all $i$, with equality for exactly one value of $i$. So, if we put $y = f(\alpha)$, we get $\ups_E(y) = \nu$. Consider the element 
$$ 
t = x{-}y = f(\beta)-f(\alpha) = \sum_{1\le i<p^r} a_i\big((\alpha{+}z)^i{-}\alpha^i\big). 
$$ 
Expand $\big((\alpha{+}z)^i{-}\alpha^i\big)$. Any fractional $\frak a$-ideal $\frak p^k$, $k\in \Bbb Z$, is stable under conjugation by $\alpha$, so every term in the expansion of $(\alpha{+}z)^i{-}\alpha^i$ lies in $\alpha^{i-1}z\frak a = \frak p^{(1{-}i)m-\nu}$. Since $p^r\ups_F(a_i) \ge mi{+}\nu$, the term $a_i\big((\alpha{+}z)^i{-}\alpha^i\big)$ lies in $\frak p^m$, whence $t = f(\beta){-}f(\alpha) = x{-}y \in \frak p^m$. 
\par 
With this element $t$, and setting 
$$ 
u = (1{+}t)^{-1}(1{+}y)^{-1}yt, 
$$ 
we have 
$$ 
1{+}x = (1{+}y)(1{+}t)(1{-}u). 
$$ 
We use this expression to evaluate $\theta(1{+}x)$. Our choice of $z$ gives $\theta(1{+}y) = \mu_M(\beta y)$ and, since $t\in \frak p^m$, we have $\theta(1{+}t) = \mu_M(\alpha t)$. As $yt\in \frak p^{m+1}$, so $\theta(1{-}u) = 1$. Therefore 
$$ 
\theta(1{+}x) = \theta(1{+}y)\theta(1{+}t) = \mu_M(\beta y)\mu_M(\alpha t). 
$$ 
On the other hand, $zt\in \frak p^{m-\nu}$ and $m{-}\nu = (m{-}2\nu){+}\nu \ge 1$, whence $\mu_M(zt) = 1$. Altogether, 
$$ 
\mu_M(\beta x)  = \mu_M(\beta y)\,\mu_M(\alpha t)\,\mu_M(zt) = \theta(1{+}y)\,\theta(1{+}t) = \theta(1{+}x), 
$$ 
as required for (7.1.2). This completes the proof of the lemma. \qed 
\enddemo 
The proposition now follows. \qed 
\enddemo 
Note that, while the proposition is an existence statement, the proof is constructive. 
\subhead 
7.2 
\endsubhead 
To state our first result, it is convenient to have a looser concept reflecting the structure of 5.3 Theorem. We consider a datum $(E/F,m)$ in which $E/F$ is a totally ramified field extension of degree $p^r$, $r\ge 1$, and $m$ is a positive integer not divisible by $p$. 
\definition{Definition} 
Say that $(E/F,m)$ is {\it standard\/} if at least one of the following conditions holds: 
\roster 
\item"\rm (a)" 
$m> 2w_{E/F}$; 
\item"\rm (b)" 
$m\le w_{E/F}$; 
\item"\rm (c)" 
$m\le 2w_{E/F}$ and $w_{E/F} \equiv 0 \pmod p$. 
\endroster 
\enddefinition 
Case (c) can only arise when $F$ has characteristic zero (1.8). Remark that, in case (b), the function $\biP EF{m/p^r}$ has an odd number of jumps (4.2 Remark). In case (c), we actually have $m<2w_{E/F}$, since $m$ is not divisible by $p$. We can always reduce to one of these cases, as follows. 
\proclaim{Lemma}
Let $\vT \in \ewc F$ have degree $p^r$. There is a simple stratum $[\frak a,m,0,\alpha]$ in $\M{p^r}F$ such that 
\roster 
\item 
$\scr C(\frak a,\alpha)$ contains a character $\theta$ of endo-class $\vT$, and 
\item 
the datum $(F[\alpha]/F,m)$ is standard. 
\endroster 
\endproclaim 
\demo{Proof} 
Choose a simple stratum $[\frak a,m,0,\beta]$ in $\M{p^r}F$ such that $\vT \in \|\scr C(\frak a,\alpha)\|$. If $m>2w_{F[\beta]/F}$, then $(F[\beta]/F,m)$ is standard. Otherwise, the lemma follows from 6.3 Corollary. \qed 
\enddemo 
We state our main results. 
\proclaim{Theorem 1} 
Let $\vT \in \ewc F$ have degree $p^r$. Let $\theta\in \scr C(\frak a,\alpha)$ be a realization of $\vT$ on a simple stratum $[\frak a,m,0,\alpha]$ in $\M{p^r}F$ for which the datum $(F[\alpha]/F,m)$ is standard. Write $E = F[\alpha]$, $l = l_E(\theta)$ and $\vs = m/p^r = \vs_\vT$. For any such realization, the following hold. 
\roster 
\item 
If $l\le \roman{max}\{0,m{-}w_{E/F}\}$, then $\Psi_\vT(x) = \upr{2\,}\Psi_{(E/F,\vs)}(x)$, $0\le x\le \vs$. 
\item 
If $l>\roman{max}\{0,m{-}w_{E/F}\}$ and $l \not\equiv m\pmod p$, then 
$$ 
\Psi_\vT(x) = \roman{max}\,\big\{\upr{2\,}\Psi_{(E/F,\vs)}(x), x-p^{-r}(m{-}l)\big\}, \quad  0\le x\le \vs.  
\tag 7.2.1 
$$ 
\item 
In part \rom{(2),} the class $\vT$ admits a parameter field $E'/F$ as follows: 
\itemitem{\rm (i)} $E' = F[\beta]$, where $\beta\in \roman P(\frak a,\alpha)$ and $\beta\equiv \alpha\pmod{\frak p^{-l}}$; 
\itemitem{\rm (ii)} $w_{E'/F} = m{-}l$ and $l_{E'}(\theta) = l$. 
\item"" 
For any such $\beta$, $\Psi_\vT(x) = \upr{2\,}\Psi_{(E'/F,\vs)}(x)$, $0\le x\le \vs$. 
\endroster 
\endproclaim  
That $\vT$ has a realization of the required form is 7.2 Lemma. We shall see in the course of the proof that (7.2.1) also holds in the situation of part (1), but says nothing new there. We comment in 7.5 below on the restrictive hypothesis in part (2) of the theorem. 
\proclaim{Theorem 2} 
Let $\vT \in \ewc F$ have degree $p^r$. Let $[\frak a,m,0,\alpha]$ be a simple stratum in $M = \M{p^r} F$ such that $\vT$ has a realization $\theta \in \scr C^\star(\frak a,\alpha)$. For any such realization, $l_{F[\alpha]}(\theta) = \roman{max}\{0,m{-}w_{F[\alpha]/F}\}$ and 
$$ 
\Psi_\vT(x) = \upr{2\,}\Psi_{(F[\alpha]/F,\vs_\vT)}(x), \quad 0\le x\le \vs_\vT. 
\tag 7.2.2 
$$ 
\endproclaim 
\remark{Remark} 
The endo-class $\vT$ has a realization of the required form, by 7.1 Proposition. When proving Theorem 2, we show that (7.2.2) holds provided only that $\vT$ has a realization $\theta \in \scr C(\frak a,\alpha)$ such that $\alpha$ is {\it weakly $\theta$-conformal\/} (7.1 Definition). We will not use that version in the rest of the paper. 
\endremark 
Before embarking on the proofs of the theorems, we give a consequence of Theorem 2. 
\proclaim{Corollary} 
Let $E/F$ be a totally ramified field extension of degree $p^r$, and let $m$ be a positive integer not divisible by $p$. There exists $\vT \in \ewc F$, of degree $p^r$, with parameter field $E/F$ and $\vs_\vT = m/p^r$, such that $\Psi_\vT(x) = \biP EF{m/p^r}(x)$, $0 \le x \le m/p^r$. 
\endproclaim 
\demo{Proof}
View $E$ as a subfield of $M = \M{p^r}F$ and take $\alpha\in E$ such that $\ups_E(\alpha) = -m$. There is a unique hereditary $\frak o_F$-order $\frak a$ in $M$ such that $[\frak a,m,0,\alpha]$ is a simple stratum in $M$. By Theorem 2, any $\vT \in \|\scr C^\star(\frak a,\alpha)\|$ has the required property. \qed 
\enddemo 
\subhead 
7.3 
\endsubhead 
We prove 7.2 Theorem 1. In part (1) of the theorem, suppose that $(E/F,m)$ is standard of type (a), resp\. (b), resp\. (c). The assertion is then equivalent to part (2), resp\. (1), resp\. (3) of 5.3 Theorem. 
\par 
The hypothesis in part (2) implies that $m < 2w_{E/F}$, so the standard datum $(E/F,m)$ is of type (b) or (c). To prove part (2), we first use 6.5 Corollary 2 to choose an element $\beta\in \roman P(\frak a,\alpha)$ such that 
$$ 
w_{F[\beta]/F} = m{-}l \quad \text{\it and\/} \quad  l_{F[\beta]}(\theta) = l. 
$$ 
Consequently, $w_{F[\beta]/F} < w_{E/F}$ and $w_{F[\beta]/F} \not\equiv 0 \pmod p$. 
\par 
Let $\theta_0$ be the unique element of $\scr C(\frak a,\alpha)$ with $l_E(\theta_0) = 0$ and let $\vT_0$ be the endo-class of $\theta_0$. The hypotheses of part (1) of the theorem apply to $\theta_0$ as an element of $\scr C(\frak a,\alpha)$, so 
$$ 
\Psi_{\vT_0}(x) = \upr{2\,} \Psi_{(E/F,\vs)}(x) , \quad 0\le x\le\vs. 
\tag 7.3.1 
$$ 
We compare $\theta$ and $\theta_0$ from the standpoint of the element $\beta$. From 6.5 Theorem 1(1), we have $l_{F[\beta]}(\theta_0) \le l$, with equality if and only if $l \not\equiv 0 \pmod p$. 
\par  
If the function $\upr{2\,}\Psi_{(F[\beta]/F,\vs)}$ has an odd number of jumps, then 
$$ 
\biP EF\vs = \Psi_{\vT_0} = \upr{2\,}\Psi_{(F[\beta]/F,\vs)} = \Psi_\vT, 
$$ 
by (7.3.1) and 5.3 Theorem (1) applied to $\vT_0$ and to $\vT$. Moreover, 
$$ 
x{-}p^{-r}(m{-}l) = x{-}p^{-r}w_{F[\beta]/F} \le \upr{2\,}\Psi_{(F[\beta]/F,\vs)}(x), \quad 0\le x\le \vs, 
$$ 
so we are done in this case. 
\par 
Assume therefore that $\upr{2\,}\Psi_{(F[\beta]/F,\vs)}$ has an even number of jumps. Let $I$ be the set of points $x$ such that $\biP{F[\beta]}F\vs'(x) = 1$. Since $\biP{F[\beta]}F\vs$ has an even number of jumps, the set $I$ is a non-empty open interval and (4.2 Proposition) 
$$ 
\biP{F[\beta]}F\vs(x) = x-p^{-r}w_{F[\beta]/F}, \quad x\in I. 
\tag 7.3.2 
$$ 
By 5.6 Corollary 1, the functions $\Psi_{\vT_0}$, $\biP {F[\beta]}F\vs$ agree outside $I$. By 5.6 Corollary 2, the only possibilities are that $\Psi_{\vT_0}(x) = \biP{F[\beta]}F\vs(x)$ for all $x\in I$, or else $\Psi_{\vT_0}(x) < \biP{F[\beta]}F\vs(x)$ for all $x\in I$. The first alternative is untenable: if $\Psi'_{\vT_0}=1$ on an interval $I'$, then (by (7.3.1)) $\Psi_{\vT_0}(x) = x{-}p^{-r}w_{E/F}$, $x\in I'$. But, if $\Psi_{\vT_0}(x)$ equalled $\biP{F[\beta]}F\vs(x)$ on $I$, we would have $\Psi_{\vT_0}(x) = x{-}p^{-r}w_{F[\beta]/F}$ there. Since $w_{E/F} > w_{F[\beta]/F}$, this is impossible. Therefore 
$$ 
\biP{F[\beta]}F\vs(x) > \Psi_{\vT_0}(x), \quad x\in I, 
\tag 7.3.3 
$$ 
and 
$$ 
\biP{F[\beta]}F\vs(x) = \roman{max}\,\big\{\Psi_{\vT_0}(x),x-p^{-r}w_{F[\beta]/F}\big\}, \quad 0\le x\le\vs. 
\tag 7.3.4 
$$ 
In terms of the ultrametric $\Bbb A$ on $\Scr E(F)$, we have $\Bbb A(\vT,\vT_0) = l/p^r > 0$. It follows that the characters $\theta$, $\theta_0$ do not agree on $U^l_{F[\beta]}$. Theorem 5.3 (4) now implies $\Psi_\vT = \biP{F[\beta]}F\vs$ and Part (2)  follows from (7.3.1) and (7.3.4). 
\par 
Part (3) holds relative to the same choice of $\beta$, so we have completed the proof of 7.2 Theorem 1. \qed 
\remark{Remark} 
The argument following (7.3.4) shows that the character $\phi$ of 5.3 Theorem (4), relative to $l$ and $\beta$, is $\theta_0\Mid U^l_{F[\beta]}$. It is trivial if and only if $l\equiv 0\pmod p$. 
\endremark 
\subhead 
7.4 
\endsubhead 
We prove 7.2 Theorem 2. Let $\theta\in \scr C(\frak a,\alpha)$ be a realization of $\vT$ for which $\alpha$ is {\it weakly $\theta$-conformal\/} and set $E = F[\alpha]$. Thus $l = l_E(\theta) = m{-}w_{E/F}$ or $0$. 
\par 
If either $m> 2w_{E/F}$ or $w_{E/F} \equiv 0 \pmod p$, the desired relation $\Psi_\vT = \biP EF\vs$ is given by 5.3 Theorem (2) or (3) respectively. We therefore assume assume that $m\le 2w_{E/F}$ and that $w_{E/F} \not\equiv 0 \pmod p$. In particular, $l \not\equiv m\pmod p$. If $\biP EF{m/p^r}$ has an odd number of jumps, then $\Psi_\vT = \biP EF\vs$ by 5.3 Theorem (1).  
\par 
We therefore assume that $\biP EF\vs$ has an even number of jumps (whence $(E/F,m)$ is not standard). Let $I$ be the open interval on which $\biP EF\vs'(x) = 1$, $0<x<\vs$. For $0\le x\le \vs$, we have (5.3 Theorem (4))
$$ 
\align 
\Psi_\vT(x) &= \biP EF\vs(x), \quad x\notin I, \\ 
\Psi_\vT(x) &\le \biP EF\vs(x) = x-p^{-r}(m{-}l), \quad x\in I. 
\endalign 
$$ 
Consequently, $\Psi_\vT(x) = \biP EF\vs(x)$ at any point where these functions are smooth and have derivative other than $1$. 
\par 
We use 6.3 Corollary to construct from $\alpha$ a standard datum $(F[\beta]/F,m)$: this will be of type (b) or (c) in the scheme of 7.2 Definition. Since $(E/F,m)$ is not standard, $w_{F[\beta]/F} > w_{E/F}$. By 6.4 Proposition, $l_{F[\beta]}(\theta) = l$. By 7.2 Theorem 1(2), 
$$ 
\Psi_\vT(x) = \roman{max}\,\big\{\biP{F[\beta]}F\vs(x), x-p^{-r}(m{-}l)\big\}, \quad 0\le x\le \vs. 
$$ 
So, if $\Psi_\vT$ is smooth at $x$ and $\Psi_\vT'(x) \neq 1$, then $\Psi_\vT(x) = \biP{F[\beta]}F\vs(x) = \biP EF\vs(x)$. Suppose, on the other hand, that $\Psi_\vT'(x) = 1$. If $\biP{F[\beta]}F\vs'(x) = 1$, then 
$$ 
\biP{F[\beta]}F\vs(x) = x - p^{-r}w_{F[\beta]/F} < x-p^{-r}(m{-}l). 
$$ 
Therefore $\Psi_\vT(x) = x-p^{-r}(m{-}l) = \biP EF\vs(x)$ at such points. Altogether, $\Psi_\vT(x) = \biP EF\vs(x)$ for $0\le x \le \vs$. We have proved 7.2 Theorem 2. \qed 
\subhead 
7.5 
\endsubhead 
Now that Theorem 2 has been proved, Theorem 1 has no further direct role in the paper. However, Theorem 2 gives no idea of how $\Psi_\vT$ varies as $\vT$ ranges over $\|\scr C(\frak a,\alpha)\|$ while Theorem 1 does just that, modulo some limitations in part (2). For the sake of tidiness, we show that all Herbrand functions $\Psi_\vT$, $\vT\in \|\scr C(\frak a,\alpha)\|$, are captured within the scheme of Theorem 1. 
\proclaim{Proposition} 
Suppose that $m \le 2w_{F[\alpha]/F}$. Let $\vT\in \|\scr C(\frak a,\alpha)\|$. There exists $\beta\in \roman P(\frak a,\alpha)$, say $F[\beta] = E$, with the following properties.  
\roster 
\item The datum $(E/F,m)$ is standard and 
\item either 
\itemitem{\rm (a)} $l_E(\theta) \le \roman{max}\,\{0,m{-}w_{E/F}\}$ or 
\itemitem{\rm (b)} $l_E(\theta) \not\equiv m \pmod p$.  
\endroster 
\endproclaim  
\demo{Proof} 
We first choose $\beta$ so that $\vT$ is the endo-class of some $\theta\in \scr C^\star(\frak a,\beta)$, as we may by 7.1 Proposition. Writing $E = F[\beta]$, suppose $w_{F[\beta]/F} \equiv 0 \pmod p$. Thus $(E/F,m)$ is standard and  $l_E(\theta) = \roman{max}\,\{0,m{-}w_{E/F}\}$, so option (a) applies. 
\par 
Suppose then that $w_{E/F} \not\equiv 0\pmod p$. Thus $l_E(\theta) = \roman{max}\,\{0,m{-}w_{E/F}\}$, so $l_E(\theta) \not \equiv m \pmod p$. If $(E/F,m)$ is standard, there is nothing to do, so suppose otherwise. We use 6.2 Corollary to find $\gamma \in \roman P(\frak a,\beta)$ such that, if $L = F[\gamma]$, then either $w_{L/F} \ge m$ or $w_{L/F} \equiv 0 \pmod p$ and $w_{L/F} > w_{E/F}$. In all cases, $(L/F,m)$ is standard. By 6.3 Proposition, $l_L(\theta) = l_E(\theta) \not\equiv m \pmod p$, so option (b) applies. \qed 
\enddemo 
Recall that, in the proposition, there is nothing to say when $m>2w_{E/F}$ (5.3 Theorem (2)). Otherwise, $\Psi_\vT$ is given by 7.2 Theorem 1. 
\remark{Remark} 
The theorems of 7.2 and the proposition above leave open the following question. What are the functions $\biP{F[\beta]}F\vs$, where $\beta$ ranges over elements of $\roman P(\frak a,\alpha)$ subject to the condition that the datum $(F[\beta]/F,m)$ is standard? 
\endremark 
\head\Rm 
8. Representations with a single jump 
\endhead 
We consider here representations $\sigma\in \wwr F$ for which the decomposition function $\vS_\sigma$ of (2.2.2) has {\it a unique jump:\/} these play a central role in what follows. 
\subhead 
8.1 
\endsubhead 
For the moment, let $G$ be a finite $p$-group with centre $Z \neq G$. Say that $G$ is {\it H-cyclic\/} if $Z$ is cyclic and $G/Z$ is elementary abelian. Equivalently, $G$ is an extra-special $p$-group of class $2$. We introduce this new terminology to avoid ambiguous usage that has accumulated here. In particular, we do not need to specify $G$ from among the various possibilities listed in, for instance, \cite{19} p.~203 {\it et seq.} The material of this sub-section is generally familiar, but we choose to give a complete, albeit brief, account. 
\par 
If $G$ is H-cyclic, the commutator group $[G,G]$ is the subgroup $Z_p$ of $Z$ of order $p$. We may view the pairing $G/Z\times G/Z \to Z_p$, induced by the commutator $(x,y) \mapsto [x,y]$, as an alternating form on the $\Bbb F_p$-vector space $G/Z$. If $x,y\in G$, then $[x,y] = 1$ if and only if $x$ centralizes $y$. The alternating form is therefore {\it nondegenerate:\/} if $[x,y] = 1$ for all $y\in G$, then $x\in Z$. 
\par 
We first give a technical result, needed in 8.4. 
\proclaim{Lemma} 
Let $G$ be an $H$-cyclic finite $p$-group with centre $Z$. Let $\alpha$ be an automorphism of $G$ which is trivial on $Z$ and induces the trivial automorphism of $G/Z$. The automorphism $\alpha$ is then inner. 
\endproclaim 
\demo{Proof} 
Consider the map $G\to Z_p$ given by $x\mapsto x^\alpha x^{-1}$. This induces a map $G/Z \to Z_p$ which is a homomorphism: $(xy)^\alpha y^{-1}x^{-1} = x^\alpha x^{-1}y^\alpha y^{-1}$. The nondegeneracy property of the commutator pairing gives a unique $y\in G/Z$ such that $x^\alpha x^{-1} = [y,x]$, for all $x$. This relation says $x^\alpha = yxy^{-1}$, as required. \qed 
\enddemo 
\proclaim{Proposition} 
Let $G$ be an H-cyclic finite $p$-group with centre $Z$, and let $\chi$ be a faithful character of $Z$. 
\roster 
\item There exists a unique irreducible representation $\sigma$ of $G$ such that $\sigma\big|Z$ contains $\chi$. The representation $\sigma$ is faithful of dimension $(G{:}Z)^{\frac12}$ and $\sigma\,\big|\,Z$ is a multiple of $\chi$. 
\item A character $\xi$ of $G$ satisfies $\xi \otimes \sigma \cong \sigma$ if and only if $\xi$ is trivial on $Z$. If $D(\sigma)$ denotes the group of such characters, then 
$$
\check\sigma\otimes \sigma = \sum_{\xi\in D(\sigma)} \xi. 
\tag 8.1.1 
$$ 
\endroster 
\endproclaim 
\demo{Proof} 
Denote by $h$ the alternating form on $G/Z$ induced by the commutator pairing $(x,y) \mapsto \chi[x,y]$, $x,y\in G$. The nondegenerate  alternating $\Bbb F_p$-space $G/Z$ has even dimension $2r$, say. Let $L$ be a {\it Lagrangian subspace\/} of $G/Z$, that is, a subspace on which $h$ is null and is maximal for this property. Thus $L$ has dimension $r$. 
\par 
Let $\wt L$ be the inverse image of $L$ in $G$. As $h$ is null on $L$, the subgroup $\wt L$ of $G$ is abelian and maximal for this property. The character $\chi$ therefore admits extension to a character $\chi_L$ of $\wt L$. Let $y\in G\smallsetminus L$. There exists $x\in L$ such that $[x,y] \neq 1$. This implies that $\chi_L^y(x) \neq \chi_L(x)$, whence $\rho_\chi = \Ind_{\wt L}^G\chi_L$ is irreducible. We form the usual inner product of characters, 
$$ 
1 = \langle\roman{tr}\,\rho_\chi, \roman{tr}\,\rho_\chi\rangle = |G|^{-1}\sum_{g\in G} |\roman{tr}\,\rho_\chi(g)|^2. 
$$ 
As $\roman{tr}\,\rho_\chi(z) = p^r\chi(z)$, for $z\in Z$, it follows that $\roman{tr}\,\rho_\chi(g) = 0$, for all $g\in G\smallsetminus Z$. Therefore $\rho_\chi$ is independent of the choice of $\chi_L$. The function $\roman{tr}\,\rho_\chi$ takes the value $p^r = \dim\rho_\chi$ only at the identity, so $\rho_\chi$ is faithful. 
\par 
Let $\sigma$ be an irreducible representation of $G$ that contains $\chi$. With $L$ as before, the restriction $\sigma \Mid \wt L$ is a sum of characters $\phi$ (since $\wt L$ is abelian), each of which satisfies $\phi\Mid Z = \chi$. However, any such character induces the representation $\rho_\chi$, so $\sigma \cong \rho_\chi$, as asserted. This deals with (1). 
\par 
A character $\xi$ of $G$ such that $\xi\otimes \sigma \cong \sigma$ is surely trivial on $Z$. That is, $\xi$ is the inflation of a character of $G/Z$. The trace calculation ensures that any such character $\xi$ satisfies $\xi\otimes \sigma = \sigma$. Thus $\xi$ occurs in the representation $\check\sigma \otimes \sigma$. The number of such characters $\xi$ is $p^{2r} = \dim \check\sigma\otimes \sigma$, whence follows (8.1.1). \qed 
\enddemo 
\subhead 
8.2 
\endsubhead 
Let $\sigma \in \wwr F$ have dimension $p^r$, and let $\bar\sigma: \scr W_F \to \roman{PGL}_{p^r}(\Bbb C)$ be the projective representation defined by $\sigma$. 
\definition{Definition 1} 
The {\it centric field\/} $Z = Z_\sigma/F$ of $\sigma$ is defined by $\scr W_Z = \roman{Ker}\,\bar\sigma$. The {\it tame centric field\/} $T_\sigma/F$ of $\sigma$ is the maximal tame sub-extension of $Z_\sigma/F$. 
\enddefinition 
Thus $\sigma$ is absolutely wild if and only if $T_\sigma = F$. Observe that, if $K/F$ is a finite tame extension and $\sigma_K = \sigma\Mid \scr W_K \in \wwr K$, then $Z_{\sigma_K} = Z_\sigma K$ and $T_{\sigma_K} = T_\sigma K$.
\par 
Define $D(\sigma)$ to be the group of characters $\chi$ of $\scr W_F$ such that $\chi\otimes\sigma\cong \sigma$. 
\par 
Since $\sigma\in \wwr F$, the restriction $\sigma_0^+ = \sigma\,\big|\,\scr P_F$ is irreducible. Let $D^+(\sigma)$ be the group of characters $\phi$ of $\scr P_F$ such that $\phi\otimes \sigma_0^+ \cong \sigma_0^+$. Since $\sigma_0^+$ factors through a representation of a finite $p$-group, the group $D^+(\sigma)$ is non-trivial. A character $\phi$ of $\scr P_F$ lies in $D^+(\sigma)$ if and only if it is a component of $\check\sigma_0^+\otimes \sigma_0^+$, whence $D^+(\sigma)$ has order at most $p^{2r}$. The group $\scr W_F$ acts on $D^+(\sigma)$ in a natural way, with $\scr P_F$ acting trivially. 
\par 
If $K/F$ is a finite tame extension, then $\scr P_K = \scr P_F$. We may identify $(\sigma_K)_0^+$ with $\sigma^+_0$ and $D^+(\sigma_K)$ with $D^+(\sigma)$. 
\proclaim{Lemma} 
Let $\sigma\in \wwr F$. 
\roster 
\item 
If $K/F$ is a finite, tamely ramified field extension, the restriction map $D(\sigma_K) \to D^+(\sigma)$ is an isomorphism of $D(\sigma_K)$ with the group of $\scr W_K$-fixed points in $D^+(\sigma)$. 
\item There is a unique minimal tame extension $T_I(\sigma)/F$ such that the map $D(\sigma_{T_I(\sigma)}) \to D^+(\sigma)$ is an isomorphism. 
\item 
The field extension $T_I(\sigma)/F$ is Galois and contained in $T_\sigma$. 
\endroster 
\endproclaim
\demo{Proof} 
The lemma summarizes the discussion in \cite{13} 8.2. \qed 
\enddemo 
We refer to $T_I(\sigma)$ as the {\it imprimitivity field\/} of $\sigma$. 
\definition{Definition 2} 
Let $\sigma\in \wwr F$. Say that $\sigma$ is {\it H-cyclic\/} if the finite $p$-group $\sigma(\scr P_F)$ is H-cyclic. 
\enddefinition 
\proclaim{Proposition} 
If $\sigma\in \wwr F$ is H-cyclic then $T_I(\sigma) = T_\sigma$. 
\endproclaim 
\demo{Proof} 
Let $Z_\sigma/F$ be the centric field of $\sigma$. Since $T_\sigma$ contains $T_I(\sigma)$, nothing is changed if we extend the base field to $T_I(\sigma)$ and assume $T_I(\sigma) = F$. According to the lemma, the group $D(\sigma)$ is then isomorphic to $D^+(\sigma)$ and so has order $p^{2r}$, where $p^r=\dim\sigma$. The non-trivial characters in $D(\sigma)$ are wildly ramified of order $p$. The sum $\sum_{\phi\in D(\sigma)}\phi$ is a sub-representation of $\sigma\otimes\check\sigma$, of the same dimension, so $\check\sigma \otimes \sigma = \sum_{\phi\in D(\sigma)} \phi$. We show that $\check\sigma\otimes \sigma$ provides a {\it faithful\/} representation of $\Gal{Z_\sigma}F$. Let $\sigma$ act on the vector space $V$. So, if $x\in \roman{Ker}\,\check\sigma\otimes\sigma$, the operator $1 = \check\sigma(x)\otimes \sigma(x)  \in \End{\Bbb C}{\check V\otimes V}$ is, in particular, a non-zero scalar. Each of the operators $\sigma(x)\in  \End{\Bbb C}V$, $\check\sigma(x) \in \End{\Bbb C}{\check V}$, is therefore scalar. In particular, $x\in \roman{Ker}\,\bar\sigma = \scr W_{Z_\sigma}$, as asserted. 
\par 
Define $K/F$ by $\scr W_K = \bigcap_{\phi\in D(\sigma)} \roman{Ker}\,\phi$. The extension $K/F$ is totally wildly ramified, and elementary abelian of degree $p^{2r}$. By definition, every $\phi\in D(\sigma)$ is trivial on $\Gal {Z_\sigma}K$, whence $K = Z_\sigma$ and so $T_\sigma = F$. \qed 
\enddemo 
\remark{Remark} 
Following the proposition, it is natural to ask whether there exists a representation $\sigma\in \wwr F$ for which the tame centric field and the imprimitivity field are distinct. We produce an example of such a representation $\sigma$ in 9.7 below. 
\endremark 
The following device is not central to our current concerns, but we include it here for its utility in constructing examples (as in 8.4 below). 
\proclaim{Example} 
Let $\sigma, \sigma'\in \wwr F$ be H-cyclic. The following are equivalent: 
\roster 
\item $D^+(\sigma)\cap D^+(\sigma') = \{1\}$; 
\item $\sigma\otimes\sigma'$ is irreducible and totally wild. 
\endroster 
When these conditions hold, the representation $\sigma\otimes\sigma'$ is H-cyclic. 
\endproclaim 
\demo{Proof} 
If $\tau$ is a smooth, finite-dimensional representation of $\scr P_F$, then $\tau$ is irreducible if and only if the space $\Hom {\scr P_F}1{\tau\otimes \check\tau}$ has dimension one. Here, $\sigma\otimes \check \sigma\Mid \scr P_F = \sum_{\phi\in D^+(\sigma)} \phi$, and similarly for $\sigma'$. Therefore 
$$
(\sigma\otimes\sigma') \otimes (\check\sigma\otimes\check\sigma')\Mid\scr P_F = \sum\Sb \phi\in D^+(\sigma), \\ \phi'\in D^+(\sigma') \endSb \phi\phi'. 
$$ 
The trivial character occurs exactly once in the sum if and only if $D^+(\sigma)\cap D^+(\sigma') = \{1\}$, so (1) is equivalent to $\sigma\otimes \sigma'$ being irreducible on $\scr P_F$: this is the same as (2).  
\par 
Abbreviate $\tau = \sigma\otimes\sigma'$, and assume $\tau$ to be irreducible. Let $C$ and $C'$ be respectively the centres of $\sigma(\scr P_F)$ and $\sigma'(\scr P_F)$. For $x\in \scr P_F$, the operator $\tau(x)^p = \sigma(x)^p \otimes \sigma'(x)^p$ is scalar and lies in $CC' = \{z\otimes z':z\in C, z'\in C'\}$. In particular, $CC'$ consists of scalars and is central in $\tau(\scr P_F)$. Thus $\tau(\scr P_F)$ is of exponent $p$ modulo its centre. Since $\tau$ is irreducible on $\scr P_F$, this centre is cyclic. \qed 
\enddemo 
\subhead 
8.3 
\endsubhead 
Let $\chi$ be a character of $\scr P_F$. Define the {\it $F$-slope $\roman{sl}_F(\chi)$ of\/} $\chi$ by 
$$ 
\roman{sl}_F(\chi) = \roman{inf}\,\{x>0: \scr R_F(x) \subset \roman{Ker}\,\chi\}. 
\tag 8.3.1 
$$ 
If $\chi$ extends to a character $\tilde\chi$ of $\scr W_F$, then $\roman{sl}_F(\chi) = \sw(\tilde\chi) = \vs_{\tilde\chi}$. 
\par 
Let $\sigma\in \wwr F$ be H-cyclic, with $\dim\sigma >1$. Say that $\sigma$ is {\it H-singular\/} if there exists $a>0$ such that $\roman{sl}_F(\chi) = a$, for all non-trivial $\chi\in D^+(\sigma)$. 
\proclaim{Proposition} 
Let $\sigma \in \wwr F$ be H-singular and let $a = \roman{sl}_F(\chi)$, for $\chi\in D^+(\sigma)$, $\chi\neq 1$. The function $\vS_\sigma$ has a unique jump, lying at the point $a$. 
\endproclaim 
\demo{Proof} 
This is immediate, on applying (2.2.2) and (8.1.1) to $\sigma$. \qed 
\enddemo 
\subhead 
8.4 
\endsubhead 
The converse of 8.3 Proposition is more interesting. 
\proclaim{Proposition} 
Let $\sigma\in \wwr F$ have dimension $p^r$, $r\ge 1$. Suppose that the decomposition function $\vS_\sigma$ has exactly one jump, at the point $a$, say. The following properties then hold:  
\roster 
\item 
the representation $\sigma$ is H-singular and $\roman{sl}_F(\chi) = a$, for every $\chi\in D^+(\sigma)$, $\chi\neq 1$; 
\item 
$\sw(\check\sigma\otimes \sigma) = p^{2r}\vS_\sigma(0) = (p^{2r}{-}1)a$;
\item 
if $\sigma$ is of Carayol type, then $a = \sw(\sigma)/(1{+}p^r)$. 
\endroster 
\endproclaim 
\demo{Proof} 
The definition (2.2.2) of $\vS_\sigma$ implies that 
$$
\vS_\sigma'(x) = \left\{\, \alignedat3 &p^{-2r}, &\quad &0< x < a, \\ &1, &\quad &x>a. \endalignedat \right. 
\tag 8.4.1 
$$ 
Consequently, the restriction of $\sigma$ to $\scr R_F(a)$ is irreducible and its restriction to $\scr R_F^+(a)$ is a multiple of a character. The group $S^+ = \sigma(\scr R_F^+(a))$ is therefore cyclic and central in $S = \sigma(\scr R_F(a))$. The finite $p$-group $S/S^+$ is a quotient of $\scr R_F(a)/\scr R_F^+(a)$, so it is elementary. Since $\sigma$ is irreducible on $\scr R_F(a)$, the centre of $S$ is cyclic. Consequently, the group $S = \sigma(\scr R_F(a))$ is H-cyclic with centre containing $S^+$. 
\par 
Let $C$ be the centralizer of $S$ in $P = \sigma(\scr P_F)$. Again, $C$ is finite cyclic. Let $y\in \scr P_F$. The representations $\sigma$, $\sigma^y$ are equivalent, particularly on $\scr R_F(a)$. The element $y$ must therefore act trivially on the centre of $S$. The commutator group $[y,\scr R_F(a)]$ is contained in $[\scr P_F,\scr R_F(a)] \subset \scr R_F^+(a)$, so $y$ acts trivially on $S$ modulo its centre. By 8.1 Lemma, there exists $x\in \scr R_F(a)$ such that $\sigma(xy)$ centralizes $S$. Therefore $P = SC$, implying that $\sigma$ is H-cyclic. 
\par 
It follows from (8.1.1) that $\check\sigma\otimes\sigma\Mid \scr P_F = \sum_{\chi\in D^+(\sigma)} \chi$. A non-trivial character $\chi\in D^+(\sigma)$ is non-trivial on $\scr P_F$ but it is trivial on the centre $C$ of $\sigma(\scr P_F)$, so $\chi$ is determined by its restriction to $\scr R_F(a)$. It is certainly trivial on $\scr R_F^+(a)$ so it has $F$-slope $a$. Thus $\sigma$ is H-singular and (1) is proven. Part (2) now follows from (8.1.1). Part (3) is 3.8 Proposition. \qed 
\enddemo 
We exhibit some implications of the preceding argument. 
\proclaim{Corollary} 
Let $Z = Z_\sigma$, $T = T_\sigma = T_I(\sigma)$ and $\sigma_T = \sigma \Mid \scr W_T$. 
\roster 
\item 
The field $Z$ is given by 
$$ 
\scr W_Z = \bigcap_{\chi\in D(\sigma_T)} \roman{Ker}\,\chi. 
$$ 
\item 
The Herbrand function $\psi_{Z/T}$ has a unique jump, lying at $e(T|F)a$. Moreover, 
\itemitem{\rm (a)} $\scr R_F^+(a) \i \scr W_Z$ and 
\itemitem{\rm (b)} $\scr W_T = \scr R_F(a)\scr W_Z$. 
\item 
The group $\scr W_T$ is the $\scr W_F$-centralizer of $\bar\sigma(\scr R_F(a))$. 
\endroster 
\endproclaim 
\demo{Proof}  
Define a field extension $Y/T$ by $\scr W_Y = \bigcap_\chi \roman{Ker}\,\chi$, with $\chi$ ranging over $D(\sigma_T)$. It follows from 8.1 Proposition that $Y/T$ is the centric field for the representation $\sigma_T$ and hence that $Y\supset Z$. We have to check that $\scr W_F$ acts trivially on $\sigma(\scr W_Y)$. However, $\sigma\Mid \scr W_Y = \sigma_T\Mid \scr W_Y$ is a multiple of a character, so that character is necessarily stable under $\scr W_F$. Therefore $Y = Z$ as required for (1).
\par 
Every non-trivial element of $D(\sigma_T)$ has Swan exponent $e(T|F)a$, whence follows the first assertion of (2). The same observation proves (a) while (b) follows from the definition of $Z$ via the group $D(\sigma_T)$. In (3), the group $\bar\sigma(\scr R_F(a))$ is the quotient of the H-cyclic group by its centre. The dual of this quotient is the character group $D(\sigma_T)$. Under the natural action of $\scr W_F$, the centralizer of this dual is $\scr W_T$, by 8.2 Lemma, implying the result. \qed 
\enddemo  
We finish with an example derived from \cite{11} and 8.2 Example. 
\subhead 
Example 
\endsubhead 
Take $p=2$, and suppose that $F$ contains a primitive cube root of unity. For $i=1,2$, let $\sigma_i \in \wwr F$ have dimension $2$ and satisfy $\sw(\sigma_i) = 1$. Theorem 5.1 of \cite{11} gives the recipe for $T_I(\sigma_i)$ and $D^+(\sigma_i)$. From that information and 8.2 Example, one sees it is possible to choose $\sigma_1$, $\sigma_2$ so that $\sigma = \sigma_1\otimes \sigma_2$ is irreducible and H-singular. It is not of Carayol type, as $\sw(\sigma) = 2$. If $[\sigma]_0^+ = \uL\vT$, $\vT \in \Scr E(F)$, then $\Psi_\vT$ has two jumps and is not convex: see 8.5 Example 1 of \cite{13} for the formula. 
\head \Rm 
9. Ramification structure 
\endhead 
Let $\sigma \in \wwr F$ be of Carayol type. We return to the methods of section 3 to work out the structure of $\sigma$ when restricted to an arbitrary ramification group of $\scr W_F$. If $[\sigma]_0^+ = \uL\vT$, $\vT \in \ewc F$, we get a formula for $\Psi_\vT$ to set against those of section 7. Despite appearances to the contrary, everything in this section relies on the local Langlands correspondence and the conductor formula of \cite{14}, since we use the main results of section 3. 
\subhead 
9.1 
\endsubhead 
To avoid carrying an irrelevant variable, we make a minor adjustment to our notation. If $\sigma\in \wwr F$ and if $\vT \in \Scr E(F)$ satisfies $[\sigma]^+_0 = \uL\vT$, we now write $\Psi_\sigma = \Psi_\vT$. 
\par 
Let $\sigma \in \wwr F$ be of Carayol type, and set $\vs = \vs_\sigma$. If $0<x<\vs$, define 
$$ 
w_\sigma(x) = \underset{\eps \to 0}\to{\text{\rm lim}}\, \Psi_\sigma'(x{+}\eps)/\Psi_\sigma'(x{-}\eps). 
\tag 9.1.1
$$ 
Thus $w_\sigma(x)$ is a non-negative power of $p$, and $w_\sigma(x) > 1$ if and only if $x$ is a jump of $\Psi_\sigma$. We then call $w_\sigma(x)$ the {\it height\/} of the jump $x$. 
\par 
Symmetry, as in 4.1 Proposition, gives an order-reversing involution $j\mapsto \bar\jmath$ on the set of jumps of $\Psi_\sigma$. If $\Psi_\sigma$ has an even number of jumps, this involution has no fixed point. If the number of jumps is odd, it fixes the middle one. In the notation of (9.1.1), the symmetry property of $\Psi_\sigma$ gives 
$$ 
w_\sigma(\bar\jmath) = w_\sigma(j). 
\tag 9.1.2 
$$ 
We will occasionally have to deal with the case of a one-dimensional representation $\sigma$. There, $\vS_\sigma(x) = \Psi_\sigma(x) = x$ and the functions $\vS_\sigma$, $\Psi_\sigma$ have no jumps. Indeed, the converse also holds \cite{13} 7.7. 
\subhead 
9.2 
\endsubhead 
Throughout the section, we use the following notation. 
\definition{Notation}
Let $\sigma \in \wwr F$ be of Carayol type and dimension $p^r$, $r\ge1$. Define $c_\sigma$ by $c_\sigma{+}\Psi_\sigma(c_\sigma) = \vs_\sigma$. Let 
$$ 
j_1<j_2< \dots < j_s < (c_\sigma) <\bar\jmath_s< \bar\jmath_{s-1} < \dots < \bar\jmath_1 
\tag 9.2.1
$$ 
be the jumps of $\Psi_\sigma$ with the understanding that 
\roster 
\item"\rm (a)" the term $c_\sigma$ is included only if $\Psi_\sigma$ has an odd number of jumps and 
\item"\rm (b)" $s=0$ when $\Psi_\sigma$ has only one jump. 
\endroster 
\enddefinition 
For the first version of the main result, we assume that $\sigma$ is {\it absolutely wild,\/} written $\sigma \in \awr F$, in the sense of 3.2 Definition. We deduce the final version, for $\sigma\in \wwr F$, in 9.5. 
\proclaim{Theorem} 
Let $\sigma \in \wwr F$ be\/ {\rm absolutely wild} of Carayol type and dimension $p^r$, $r\ge1$. 
\roster 
\item 
The restriction $\sigma\Mid \scr R_F^+(c_\sigma)$ is a direct sum of characters. 
\item 
Let $\xi$ be a character of $\scr R_F^+(c_\sigma)$ occurring in $\sigma$ and let $\scr W_{L_\xi}$ be the $\scr W_F$-stabilizer of $\xi$. Let $\sigma_\xi$ be the natural representation of $\scr W_{L_\xi}$ on the $\xi$-isotypic subspace of $\sigma$. 
\itemitem{\rm (a)} 
The field extension $L_\xi/F$ is absolutely wildly ramified (cf\. \rom{1.2}) of degree $p^rw_\sigma(c_\sigma)^{-\frac12}$ and $\scr W_{L_\xi}$ contains $\scr R_F^+(c_\sigma)$. 
\itemitem{\rm (b)} The representation $\sigma_\xi$ is irreducible, absolutely wild and 
$$ 
\sigma = \Ind_{L_\xi/F}\,\sigma_\xi. 
$$ 
\itemitem{\rm (c)} If $c_\sigma$ is not a jump of $\Psi_\sigma$, then $\sigma_\xi$ is a character. Otherwise, $\sigma_\xi$ is H-singular, of Carayol type and dimen\-sion $w_\sigma(c_\sigma)^{\frac12}$. The unique jump of\/ $\Psi_{\sigma_\xi}$ lies at $\psi_{L_\xi/F}(c_\sigma)$. 
\endroster 
\endproclaim 
\remark{Remarks} 
\roster 
\item 
The triple $(\xi,L_\xi,\sigma_\xi)$ is uniquely determined by $\sigma$, up to $\scr W_F$-conjugation. 
\item 
The function $\Psi_\sigma$ has no jump lying strictly between $c_\sigma$ and $\bar\jmath_s$. So, if $\xi$ and $\xi'$ are components of $\sigma\Mid \scr R_F^+(c_\sigma)$, then $\xi = \xi'$ if and only if $\xi\Mid \scr R_F(\bar\jmath_s) = \xi'\Mid \scr R_F(\bar\jmath_s)$. If $c_\sigma$ is not a jump then, in the same way, $\sigma\Mid \scr R_F^+(j_s)$ is a sum of characters, two of which are equal if and only if their restrictions to $\scr R_F(\bar\jmath_s)$ are equal. 
\endroster 
\endremark 
As we prove the theorem, we uncover further features of interest that we now list. 
\proclaim{Complement 1} 
Let $1\le k\le s$. 
\roster 
\item 
The restriction $\sigma\Mid \scr R_F^+(j_k)$ is a multiplicity-free direct sum of irreducible representations. 
\item 
The restriction $\sigma\Mid \scr R_F(\bar\jmath_k)$ is a direct sum of characters. The isotypic components of $\sigma\Mid \scr R_F(\bar\jmath_k)$ are the subspaces $\tau\Mid \scr R_F(\bar\jmath_k)$, as $\tau$ ranges over the irreducible components of $\sigma\Mid \scr R_F^+(j_k)$. 
\endroster 
\endproclaim 
In light of Remark (2) above, one can equally relate the decompositions of $\sigma\Mid \scr R_F(j_k)$ and $\sigma\Mid \scr R^+_F(\bar\jmath_k)$. In the next result, we use the concept of {\it elementary resolution\/} from 1.9.  
\proclaim{Complement 2} 
For $1\le k\le s$, choose an irreducible component $\tau_k$ of the restriction $\sigma\Mid \scr R_F^+(j_k)$ so that $\tau_{k+1}$ is a component of $\tau_k\Mid \scr R_F^+(j_{k+1})$, $1\le k <s$. Let $\scr W_{E_k}$ be the $\scr W_F$-stabilizer of $\tau_k$. 
\roster 
\item"\rm (a)" If $\xi$ is a character of $\scr R_F^+(c_\sigma)$ occurring in $\tau_s\Mid \scr R_F^+(c_\sigma)$, then $E_s = L_\xi$. 
\item"\rm (b)" 
The Herbrand function satisfies 
$$ 
\Psi_\sigma(x) = p^{-r}\psi_{L_\xi/F}(x), \quad 0\le x\le c_\sigma. 
\tag 9.2.2 
$$  
\item"\rm (c)" 
The tower of fields 
$$ 
F\i E_1\i E_2 \i \dots \i E_s = L_\xi 
\tag 9.2.3 
$$ 
is the elementary resolution of $L_\xi/F$. 
\endroster 
\endproclaim 
Following 4.1 Proposition, symmetry implies that the relation (9.2.2) determines $\Psi_\sigma(x)$ for $0\le x\le \vs_\sigma$. The tower of fields (9.2.3) is uniquely determined by $\sigma$, up to $\scr W_F$-conjugacy. 
\par 
The proofs of these results occupy 9.3 and 9.4. 
\subhead 
9.3 
\endsubhead 
The theorem of 9.2 is proved by induction on the number of jumps. If $\Psi_\sigma$ has no jumps, then $\dim\sigma = 1$ and this case has been excluded. If $\Psi_\sigma$ has just one jump, the theorem and its complements follow from 8.4 Proposition with $L_\xi = F$. 
\par In this sub-section, we assume that $\Psi_\sigma$ has at least two jumps and give a reduction step concerned only with the outermost jumps. As in 8.2, let $D(\sigma)$ be the group of characters $\chi$ of $\scr W_F$ such that $\chi\otimes\sigma \cong \sigma$. 
\proclaim{Proposition} 
Let $\sigma\in \awr F$ be of Carayol type. Suppose that $\Psi_\sigma$ has at least two jumps, of which $a$ is the first and $z$ the last. Let $D_a(\sigma)$ be the group of $\chi\in D(\sigma)$ for which $\sw(\chi)\le a$. Let $E_1/F$ be class field to $D_a(\sigma)$. 
\roster 
\item 
The group $D_a(\sigma)$ is elementary abelian of order $w_\sigma(a)$. 
\item 
The group $\scr R_F^+(a)$ is contained in $\scr W_{E_1}$ and $\scr W_F = \scr R_F(a)\scr W_{E_1}$. 
\item 
There exists $\sigma_1\in \awr{E_1}$ such that $\sigma =  \Ind_{E_1/F}\,\sigma_1$. Moreover, 
$$ 
\sigma\,\big|\,\scr R_F^+(a) = \sum_{\gamma\in \Gal{E_1}F} \sigma_1^\gamma\,\big|\,\scr R_F^+(a). 
\tag 9.3.1 
$$ 
The representations $\sigma_1^\gamma\,\big|\,\scr R_F^+(a)$, $\gamma \in \Gal{E_1}F$, are distinct and irreducible. The $\scr W_F$-stabilizer of $\sigma_1\,\big|\, \scr R_F^+(a)$ is $\scr W_{E_1}$. 
\item 
The jumps of $\Psi_{\sigma_1}$ are $\psi_{E_1/F}(j)$, as $j$ ranges over the jumps of $\Psi_\sigma$, $j\neq a,z$. Indeed, $w_\sigma(y) = w_{\sigma_1}(\psi_{E_1/F}(y))$, for $y\neq a,z$. 
\item 
The restriction $\sigma \Mid \scr R_F(z)$ is a direct sum of characters $\xi$. The $\scr W_F$-stabilizer of any such $\xi$ is $\scr W_
{E_1}$. 
\endroster 
\endproclaim  
\demo{Proof} 
The group $D_a(\sigma)$ is non-trivial (3.3 Lemma 1), so choose $\chi\in D_a(\sigma)$, $\chi\neq 1$. Set $\scr W_K = \roman{Ker}\,\chi$. The extension $K/F$ is cyclic of degree $p$, and $\psi_{K/F}$ has a unique jump, lying at $a$. As in 3.3 Lemma 2, $\scr W_K\cap \scr R_F(a) = \scr R_K(a)$ is of index $p$ in $\scr R_F(a)$ and $\scr R_K^+(a) = \scr R_F^+(a)$. There exists $\tau\in \wwr K$ with $\sigma = \Ind_{K/F}\,\tau$, the representation $\tau$ being either absolutely wild of Carayol type (3.3 Lemma 1) or a character. By 3.4 Theorem (2) and 3.5 Theorem, $w_\sigma(y) = w_\tau(\psi_{K/F}(y))$, provided $y\neq a,z$. On the other hand, $w_\sigma(a) = pw_\tau(a)$ and $w_\sigma(z) = pw_\tau(\psi_{K/F}(z))$ {\it loc\. cit.} 
\remark{Note} 
Since $\sigma$ has at least two jumps, 3.5 Corollary 1 shows that the case of 3.4 Theorem (3) need not be considered here. 
\endremark 
\proclaim{Lemma} 
\roster 
\item 
The restriction $\tau\,\big|\,\scr R_K(a)$ is irreducible and 
$$ 
\sigma\,\big|\,\scr R_F^+(a) = \sum_{\delta\in \Gal KF} \tau^\delta\,\big|\,\scr R_F^+(a). 
\tag 9.3.2 
$$ 
\item 
The representations $\tau^\delta\, \big|\,\scr R_F^+(a)$, $\delta\in \Gal KF$, are disjoint.  
\endroster 
\endproclaim 
\demo{Proof} 
Since $a$ is the first jump of $\vS_\sigma$, the restriction $\sigma\,\big|\,\scr R_F(a)$ is irreducible. The Mackey formula implies that this restriction is $\Ind_{\scr R_K(a)}^{\scr R_F(a)}\,\tau\Mid \scr R_K(a)$, whence the first assertion follows. The relation (9.3.2) again follows from the Mackey formula. 
\par 
Since $\sigma\Mid \scr R_F(a)$ is irreducible, the irreducible components of $\sigma \Mid \scr R_F^+(a)$ are all conjugate and occur with the same multiplicity. So, in (2) the representations $\tau^\delta \Mid \scr R_F^+(a)$ are either disjoint or identical. We show they are disjoint. 
\par 
Let $\Delta_K$ be the canonical ultrametric on $\scr W_K\backslash \wP K$. Let $\delta \in \Gal KF$, $\delta \neq 1$. By 3.5 Theorem (and recalling the definition (3.4.2)) we have 
$$ 
\Delta_K(\tau,\tau^\delta) = \psi_{K/F}(z) > \psi_{K/F}(a) = a. 
\tag 9.3.3 
$$ 
The representations  $\tau^\delta \Mid\scr R_F^+(a)$,  $\tau\Mid\scr R_F^+(a)$ are therefore disjoint, as asserted. \qed 
\enddemo 
\remark{Remark} 
The relation (9.3.3) implies that $\tau$ and $\tau^\delta$ are disjoint on $\scr R_F(z)$. 
\endremark 
We proceed by induction on the integer $w_\sigma(a)$. Suppose first that $w_\sigma(a) = p$, whence $w_\tau(a) =1$. We prove the proposition with $E_1 = K$ and $\sigma_1 = \tau$. As $a$ is not a jump of $\vS_\tau$ (giving point (4)), so $\tau$ is irreducible on $\scr R_F^+(a) = \scr R_K^+(a)$. It follows that $D(\tau)$ has no element of Swan exponent $a$. The conjugates $\tau^\delta$, $\delta\in \Gal KF$, are disjoint on $\scr R_K^+(a)$, by the lemma.  Consequently, $D_a(\sigma)$ has order $p = w_\sigma(a)$. The point $\psi_{K/F}(z)$ is not a jump of $\vS_\tau$, by 3.4 Theorem again. It follows that $\tau\,\big|\,\scr R_F(z)$ is a multiple of a character. Thus 
$$ 
\sigma\,\big|\,\scr R_F(z) = \sum_{\delta\in \Gal KF} \tau^\delta\,\big|\,\scr R_F(z) 
$$ 
is a sum of characters. Since $z$ is a jump of $\Psi_\sigma$, these characters cannot all be the same: they fall into $p$ distinct orbits under $\scr W_F$. Assertion (5) follows, and the proof is complete in the case $w_\sigma(a) = p$. 
\par 
Suppose next that $w_\sigma(a)$ is divisible by $p^2$. In particular, $\tau$ is not a character. Inductively, we may assume that the result holds for the representation $\tau\in \awr K$. For convenience, we expand this assumption and fix some notation.   
\proclaim{Inductive hypothesis} 
Let $E/K$ be class field to the group $D_a(\tau)$. Let $\zeta\in \awr E$ satisfy $\Ind_{E/K}\,\zeta = \tau$. 
\roster 
\item 
The group $D_a(\tau)$ is elementary abelian of order $w_\tau(a)$. 
\item
The group $\scr R_K^+(a)$ is contained in $\scr W_E$ and $\scr W_K = \scr R_K(a)\scr W_E$. 
\item 
In the expansion 
$$ 
\tau\,\big|\,\scr R_K^+(a) = \sum_{\gamma\in \Gal EK} \zeta^\gamma\,\big|\,\scr R_K^+(a), 
\tag 9.3.4 
$$ 
the terms $\zeta^\gamma\,\big|\,\scr R_F^+(a)$, $\gamma \in \Gal EK$, are distinct and irreducible. The $\scr W_K$-stabilizer of $\zeta\,\big|\,\scr R_K^+(z)$ is $\scr W_E$. 
\item 
The jumps of $\Psi_{\zeta}$ are $\psi_{E/K}(k)$, as $k$ ranges over the jumps of $\Psi_\tau$, $k\neq a, \psi_{K/F}(z)$. Indeed, $w_\tau(y) = w_\zeta(\psi_{E/K}(y))$, for $y\neq a, \psi_{K/F}(z)$. 
\item 
The restriction of $\tau$ to $\scr R_F(z) = \scr R_K(\psi_{K/F}(z))$ is a direct sum of characters $\xi$. The $\scr W_K$-stabilizer of any such $\xi$ is $\scr W_E$. 
\endroster 
\endproclaim 
We prove that $E/F$ is class field to $D_a(\sigma)$. Each of the functions $\psi_{K/F}$, $\psi_{E/K}$ has a unique jump, lying at $a$. The same therefore applies to $\psi_{E/F}$. The field $E$ appears as a subfield of the centric field of $\sigma$, so $E/F$ is absolutely wild. As $\psi_{E/F}$ has a unique jump, lying at $a$, the case $k=1$ of 1.9 Corollary 1 implies that $E/F$ is elementary abelian and so every element $\phi$ of $D_{(1)}(E|K)$ (notation of 1.9 Proposition) extends to a character $\tilde\phi$ of $\scr W_F$ lying in $D_{(1)}(E|F)$. We have $\tilde\phi\otimes\sigma = \Ind_{K/F}\,\phi\otimes\tau = \Ind_{K/F}\,\tau = \sigma$. That is, $\tilde\phi \in D_a(\sigma)$, whence $D_a(\sigma) = D_{(1)}(E|F)$ and this group has order $pw_\tau(a) = w_\sigma(a)$. 
\par 
We have proved part (1) of the proposition, with $E_1 = E$. Part (2) follows from the relation $\psi_{E/F}(a) = a$. The lemma applies equally here, so the irreducible representations 
$$ 
 \zeta^{\gamma\delta} \Mid \scr R_F^+(a), \quad \gamma\in \Gal EK,\ \delta \in \Gal KF, 
$$ 
are disjoint. Part (3) of the proposition now follows by induction. 
\par 
Part (4) of the proposition follows directly from part (4) of the inductive hypothesis. It remains to prove part (5). By part (5) of the inductive hypothesis, $\sigma\,\big|\,\scr R_F(z)$ is a sum of characters. The representations $\tau^\delta$, $\delta \in \Gal KF$, are disjoint on $\scr R_F(z)$ by the remark following (9.3.3). The result follows from the inductive hypothesis, with $E_1 = E$ and $\sigma_1 = \zeta$. \qed 
\enddemo 
\subhead 
9.4  
\endsubhead 
We prove 9.2 Theorem and its complements. We proceed by induction on the number of jumps of $\Psi_\sigma$. 
\demo{Proof of Theorem}
When $\Psi_\sigma$ has at most one jump, there is nothing more to say. We therefore assume, in the notation (9.2.1), that $s\ge 1$. We apply 9.3 Proposition to get a Galois extension $E_1/F$ and a representation $\sigma_1\in \awr{E_1}$ such that $\sigma = \Ind_{E_1/F}\,\sigma_1$. The extension $E_1/F$ has a unique jump, lying at $j_1$, so $\scr R^+_F(x) \i \scr W_{E_1}$ for $x\ge j_1$. The function $\Psi_{\sigma_1}$ has jumps at $\psi_{E_1/F}(j)$, where $j$ ranges over all jumps of $\Psi_\sigma$, subject to $j\neq j_1, \bar\jmath_1$. 
\par 
Suppose the number of jumps to be even. Assume to start with that this number is two, that is, $s=1$. In 9.3 Proposition, the representation $\sigma_1$ is a character. The conjugates $\sigma_1^\gamma$, $\gamma\in \Gal{E_1}F$, agree on $\scr R_F^+(\bar \jmath_1)$ but are distinct on $\scr R_F(\bar \jmath_1)$. All assertions of the theorem follow readily in this case. We therefore assume that $s\ge2$. By inductive hypothesis, $\sigma_1\Mid \scr R_F^+(j_s)$ is a sum of characters, so the same applies to $\sigma\Mid \scr R_F^+(j_s)$. Part (1) is done in this case. The field $L = L_\xi$ appears as a subfield of the centric field of $\sigma$, so $L/F$ is absolutely wild. The inductive hypothesis gives a character $\rho_1$ of $\scr W_L$ which induces $\sigma_1$. It follows that $\Ind_{L/F}\,\rho_1 = \sigma$, and $\rho_1$ has the necessary properties relative to $\sigma$. This proves part (2) of the theorem when the number of jumps is even. 
\par 
Suppose that the number of jumps is odd. Thus, by inductive hypothesis, $\sigma_1\Mid \scr R_{E_1}(c_{\sigma_1})$ is not a sum of characters, while $\sigma_1\Mid \scr R^+_{E_1}(c_{\sigma_1})$ is such. Since $\scr R_{E_1}(c_{\sigma_1}) = \scr R_F(\vf_{E_1/F}(c_{\sigma_1}))$, the point $\vf_{E_1/F}(c_{\sigma_1})$ is a jump of $\Psi_\sigma$. That is, $c_\sigma = \vf_{E_1/F}(c_{\sigma_1})$ and we have proved part (1) of the theorem. Assertions (2)(a)--(c) now follow by induction, exactly as in the first case, on noting that $\dim\sigma_\xi = w_\sigma(c_\sigma)^{\frac12}$ by 8.1 Proposition. \qed 
\enddemo  
\demo{Proof of Complement 1} 
We follow 9.3 Proposition to write $\sigma = \Ind_{E_1/F}\,\sigma_1$. That result also shows that $\sigma\Mid \scr R_F^+(j_1)$ is multiplicity-free. For $2\le k\le s$, the restriction $\sigma_1\Mid \scr R_F^+(j_k)$ is multiplicity-free by the inductive hypothesis. The relation $w_{\sigma_1}(\psi_{E_1/F}(j_k)) = w_\sigma(j_k)$ shows that $\sigma\Mid \scr R_F^+(j_k)$ is multiplicity-free, and we have proved part (1). 
\par 
The first assertion of (2) follows from part (1) of the theorem. Since $\bar\jmath_1$ is the last jump of $\Psi_\sigma$, the restriction $\sigma_1\Mid \scr R_F(\bar\jmath_1)$ is a multiple of a character while $\sigma\Mid \scr R_F(\bar\jmath_1)$ is a direct sum of characters. The number of isotypic components in $\sigma\Mid \scr R_F(\bar\jmath_1)$ is $w_\sigma(\bar\jmath_1) = w_\sigma(j_1) = [E_1{:}F]$, by 9.3 Proposition, whence the result follows. \qed 
\enddemo 
\demo{Proof of Complement 2} 
Recall that $E_1/F$ was defined in 9.3 as class field to the group $D_{j_1}(\sigma)$ of characters $\chi$ of $\scr W_F$ such that $\chi\otimes\sigma \cong \sigma$ and $\sw(\chi) \le j_1$. Thus $E_1/F$ is Galois and, by 9.3 Proposition (3), $\scr W_{E_1}$ is the $\scr W_F$-stabilizer of any irreducible component of $\sigma\Mid \scr R^+_F(j_1)$. In the first instance, we may therefore choose the extension $L = L_\xi/F$ of the theorem, within its conjugacy class, so that $E_1 \i L_\xi$. Since all choices of $\xi$ are $\scr W_F$-conjugate and $E_1/F$ is Galois, we have $E_1\i L_\xi$ for all $\xi$. That is, $E_1 \i L$. 
\par 
Because of the relation $\sigma = \Ind_{L/F}\,\rho_\xi$, a character $\phi$ of $\scr W_F$ with $\phi\Mid \scr W_L$ trivial must satisfy $\phi\otimes \sigma\cong\sigma$. The definition of $E_1$ in 9.3 implies that $j_1$ is the least jump of $\psi_{L/F}$. By 1.9 Proposition (3), $E_1/F$ is the first step in the elementary resolution of $L/F$. Parts (a) and (c) of Complement 2 now follow by induction. 
\par 
In the proof of the theorem, we showed that $c_{\sigma_1} = \psi_{E_1/F}(c_\sigma)$. From 3.4 Theorem we conclude that 
the jumps of $\Psi_{\sigma_1}$ are 
$$ 
\multline 
\psi_{E_1/F}(j_2) < \psi_{E_1/F}(j_3) < \dots < \psi_{E_1/F}(j_s) \\ < (\psi_{E_1/F}(c_{\sigma_1})\big) 
< \psi_{E_1/F} (\bar\jmath_s) < \dots < \psi_{E_1/F}(\bar\jmath_2), 
\endmultline 
$$ 
with the same convention regarding the central entry in the list. Moreover, 
$$ 
w_{\sigma_1}(\psi_{E_1/F}(j_k)) = w_\sigma(j_k), \quad 2\le k\le s, 
\tag 9.4.1 
$$ 
and similarly relative to the central jump. Let $w_1 = w_\sigma(j_1)$, so that $w_1 = [E_1{:}F]$. The functions $\Psi_\sigma(x)$, $w_1^{-1}\Psi_{\sigma_1}(\psi_{E_1/F}(x))$ have the same jumps in the region $0\le x\le c_\sigma$. The heights (9.1) of these jumps are the same, and the functions agree on a region $0\le x < \ve$. We conclude by induction that 
$$ 
\align 
\Psi_\sigma(x) &= w_1^{-1}\Psi_{\sigma_1}(\psi_{E_1/F}(x)) \\ 
&= p^{-r}\psi_{L/F}(x), \quad 0\le x\le c_\sigma. 
\endalign 
$$ 
This proves part (b). \qed 
\enddemo 
\subhead 
9.5 
\endsubhead 
We extend the results of 9.2 to representations of Carayol type that are totally, but not necessarily absolutely, wild. The notational conventions of 9.1, 9.2 remain in force. 
\proclaim{Corollary} 
Let $\sigma\in \wwr F$ be of Carayol type and dimension $p^r$, $r \ge 1$. Define $c_\sigma$ by the equation $c_\sigma{+} \Psi_\sigma(c_\sigma) = \vs_\sigma$. 
\roster 
\item 
The representation $\sigma\Mid \scr R_F^+(c_\sigma)$ is a direct sum of characters. 
\endroster 
Let $\xi$ be a character of $\scr R_F^+(c_\sigma)$ occurring in $\sigma$. Let $\scr W_{L_\xi}$ be the $\scr W_F$-stabilizer of $\xi$ and let $\sigma_\xi$ be the natural representation of $\scr W_{L_\xi}$ on the $\xi$-isotypic subspace of $\sigma\Mid \scr R_F^+(c_\sigma)$. 
\roster 
\item"\rm (2)" The representation $\sigma_\xi$ is irreducible and $\Ind_{L_\xi/F}\,\sigma_\xi = \sigma$. Moreover, 
\itemitem{\rm (a)} $\dim\sigma_\xi = w_\sigma(c_\sigma)^{1/2}$, and 
\itemitem{\rm (b)} if $\dim\sigma_\xi>1$, then $\sigma_\xi$ is totally wild, H-singular and of Carayol type. 
\item"\rm (3)" 
The field extension $L_\xi/F$ is totally ramified of degree $p^r/\dim\sigma_\xi$ and 
$$ 
\Psi_\sigma(x) = p^{-r}\psi_{L_\xi/F}(x), \quad 0\le x\le c_\sigma. 
\tag 9.5.1 
$$ 
\endroster 
\endproclaim 
\demo{Proof} 
Let $T = T_\sigma/F$ be the tame centric field of $\sigma$. Thus $\tau = \sigma\Mid \scr W_T$ is absolutely wild of Carayol type. If $e = e(T|F)$, then $\Psi_\sigma(x) = \Psi_\tau(ex)/e$ and $\vs_\tau = e\vs_\sigma$, so $c_\tau = ec_\sigma$. 
\par
Consequently, $\scr R^+_F(c_\sigma) = \scr R^+_T(c_\tau)$ and part (1) follows from part (1) of 9.2 Theorem. All choices of $\xi$ are $\scr W_F$-conjugate so let us fix one and write $L_\xi = L$. The $\scr W_T$-stabilizer of $\xi$ is $\scr W_T \cap \scr W_L = \scr W_{LT}$. The natural representation of $\scr W_{LT}$ on the $\xi$-isotypic subspace of $\tau\Mid \scr R_T^+(c_\tau)$ is $\sigma_\xi\Mid \scr W_{LT}$, which is irreducible. It follows that $\sigma_\xi$ is irreducible and has properties (2)(a), (2)(b). Moreover, $\scr R_F(c_\sigma)$ is contained in $\scr W_L$ and $\sigma_\xi \Mid \scr R_F(c_\sigma)$ is irreducible. 
\par 
The degree $[L{:}F]$ is the number of distinct characters occurring in the representation $\sigma\Mid \scr R_F^+(c_\sigma) = \tau\Mid \scr R_F^+(c_\sigma)$, so $[L{:}F] = [LT{:}T]$ and $L/F$ is totally wildly ramified. Further, 
$$ 
\Ind_{L/F}\,\sigma_\xi \Mid \scr W_T = \Ind_{LT/T}\,(\sigma_\xi\Mid \scr W_{LT}). 
$$ 
This restriction is irreducible, whence $\Ind_{L/F}\,\sigma_\xi$ is irreducible and equivalent to $\sigma$. Finally, for $0\le x\le c_\sigma$, 
$$ 
\Psi_\sigma(x) = \Psi_\tau(ex)/e = p^{-r} \psi_{LT/T}(ex)/e = p^{-r}\psi_{L/F}(x), 
$$ 
by (2.2.3), 9.2 Complement 1 and 1.1 Lemma. \qed 
\enddemo 
\proclaim{Complement} 
If $\sigma \in \wwr F$, the assertions of 9.2 Complement 1 apply unchanged. 
\endproclaim 
\demo{Proof} 
Take $T/F$, $e = e(T|F)$ and $\tau = \sigma \Mid \scr W_T$ as in the proof of the corollary. Thus $\scr R_F(x) = \scr R_T(ex)$, $\scr R^+_F(x) = \scr R_T^+(ex)$, for all $x>0$. So, for $x>0$, the decomposition structures of $\sigma\Mid \scr R_F(x)$ and $\sigma\Mid \scr R_F^+(x)$ are identical to those of $\scr R_T(ex)$ and $\tau\Mid \scr R_T^+(ex)$. \qed 
\enddemo 
\remark{Remark} 
Let $K/F$ be a finite tame extension and set $e = e(K|F)$. We may view $\xi$ as a character of $\scr R_K^+(ec_\sigma) = \scr R_K^+(c_{\sigma_K})$, where $\sigma_K = \sigma\Mid \scr W_K$. The arguments in the proof of 9.5 Corollary show that $L_{\sigma_K,\xi} = KL_{\sigma,\xi}$, in the obvious notation. 
\endremark 
\subhead 
9.6 
\endsubhead 
We continue with the notation of 9.5 Corollary, and look into the structure of the inducing representation $\sigma_\xi$. This is in preparation for a more detailed discussion in the next section. 
\definition{Definition} 
Let $\widetilde L_{\sigma,\xi}/L_\xi$ be the centric field of the representation $\sigma_\xi \in \wwr{L_\xi}$. 
\enddefinition 
The extension $\widetilde L_{\sigma,\xi}/L_\xi$ is Galois and $\widetilde L_{\sigma,\xi}/F$ is uniquely determined by $\sigma$, up to conjugation in $\scr W_F$. 
\proclaim{Proposition} 
Suppose $\sigma$ is absolutely wildly ramified. The extension $\widetilde L_{\sigma,\xi}/L_\xi$ is totally ramified and elementary abelian of degree $(\dim\sigma_\xi)^2$. If $\widetilde L_{\sigma,\xi} \neq L_\xi$, the extension $\widetilde L_{\sigma,\xi}/L_\xi$ has a unique ramification jump, lying at $\psi_{L_\xi/F}(c_\sigma)$. In particular, $\scr R_F^+(c_\sigma) \subset \scr W_{\widetilde L_{\sigma,\xi}}$. 
\endproclaim 
\demo{Proof} 
All assertions are trivial if $\sigma_\xi$ is a character, so assume otherwise. By 9.2 Theorem, the representation $\sigma_\xi$ of $\scr W_{L_\xi}$ is absolutely wild and H-singular. Thus $\wt L_{\sigma,\xi}/L_\xi$ is totally ramified and elementary abelian of degree $(\dim\sigma_\xi)^2$. By 8.1 Proposition, it is class field to the character group $D(\sigma_\xi)$. The unique ramification jump of $\sigma_\xi$ lies at $\psi_{L_\xi/F}(c_\sigma)$ (9.2 Theorem again), so every non-trivial element of $D(\sigma_\xi)$ has Swan exponent $\psi_{L_\xi/F}(c_\sigma)$ (8.3 Proposition). Therefore $\scr W_{\wt L_{\sigma,\xi}} \supset \scr R_{L_\xi}^+(\psi_{L_\xi/F}(c_\sigma)) = \scr R_F^+(c_\sigma)$. \qed 
\enddemo 
In the general case $\sigma \in \wwr F$, the extension $\widetilde L_{\sigma,\xi}/L_\xi$ is not totally wildly ramified. We recall the standard example. 
\example{Example} 
For this example, we adhere to the classical framework of the exposition in section 41 of \cite{8}. Take $p=2$, and let $\sigma\in \wW F$ be {\it primitive\/} of dimension $2$. The representation $\sigma$ is then totally ramified and H-singular. After twisting with a character, if necessary, we may assume that $\sigma$ is of Carayol type. In terms of the preceding discussion, $\Psi_\sigma$ has one jump and $L_\xi = F$. Using standard notation for permutation groups, $\bar\sigma(\scr W_F)$ is either $A_4$ (if $F$ has a primitive cube root of unity) or $S_4$ (otherwise). The tame centric field $T_\sigma/F$ is cyclic of degree $3$ in the first case and, in the second, $\Gal {T_\sigma}F \cong S_3$. 
\endexample 
\subhead 
9.7 
\endsubhead 
As an application of the methods of this section, we return to the question posed in 8.2 Remark. If $\sigma \in \wwr F$, we use the notation $D(\sigma)$, $D^+(\sigma)$, $T_I(\sigma)$ introduced in section 8. In addition, $T(\sigma)$ shall be the tame centric field of $\sigma$. 
\proclaim{Application} 
There exist a field $F$, of residual characteristic $2$, and a representation $\sigma\in \wwr F$ such that $T_I(\sigma) \neq T(\sigma)$. One may take $\sigma$ to be of Carayol type and dimension $4$. 
\endproclaim 
\demo{Proof} 
Let $F$ have residual characteristic $2$. Let $K/F$ be totally ramified of degree $4$, such that $\psi_{K/F}$ has two jumps $a<b$, of which $a$ is an odd integer. Replacing $F$ by $E$ and $K/F$ by $EK/E$, where $E/F$ is finite and tamely ramified, we may assume $b{-}a$ to be as large as necessary without affecting the parity of $a$. 
\par 
Let $m$ be a positive integer and define $c = c_m$ by the equation $4c+\psi_{K/F}(c) = m$. 
\proclaim{Lemma} 
If $b{-}a$ is sufficiently large, one may choose $m$ so that 
\roster 
\item $a<c_m<b$, 
\item $m\not\equiv 2a  \pmod 3$, 
\item $m \equiv a {+} 2 \pmod 4$
\endroster 
\endproclaim 
This is clear. Assume it has been done, and note that $m$ is odd. We get 
$$ 
c = c_m = (m{-}2a)/6. 
$$ 
The bi-Herbrand function $\Psi = \biP KF{m/4}$ has three jumps, namely $a$, $c$ and $z$, satisfying $a<c<z$. By 7.2 Corollary, there exists $\vT \in \ewc F$ such that $\Psi(x) = \Psi_\vT(x)$, $0\le x\le m/4$. Choose $\sigma \in \wwr F$ such that $[\sigma]_0^+ = \uL\vT$. We show that $\sigma$ has the desired properties. 
\par 
Let $\phi\in D^+(\sigma)$, $\phi \neq 1$. The $F$-slope $\roman{sl}_F(\phi)$ of $\phi$, as in (8.3.1), can only take a value $a,c,z$ ({\it cf\.} 8.1 Proposition of \cite{13}). Suppose $\roman{sl}_F(\phi) = a$. The jump $a$ has height $2$, so there is only one possibility for $\phi$. Since $a$ is an integer, the $\scr W_F$-stabilizer of $\phi\Mid \scr R_F(a)$ is of the form $\scr W_E$, where $E/F$ is unramified. The character $\phi\in D^+(\sigma)$ is completely determined by its restriction to $\scr R_F(a)$ so $\scr W_E$ is the $\scr W_F$-stabilizer of $\phi$. So, writing $\sigma_E = \sigma\Mid \scr W_E$, there exists a unique character $\tilde\phi\in D(\sigma_E)$ such that $\tilde\phi\Mid \scr P_F = \phi$ (8.2 Lemma). Thus $D_a(\sigma_E)$ has order $2$. 
\par 
Suppose next that $\roman{sl}_F(\phi) = c = (m{-}2a)/6$. The conditions imposed on $m$ imply $3c\in \frac12\Bbb Z \smallsetminus \Bbb Z$. We conclude that there is no finite tame extension $E/F$ for which $\phi$ extends to a character of $\scr W_E$. Finally, consider the case where $\roman{sl}_F(\phi) = z$. By 3.5 Theorem, $z = (m{-}a)/4 \in \frac12\Bbb Z\smallsetminus \Bbb Z$ and the same conclusion holds. We have shown: 
\proclaim{Proposition} 
The group $D^+(\sigma)$ has order $2$ and there is a finite unramified extension $E/F$ such that every character $\phi\in D^+(\sigma)$ is fixed by $\scr W_E$. Further, 
\roster 
\item 
$D(\sigma_E) =  D_a(\sigma_E)$, where $\sigma_E = \sigma \Mid \scr W_E$ and 
\item  
$T_I(\sigma)/F$ is unramified. 
\endroster 
\endproclaim 
We now follow the procedure of 9.5 to choose a character $\xi$ of $\scr R_F^+(c)$ occurring in $\sigma\Mid \scr R_F^+(c)$. We set $L = L_\xi$ and $\tau = \sigma_\xi$. We have $\sigma  = \Ind_{L/F}\,\tau$. Since $\sw(\sigma) =  m$ and $w_{L/F} = a$, we get $\sw(\tau) = m{-}2a \not\equiv 0 \pmod 3$. The Herbrand function $\Psi_\tau$ has a unique jump, which lies at $(m{-}2a)/3$ (8.4 Proposition). It follows that $e(T_I(\tau)|L)$ is divisible by $3$. This implies that $e(T(\sigma)|F)$ is divisible by $3$, whence  $T(\sigma) \neq T_I(\sigma)$. \qed 
\enddemo 
\remark{Remark} 
The choice of $p=2$ in the example is for simplicity only. There is nothing special about the case $p=2$ in this context. 
\endremark 
\head\Rm 
10. Parameter fields  
\endhead 
Let $[\frak a,m,0,\alpha]$ be a simple stratum in $\M{p^r}F$, $r\ge 1$, with the usual properties: $F[\alpha]/F$ is totally ramified of degree $p^r$ and $m = -\ups_{F[\alpha]}(\alpha)$ is not divisible by $p$. Let 
$$ 
\scr G^\star(\alpha) = \{\sigma\in \wwr F: [\sigma]^+_0 \in \uL\|\scr C^\star(\frak a,\alpha)\|\}. 
$$ 
Observe that every $\sigma \in \scr G^\star(\alpha)$ has dimension $p^r$. 
\par 
If $\sigma \in \scr G^\star(\alpha)$ and $[\sigma]_0^+ = \uL\vT$, we have two determinations of $\Psi_\vT$, from 7.2 Theorem 2 and 9.5 Corollary respectively. In 9.5 and 9.6, we attached to $\sigma$ a tower of fields $F\subset L_\xi \subset \widetilde L_{\sigma,\xi}$, given by a character $\xi$ of $\scr R_F^+(c_\alpha)$ occurring in $\sigma$. This configuration is determined by $\sigma$ up to $\scr W_F$-conjugation. We now examine how it varies when $\vT$ ranges over $\|\scr C^\star(\frak a,\alpha)\|$. 
\subhead 
10.1 
\endsubhead 
We fix notation for the rest of the section. With $[\frak a,m,0,\alpha]$ as above, we abbreviate  
$$ 
\alignedat3 
\vs_\alpha &= m/p^r, &\quad w_\alpha &= w_{F[\alpha]/F},\\ 
l_\alpha &= \roman{max}(0,m{-}w_\alpha), &\quad \lambda_\alpha &= [l_\alpha/2]. 
\endalignedat 
\tag 10.1.1 
$$ 
By 7.2 Theorem 2, $\Psi_\vT(x) = \biP {F[\alpha]}F{\vs_\alpha}(x)$, for all $\vT \in \|\scr C^\star(\frak a,\alpha)\|$ and $0\le x\le \vs_\alpha$. We use the notation 
$$ \aligned
\biP {F[\alpha]}F{\vs_\alpha} &= \Psi_\alpha, \\ 
c_\alpha +\Psi_\alpha(c_\alpha) &= \vs_\alpha, \\ 
\Psi_\alpha(\eps_\alpha) &= \lambda_\alpha/p^r. 
\endaligned 
\tag 10.1.2 
$$ 
\indent 
Let $\scr G^\star_0(\alpha)$ be the subset $\uL\|\scr C^\star(\frak a,\alpha)\|$ of $\scr W_F\backslash\wP F$. Every element of $\scr G^\star_0(\alpha)$ is a singleton orbit, so we may treat such orbits as individual representations of $\scr P_F$. Restriction to $\scr P_F$ gives a surjective map $\scr G^\star(\alpha) \to \scr G^\star_0(\alpha)$. Each fibre of this map is a principal homogeneous space over the group of tamely ramified characters of $\scr W_F$, as in \cite{12} 1.3 Proposition. 
\subhead 
10.2 
\endsubhead 
We give a relative characterization of the elements of $\scr G^\star(\alpha)$ in terms of the ultrametric pairing $\Delta$ on $\wW F$. 
\proclaim{Proposition} 
Let $\sigma\in \scr G^\star(\alpha)$ and $\tau\in \wwr F$. The following conditions are equivalent:   
\roster 
\item$\tau \in \scr G^\star(\alpha)$; 
\item $\dim\tau \le p^r$ and $\Delta(\sigma,\tau) \le \eps_\alpha$; 
\item $\dim\tau \le p^r$ and $\Hom{\scr R_F^+(\eps_\alpha)}\sigma\tau \neq 0$. 
\endroster 
\endproclaim 
\demo{Proof} 
We first work on the GL-side. 
\proclaim{Lemma} 
Let $\vT\in \|\scr C^\star(\frak a,\alpha)\|$ and $\vF \in \Scr E(F)$. The following are equivalent: 
\roster 
\item 
$\vF\in \|\scr C^\star(\frak a,\alpha)\|$; 
\item 
$\deg\vF \le p^r$ and $\Bbb A(\vF,\vT) \le \lambda_\alpha/p^r$. 
\endroster 
\endproclaim 
\demo{Proof} 
Let $\theta\in \scr C^\star(\frak a,\alpha)$ have endo-class $\vT$. If $\vF \in \|\scr C^\star(\frak a,\alpha)\|$, then $\deg\vF = p^r$ and $\vF$ is the endo-class of some $\phi\in \scr C^\star(\frak a,\alpha)$. By definition, $\phi$ agrees with $\theta$ on $H^{1+\lambda_\alpha}(\alpha,\frak a)$, whence $\Bbb A(\vF,\vT) \le \lambda_\alpha/p^r$. Thus (1) implies (2). 
\par 
Assume (2) holds. Since $\Bbb A(\vF,\vT) \le l_\alpha/2p^r < m/p^r$, we conclude that $\vs_\vF = m/p^r$: this follows from the definition of $\Bbb A$. As $\deg\vF \le  p^r$ and $p$ does not divide $m$, so $\deg\vF = p^r$ and $\vF$ has a realization $\phi\in \scr C(\frak a,\beta)$, for a simple stratum $[\frak a,m,0,\beta]$ in which $F[\beta]/F$ is totally ramified of degree $p^r$. The characters $\phi\Mid H^{1+\lambda_\alpha}(\beta,\frak a)$, $\theta\Mid H^{1+\lambda_\alpha}(\alpha,\frak a)$ intertwine in $\GL{p^r}F$ by hypothesis. Since $\lambda_\alpha < m/2$, (3.5.11) Theorem of \cite{15} allows us to replace $\phi$ by a conjugate to achieve $H^1(\beta,\frak a) = H^1(\alpha,\frak a)$ and $\phi\in \scr C(\frak a,\alpha)$. The characters $\phi\Mid H^{1+\lambda_\alpha}(\alpha,\frak a)$, $\theta\Mid H^{1+\lambda_\alpha}(\alpha,\frak a)$ intertwine and so are equal \cite{15} (3.3.2). That is, $\phi \in \scr C^\star(\frak a,\alpha)$, as required. \qed 
\enddemo 
In the proposition, the equivalence of (2) and (3) is the definition of $\Delta$. Write $[\sigma]_0^+ = \uL\vT$, $[\tau]_0^+ = \uL\vF$. In particular, $\vT \in \ewc F$ while $\vF \in \Scr E(F)$ is totally wild of degree at most $p^r$. We have $\Psi_\vT(\Delta(\sigma,\tau)) = \Bbb A(\vT,\vF)$. The definition (10.1.2) shows that $\Bbb A(\vT,\vF) \le \lambda_\alpha/p^r$ if and only if $\Delta(\sigma,\tau) \le \eps_\alpha$. The proposition thus follows from the lemma. \qed 
\enddemo 
\remark{Remark} 
In the lemma, the hypothesis $\deg\vF \le p^r$ is essential. For, the Density Lemma of \cite{13} 5.3 shows that the set of values $\Bbb A(\vT,\vF)$, $\vF \in \Scr E(F)$, is dense on the positive real axis. Indeed, the same proof shows that the set of $\Bbb A(\vT,\vF)$ is dense when $\vF$ is confined to the set of totally wild endo-classes. In the proposition, the hypothesis $\dim\tau \le p^r$ is likewise essential. Interpretation of the general case, with $\dim\tau$ unbounded, is the subject of \cite{13} 6.5 Corollary. 
\endremark 
\subhead 
10.3 
\endsubhead 
Let $j_\infty(\alpha) = j_\infty(F[\alpha]|F)$ be the greatest jump of the function $\psi_{F[\alpha]/F}$. 
\definition{Definition} 
Say that $[\frak a,m,0,\alpha]$ (or the element $\alpha$) is {\it $\star$-exceptional\/} if $j_\infty(\alpha) = c_\alpha$, $l_\alpha > 0$ and $l_\alpha \equiv 0 \pmod2$. Otherwise, say that $\alpha$ is {\it $\star$-ordinary.} 
\enddefinition 
Both exceptional and ordinary cases arise. If $\alpha$ is $\star$-exceptional, then $\Psi_\alpha$ has an odd number of jumps. Otherwise, both odd and even cases occur. We prove: 
\proclaim{Theorem}  
\roster 
\item 
There is a character $\xi$ of\/ $\scr R_F^+(c_\alpha)$ occurring in every representation $\sigma \in \scr G^\star(\alpha)$. This condition determines $\xi$ uniquely, up to $\scr W_F$-conjugation. In particular, each $\sigma\in \scr G^\star(\alpha)$ determines the same conjugacy class of field extensions $L_\xi/F$. 
\item
Suppose that $\alpha$ is $\star$-ordinary but that $\Psi_\alpha$ has an odd number of jumps. There is an irreducible representation $\rho_\xi$ of $\scr R_F(c_\alpha)$ that contains $\xi$ and occurs in every $\sigma \in \scr G^\star(\alpha)$. This condition determines $\rho_\xi$ uniquely, up to $\scr W_F$-conjugation. 
\item 
If $\alpha$ is $\star$-ordinary, then $\wt L_{\sigma,\xi} = \wt L_{\tau,\xi}$, for all $\sigma,\tau \in \scr G^\star(\alpha)$. 
\endroster 
\endproclaim 
\demo{Proof} 
We estimate the number $c_\alpha$ to get a more effective bound for the distance $\Delta(\sigma_1,\sigma_2)$, $\sigma_i \in \scr G^\star(\alpha)$. 
\proclaim{Lemma 1} 
Write $j_\infty(\alpha) = j_\infty(F[\alpha]|F)$. 
\roster 
\item 
If $j_\infty(\alpha) \le c_\alpha$, then $c_\alpha = (m{+}w_\alpha)/2p^r$ and $\Psi_\alpha(c_\alpha) = l_\alpha/2p^r$. 
\item 
If $j_\infty(\alpha) > c_\alpha$, then $c_\alpha < (m{+}w_\alpha)/2p^r$ and $\Psi_\alpha(c_\alpha) > l_\alpha/2p^r \ge \lambda_\alpha/p^r$. 
\endroster 
\endproclaim 
\demo{Proof} 
Suppose $j_\infty(\alpha) < c_\alpha$. The function $\Psi_\alpha$ then has an even number of jumps, its graph contains a non-empty open segment of the line $y=x{-}p^{-r}w_\alpha$, and $x = c_\alpha$ is the intersection of this line segment with $x{+}y=\vs_\alpha$ (4.2 Proposition). That is, $c_\alpha = (m{+}w_\alpha)/2p^r$ and so $\Psi_\alpha(c_\alpha) = c_\alpha{-}p^{-r}w_\alpha = l_\alpha/2p^r$. 
\par 
Suppose next that $j_\infty(\alpha) = c_\alpha$. Therefore $\Psi_\alpha(c_\alpha) = p^{-r}\psi_{F[\alpha]/F}(c_\alpha) =c_\alpha{-}p^{-r}w_\alpha$. Thus $2c_\alpha{-}p^{-r}w_\alpha = \vs_\alpha$ whence $c_\alpha = (m{+}w_\alpha)/2p^r$ and $\Psi_\alpha(c_\alpha) = l_\alpha/2p^r$as desired. 
\par 
In (2), the line $y=x{-}p^{-r}w_\alpha$ lies strictly below the graph $y = \Psi_\alpha(x)$, ({\it cf\.} 1.6 Proposition, 4.2 Proposition), giving 
$$ 
\vs_\alpha-c_\alpha = \Psi_\alpha(c_\alpha) > c_\alpha-p^{-r}w_\alpha, 
$$ 
and hence the first assertion. The three lines $y= l_\alpha/2p^r$, $y=x{-}p^{-r}w_\alpha$ and $x{+}y=\vs_\alpha$ all meet at $x = (m{+}w_\alpha)/2p^r$. As $(m{+}w_\alpha)/2p^r > c_\alpha$, so $\Psi_\alpha(c_\alpha) > l_\alpha/2p^r$. \qed 
\enddemo 
\proclaim{Lemma 2} 
\roster 
\item 
If $\sigma_1,\sigma_2 \in \scr G^\star(\alpha)$ then $\Delta(\sigma_1,\sigma_2) \le c_\alpha$. 
\item 
There exist $\sigma_1,\sigma_2 \in \scr G^\star(\alpha)$ such that $\Delta(\sigma_1,\sigma_2) = c_\alpha$ if and only if either 
\itemitem{\rm (a)} $j_\infty(\alpha) < c_\alpha$ and $l_\alpha$ is even, or 
\itemitem{\rm (b)} $\alpha$ is $\star$-exceptional. 
\endroster 
\endproclaim 
\demo{Proof} 
By 10.2 Proposition, 
$$ 
\roman{max}\,\{\Delta(\sigma_1,\sigma_2): \sigma_i \in \scr G^\star(\alpha\} = \eps_\alpha = \Psi_\alpha^{-1}(\lambda_\alpha/p^r). 
$$ 
By Lemma 1 above, $\Psi_\alpha(c_\alpha) \ge \lambda_\alpha/p^r$, whence $c_\alpha \ge \lambda_\alpha/p^r$. This proves (1). If $j_\infty(\alpha) > c_\alpha$, Lemma 1 gives $\eps_\alpha > c_\alpha$, so $\Delta(\sigma_1,\sigma_2) < c_\alpha$ in this case. If $j_\infty(\alpha)< c_\alpha$, Lemma 1 gives $c_\alpha = \Psi_\alpha^{-1}(l_\alpha/2p^r) \ge \Psi_\alpha^{-1}(\lambda_\alpha/p^r)$, with equality if and only if $l_\alpha$ is even. This accounts for option (a) in case (2). 
\par 
This leaves the case $j_\infty(\alpha) = c_\alpha$. If $l_\alpha\neq 0$, the same argument applies and gives option (b). It remains only to show that the conditions $j_\infty(\alpha) \le c_\alpha$ and $l_\alpha = 0$ are incompatible. 
\par 
Suppose these two conditions hold. We have $m\le w_\alpha$ while, by Lemma 1, $c_\alpha = (m{+}w_\alpha)/2p^r$. Now 1.6 Corollary implies 
$$  
c_\alpha = \frac{m{+}w_\alpha}{2p^r} \le \frac {w_\alpha}{p^r} \le \frac{p^r{-}1}{p^r}\,j_\infty(\alpha) < j_\infty(\alpha), 
$$ 
contrary to the hypothesis $j_\infty(\alpha) \le c_\alpha$. \qed 
\enddemo 
We prove the theorem. In part (1), choose $\sigma \in \scr G^\star(\alpha)$ and apply 9.5 Corollary. In the notation of that result, $\sigma\Mid \scr R_F^+(c_\alpha)$ is a direct sum of $\scr W_F$-conjugate characters $\xi$. If $\tau \in \scr G^\star(\alpha)$, Lemma 2 gives $\Delta(\sigma,\tau) \le c_\alpha$ whence any such $\xi$ occurs in $\tau$. The uniqueness property follows by symmetry. 
\par 
In part (2), take $\xi$ as in part (1) and set $L = L_\xi$. By definition, $\scr W_L$ is the $\scr W_F$-stabilizer of $\xi$ and we have $\scr R^+_F(c_\alpha) \subset \scr W_L$. Let $\sigma_\xi$ be the natural representation of $\scr W_L$ on the $\xi$-isotypic subspace of $\sigma$. The $\scr R_F(c_\alpha)$-normalizer of the character $\xi$ is $\scr W_L\cap \scr R_F(c_\alpha)$, by the definition of $L$. So, the representation $\rho_\xi$ of $\scr R_F(c_\alpha)$, induced by $\sigma_\xi\Mid \scr W_L \cap \scr R_F(c_\alpha)$, is irreducible. If $\tau \in \scr G^\star(\alpha)$, Lemma 2 asserts that $\Delta(\sigma,\tau) < c_\alpha$, so $\rho_\xi$ also occurs in $\tau$. The representation $\rho_\xi$ therefore has the required properties. 
\par 
Part (3) is trivial if $\Psi_\alpha$ has an even number of jumps, as then $\wt L_{\sigma,\xi} = L_\xi$. Assume otherwise. In the same notation as in the proof of part (2), $\sigma_\xi\Mid \scr W_L\cap \scr R_F(c_\alpha)$ is the natural representation on the $\xi$-isotypic subspace of $\rho_\xi$. Consequently, if $\tau\in \scr G^\star(\alpha)$, the representations $\sigma_\xi$, $\tau_\xi$ agree and are irreducible on $\scr W_L\cap \scr R_F(c_\alpha)$. Therefore $\tau_\xi \cong \chi\otimes\sigma_\xi$, for a character $\chi$ of $\scr W_L$ trivial on $\scr R_F^+(c_\alpha)$, and so $\sigma_\xi$, $\tau_\xi$ define the same projective representation of $\scr W_L$. Their centric fields are therefore the same. This proves (3) and completes the proof of the theorem. \qed 
\enddemo  
\subhead 
10.4 
\endsubhead 
We fix a character $\xi$ of $\scr R_F^+(c_\alpha)$, occurring in some, hence any, $\sigma\in \scr G^\star(\alpha)$. Let $\scr W_L$ be the $\scr W_F$-stabilizer of $\xi$. For $\sigma \in \scr G^\star(\alpha)$, let $\sigma_\xi$ denote the natural representation of $\scr W_L$ on the $\xi$-isotypic subspace of $\sigma$.  
\proclaim{Lemma 1} 
If $\Delta_L$ denotes the canonical ultrametric pairing on $\wW L$, then 
$$ 
\roman{max}\,\big\{\Delta_L(\sigma_\xi,\tau_\xi): \sigma,\tau\in \scr G^\star(\alpha)\big\} = \lambda_\alpha. 
\tag 10.4.1 
$$ 
\endproclaim 
\demo{Proof} 
By 10.2 Lemma, $\roman{max}\{\Bbb A(\vT,\vF):\vT,\vF \in \|\scr C^\star(\frak a,\alpha)\|\} = \lambda_\alpha/p^r$. So, 
$$
\roman{max}\,\big\{\Delta(\sigma,\tau):\sigma,\tau\in \scr G^\star(\alpha)\big\} = \Psi_\alpha^{-1}(\lambda_\alpha/p^r)  = \vf_{L/F}(\lambda_\alpha),  
$$ 
by (9.5.1). The relation (10.4.1) now follows from 1.4 Proposition. \qed 
\enddemo 
Let $k\ge0$ be an integer and $K/F$ a finite field extension. Let $\Gamma_k(K)$ be the group of characters of $K^\times/U^{1+k}_K$, sometimes viewed as characters of $\scr W_K$. Let $\Gamma^0_k(K)$ be the group of characters of $U^1_K/U^{1+k}_K$. 
\par
Let $\scr H(L,\xi)$ be the set of representations $\sigma_\xi \in \wwr L$, for $\sigma \in \scr G^\star(\alpha)$. The induction functor $\Ind_{L/F}$ then gives a bijection $\scr H(L,\xi) \to \scr G^\star(\alpha)$. 
\proclaim{Proposition} 
\roster 
\item 
If $\chi\in \Gamma_{\lambda_\alpha}(L)$ and $\kappa\in \scr H(L,\xi)$, then $\chi\otimes \kappa \in \scr H(L,\xi)$. 
\item 
If $\alpha$ is $\star$-ordinary, then the set $\scr H(L,\xi)$ is a principal homogeneous space over $\Gamma_{\lambda_\alpha}(L)$. 
\endroster 
\endproclaim 
\demo{Proof} 
By definition, the character $\chi$ is trivial on $\scr R^+_L(\lambda_\alpha)$, so $\Delta_L(\kappa,\chi\otimes\kappa) \le \lambda_\alpha$. The representation $\kappa^F = \Ind_{L/F}\,\kappa$ is irreducible, and lies in $\scr G^\star(\alpha)$. If $\rho$ is an irreducible component of $\Ind_{L/F}\,\chi\otimes \kappa$, it follows that $\Delta(\kappa^F,\rho) \le \phi_{L/F}(\lambda_\alpha) = \eps_\alpha$ and $\dim\rho \le p^r$. From 10.2 Proposition we deduce that $\rho\in \scr G^\star(\alpha)$ and $\dim\rho  =p^r$. That is, $\Ind_{L/F}\,\chi\otimes\kappa = \rho$ is irreducible and lies in $\scr G^\star(\alpha)$. Therefore $\chi\otimes\kappa \in \scr H(L,\xi)$. 
\par 
Let $\scr H_0(L,\xi)$ be the set of equivalence classes of  representations $\sigma_\xi\Mid \scr P_L$, $\sigma \in \scr G^\star(\alpha)$. Induction, from $\scr P_L$ to $\scr P_F$, gives a bijection $\scr H_0(L,\xi) \to \scr G_0^\star(\alpha)$. In the second part of the proposition, it is enough to show that $\scr H_0(L,\xi)$ is a principal homogeneous space over $\Gamma^0_{\lambda_\alpha}(L)$. The sets $\scr H_0(L,\xi)$, $\|\scr C^\star(\frak a,\alpha)\|$ and $\scr C^\star(\frak a,\alpha)$ are in canonical bijection, and $\scr C^\star(\frak a,\alpha)$ visibly has exactly $q^{\lambda_\alpha}$ elements, where $q$ is the cardinality of the residue field of $F$. This reduces us to showing that, if $\kappa \in \scr H_0(L,\xi)$ and $\chi\in \Gamma^0_{\lambda_\alpha}(L)$, $\chi\neq1$, then $\chi\otimes\kappa \not\cong \kappa$. 
\par 
If $j_\infty(\alpha) < c_\alpha$, the representation $\kappa$ is a character, and the result is obvious. To deal with the other cases, we need the following general fact. Recall that $j_\infty(L|F)$ denotes the largest jump of $\psi_{L/F}$. 
\proclaim{Lemma 2} 
If $\Psi_\alpha$ has an odd number of jumps, that is, if $j_\infty(\alpha) \ge c_\alpha$, then $j_\infty(L|F) < c_\alpha$. 
\endproclaim 
\demo{Proof} 
Let $d^-$ (resp\. $d^+$) be the left (resp\. right) derivative of $\Psi_\alpha$ at $c_\alpha$. Let $\dim\sigma_\xi = p^s = p^r/[L{:}F]$. The H-singular representation $\sigma_\xi$ of $\scr W_L$ is irreducible on $\scr R_F(c_\alpha)\cap \scr W_L$, but is a sum of $p^s$ copies of $\xi$ on $\scr R_F^+(c_\alpha)$ (which is contained in $\scr W_L$). Therefore $d^+/d^- = p^{2s}$. Symmetry (as in 3.1) implies that $d^+ = (d^-)^{-1}$, whence $d^- = p^{-s} = p^{-r}[L{:}F]$. So, if $\delta$ is small and positive, $\psi'_{L/F}(x) = [L{:}F]$ for $c_\alpha{-}\delta < x < c_\alpha$. It follows ({\it cf\.} 1.6 Proposition) that $j_\infty(L|F) \le c_\alpha{-}\delta < c_\alpha$, as required. \qed 
\enddemo 
Suppose that $j_\infty(\alpha) = c_\alpha$ and $l_\alpha$ is odd, or that $j_\infty(\alpha) > c_\alpha$. In either case, the function $\Psi_\alpha$ has an odd number of jumps. It follows from Lemma 2 that $j_\infty(L|F) < c_\alpha$, so $\scr R_F(c_\alpha) \subset \scr W_L$ by 1.9 Corollary 2. The restriction $\kappa\Mid \scr R_F(c_\alpha)$ is irreducible, since $\psi_{L/F}(c_\alpha)$ is the only jump of $\kappa$. If $\rho$ is a representation of $\scr P_L$ such that $\rho\Mid \scr R_F(c_\alpha) = \kappa\Mid \scr R_F(c_\alpha)$,  there is a {\it unique\/} character $\phi$ of $\scr P_L$, trivial on $\scr R_F(c_\alpha)$, such that $\rho = \phi\otimes\kappa$. By 10.3 Lemma 1, any $\chi\in \Gamma^0_{\lambda_\alpha}(L)$ is trivial on $\scr R_F(c_\alpha)$. The representations $\chi\otimes \kappa$ are therefore distinct, as $\chi$ ranges over $\big(U^1_L/U^{1+\lambda_\alpha}_L\big)\sphat$, and the proposition follows. \qed 
\enddemo 
\remark{Remarks} 
\roster 
\item 
Suppose that $\alpha$ is $\star$-ordinary. The set $\scr G^\star(\alpha)$ then inherits the structure of principal homogeneous space over $\Gamma_{\lambda_\alpha}(L)$, via the bijection $\Ind_{L/F}:\scr H(L,\xi) \to \scr G^\star(\alpha)$.
\item  
If $\alpha$ is $\star$-exceptional, there will, in many cases, exist non-trivial characters $\chi \in \Gamma_{\lambda_\alpha}(L)$ such that $\chi\otimes \sigma_\xi \cong \sigma_\xi$. This is incompatible with a principal homogeneous space structure. 
\endroster 
\endremark 
\subhead 
10.5 
\endsubhead 
We assume, in this sub-section, that $\alpha$ is {\it $\star$-exceptional.} We fix a character $\xi$ of $\scr R_F^+(c_\alpha)$ as in 10.3 Theorem and abbreviate $L = L_\xi$, $\wt L_\sigma = \wt L_{\sigma,\xi}$. Let $T_\sigma/L$ be the maximal tame sub-extension of $\widetilde L_\sigma/L$, and define the character group $D(\sigma_\xi)$ as in 8.2. 
\proclaim{Theorem} 
Suppose that $\alpha$ is $\star$-exceptional. 
\roster 
\item 
If $\sigma,\tau \in \scr G^\star(\alpha)$, then $T_\sigma = T_\tau$. 
\item 
The integer $d = |D(\sigma_\xi)|$ is independent of the choice of $\sigma\in \scr G^\star(\alpha)$. It satisfies $d^{\frac12}\le \dim\sigma_\xi = p^r/[L{:}F]$. 
\item 
There are, at most,  $d$ distinct Galois extensions of the form $\widetilde L_\sigma/L$, as $\sigma$ ranges over $\scr G^\star(\alpha)$. If $p$ does not divide $[T_\sigma{:}L]$, there are exactly $d$ such extensions. 
\endroster 
\endproclaim 
\demo{Proof} 
We gather some identities. First, $\Psi_\alpha(x) = p^{-r}\psi_{L/F}(x)$, $0\le x\le c_\alpha$, by 9.5 Corollary. Since  $j_\infty(\alpha) = c_\alpha$, 10.3 Lemma 1 gives $\Psi_\alpha(c_\alpha) = l_\alpha/2p^r$. Consequently 
$$ 
\psi_{L/F}(c_\alpha) = l_\alpha/2. 
\tag 10.5.1 
$$ 
In this situation, $j_\infty(\alpha) = c_\alpha > j_\infty(L|F)$ by 10.4 Lemma 2, so 
$$ 
\scr R_F(c_\alpha) = \scr R_L(\psi_{L/F}(c_\alpha)) = \scr R_L(l_\alpha/2) 
\tag 10.5.2 
$$ 
by 1.9 Corollary 2. Write $e_\sigma = e(T_\sigma|L)$, so that $\scr R_L(l_\alpha/2) = \scr R_{T_\sigma}(e_\sigma l_\alpha/2)$. 
The point $e_\sigma l_\alpha/2$ is the unique jump of $\widetilde L_\sigma/T_\sigma$, so 
$$ 
\scr R^+_F(c_\alpha) = \scr R^+_L(l_\alpha/2) = \scr R^+_{\widetilde L_\sigma}(e_\sigma l_\alpha/2), 
\tag 10.5.3 
$$ 
and  
$$ 
\scr W_{T_\sigma} = \scr W_{\widetilde L_\sigma}\scr R_L(l_\alpha/2). 
\tag 10.5.4 
$$ 
We prove part (1) of the theorem. 
\proclaim{Lemma 1} 
If $\sigma, \tau\in \scr G^\star(\alpha)$, then $T_\tau = T_\sigma$. 
\endproclaim 
\demo{Proof} 
By 8.2 Lemma and Proposition, the group $\scr W_{T_\sigma}$ is the common $\scr W_L$-stabilizer of the elements of the character group $D^+(\sigma_\xi)$. Dualizing (via 8.1 Proposition), the group $\scr W_{T_\sigma}$ is the $\scr W_L$-centralizer of $\sigma_\xi(\scr R_L(\psi_{L/F}(c_\alpha))$ modulo its centre ({\it cf\.} 8.4 Corollary). This centre, we assert, is independent of $\sigma$. The pairing $(x,y) \mapsto \xi([x,y])$ defines an alternating form on the $\Bbb F_p$-vector space $\scr R_L(l_\alpha/2)/\scr R_L^+(l_\alpha/2)$. Let $R$ be the inverse image, in $\scr R_L(l_\alpha/2)$, of the radical of this pairing. Since $\scr W_L$ fixes $\xi$, it normalizes $R$. The image $\sigma_\xi(R)$ is the centre of $\sigma_\xi(\scr R_L(l_\alpha/2))$. Thus $\scr W_{T_\sigma}$ is the $\scr W_L$-centralizer of the finite group $\scr R_L(l_\alpha/2)/R$ and so is  independent of $\sigma$. \qed 
\enddemo 
In part (2) of the theorem, the integer $\dim \sigma_\xi = p^r/[L{:}F]$ is certainly independent of $\sigma \in \scr G^\star(\alpha)$. By 8.2 Lemma (1), the order of the group $D(\sigma_\xi)$ is the number of fixed points for the natural action of $\scr W_L$ on $\scr R_L(l_\alpha/2)/R$, in the notation of the proof of Lemma 1. It is therefore independent of $\sigma$ and we have proved part (2) of the theorem. 
\par 
In light of part (1), we abbreviate $T = T_\sigma$. 
\proclaim{Lemma 2} 
Suppose that $T = L$. For $\tau\in \scr G^\star(\alpha)$, the following are equivalent: 
\roster 
\item 
$\widetilde L_\tau = \widetilde L_\sigma$; 
\item 
there is a character $\chi$ of $\scr W_L$, trivial on $\scr R_L^+(l_\alpha/2)$, such that $\tau_\xi \cong \chi \otimes \sigma_\xi$. 
\endroster 
\endproclaim 
\demo{Proof} 
Surely (2) implies (1), so suppose that (1) holds. The restrictions $\sigma_\xi' = \sigma_\xi\Mid \scr R_L(l_\alpha/2)$, $\tau_\xi' =\tau_\xi\Mid \scr R_L(l_\alpha/2)$ are irreducible, and each is a multiple of $\xi$ on $\scr R_L^+(l_\alpha/2)$. On the group $R$ (as in the proof of Lemma 1), each is a multiple of a character of $R$ extending $\xi$. Consequently, there is a character $\phi_R$ of $R$, trivial on $\scr R_L^+(l_\alpha/2)$, such that $\tau_\xi\Mid R = \phi_R \otimes \sigma_\xi\Mid R$. The character $\phi_R$ extends to a character $\phi$ of $\scr R_L(l_\alpha/2)$. For any such $\phi$, we have $\tau'_\xi = \phi\otimes \sigma_\xi'$. The projective representations $\bar\sigma_\xi$, $\bar\tau_\xi$ defined by $\sigma_\xi$, $\tau_\xi$ are therefore identical on $\scr R_L(l_\alpha/2)$. Each of these projective representations has $\scr W_{\widetilde L_\sigma} = \scr W_{\widetilde L_\tau}$ in its kernel, so $\bar\sigma_\xi$, $\bar\tau_\xi$ are the same on the group $\scr W_L = \scr W_T = \scr W_{\widetilde L_\sigma}\scr R_L(l_\alpha/2)$. That is, $\sigma_\xi$, $\tau_\xi$ are liftings to $\scr W_L$ of the same projective representation $\bar\sigma_\xi$. It follows that  $\tau_\xi\cong \chi\otimes \sigma_\xi$, for some character $\chi$ of $\scr W_L$ trivial on $\scr R_L^+(l_\alpha/2)$. \qed 
\enddemo 
In the case $T = L$, we have $D(\sigma_\xi) \subset \Gamma_{l_\alpha/2}(L)$, so Lemma 2 implies that the number of distinct fields $\widetilde L_\sigma/L$, $\sigma \in \scr G^\star(\alpha)$, is 
$$ 
\big|\Gamma_{l_\alpha/2}(L)\backslash \scr G^\star(\alpha)\big| = \big|\Gamma^0_{l_\alpha/2}(L)\backslash \scr G^\star_0(\alpha)\big|. 
$$ 
The set $\scr G^\star_0(\alpha)$ is in bijection with $\|\scr C^\star(\frak a,\alpha)\|$, and so has $q^{l_\alpha/2} = |\Gamma^0_{l_\alpha/2}(L)|$ elements, while each element of $\scr G_0^\star(\alpha)$ is fixed under twisting by exactly $d$ elements of $\Gamma^0_{l_\alpha/2}(L)$. Therefore $ \big|\Gamma^0_{l_\alpha/2}(L)\backslash \scr G^\star_0(\alpha)\big| =d$, as required for part (3) of the theorem in this case. 
\par 
Return to the general case and write $e = e(T|L)$. For $\sigma\in \scr G^\star(\alpha)$, write $\sigma_\xi^T = \sigma_\xi\Mid \scr W_T$. Thus $\sigma_\xi^T$ has centric field $\widetilde L_\sigma/T$. For $\sigma,\tau\in \scr G^\star(\alpha)$, Lemma 2 shows that  $\widetilde L_\sigma = \widetilde L_\tau$  if and only if there exists $\chi \in \Gamma_{el_\alpha/2}(T)$ such that $\tau^T_\xi = \chi\otimes \sigma^T_\xi$. So, if there exists  $\phi \in \Gamma_{l_\alpha/2}(L)$ such that $\tau_\xi = \phi\otimes \sigma_\xi$, then $\wt L_\tau = \wt L_\sigma$. Counting as before, there are at most $d =  |D(\sigma_\xi)|$ distinct fields $\wt L_\sigma$, as $\sigma$ ranges over $\scr G_0^\star(\alpha)$. We have proved the first assertion of part (3) of the theorem. 
\par 
In general, the relation $\tau^T_\xi = \chi\otimes \sigma^T_\xi$ implies $\chi/\chi^\gamma \in D(\sigma_\xi^T)$, for all $\gamma \in \Gal TL$. That is, $\chi$ defines a $\Gal TL$-fixed point in $\Gamma_{el_\alpha/2}(T)/D(\sigma_\xi^T)$. If $p$ does not divide $[T{:}L]$, this is equivalent to $\chi\in \Gamma_{l_\alpha/2}(L)/D(\sigma_\xi)$, since $D(\sigma_\xi)$ is the group of $\scr W_L$-fixed points in $D(\sigma_\xi^T)$ (8.2 Lemma). The final assertion follows. \qed 
\enddemo 
\remark{Remark} 
There are indeed cases of $p$ dividing $[T_\sigma{:}L]$ in the context of the theorem: we have already seen this in the example of 9.6. 
\endremark 
\subhead 
10.6 
\endsubhead 
{\it For this sub-section only, we assume that $p\neq 2$.} We outline a mild variant to our approach, following M\oe glin \cite{36}. It gives a simpler expression of the results, at the cost of a loss of generality. 
\par 
Otherwise, we use the notation from the beginning of the section. Suppose that $l_\alpha > 0$. Define $\scr C^\dagger(\frak a,\alpha)$ to be the set of $\theta\in \scr C^\star(\frak a,\alpha)$ satisfying 
$$ 
\theta(1{+}y) = \mu_M(\alpha(y-\tfrac12 y^2)), \quad y\in E,\ \ups_E(y) \ge [(1{+}l_\alpha)/2]. 
\tag 10.6.1 
$$ 
This expression does indeed define a character of $U^{[(1+l_\alpha)/2)]}_E$. Surely $\scr C^\dagger(\frak a,\alpha)$ is not empty. It is equal to $\scr C^\star(\frak a,\alpha)$ when $l_\alpha$ is odd. Let $\|\scr C^\dagger(\frak a,\alpha)\|$ be the set of endo-classes of characters $\theta\in \scr C^\dagger(\frak a,\alpha)$. In the case $l_\alpha = 0$, we may put $\scr C^\dagger(\frak a,\alpha) = \scr C^\star(\frak a,\alpha)$: remember that this set has only one element. 
\proclaim{Lemma 1}  
Let $\theta\in \scr C(\frak a,\alpha)$. There exists $\beta\in \roman P(\frak a,\alpha)$ such that $\theta\in \scr C^\dagger(\frak a,\beta)$. 
\endproclaim 
\demo{Proof} 
This follows readily from 7.1 Proposition. \qed 
\enddemo 
Let $\scr G^\dagger(\alpha)$ be the set of $\sigma \in \wwr F$ such that $[\sigma]_0^+ \in \uL\|\scr C^\dagger(\frak a,\alpha)\|$. The advantage of this approach is encapsulated in the following lemma. 
\proclaim{Lemma 2} 
For $\sigma,\tau \in \scr G^\dagger(\alpha)$, one has $\Delta(\sigma,\tau) < c_\alpha$. 
\endproclaim 
This follows from 10.3 Lemma 2. Imitating the discussion in 10.4 and 10.5, using the same notation, we find: 
\proclaim{Proposition} 
\roster 
\item 
If $\sigma, \tau \in \scr G^\dagger(\alpha)$, then $\wt L_\tau =  \wt L_\sigma$. 
\item 
Let $\lambda'_\alpha = \roman{max}\,\{[(1{+}l_\alpha)/2]{-}1,0\}$. The set $\scr G^\dagger(\alpha)$ is a principal homogen\-eous space over $\Gamma_{\lambda'_\alpha}(L)$. 
\endroster 
\endproclaim 
\subhead 
10.7 
\endsubhead 
Explicit results concerning the local Langlands correspondence fall into three areas. For essentially tame representations (which have trivial Herbrand functions), complete results are given in \cite{6}, \cite{7}, \cite{9}. A method for reducing to the totally wild case is worked out in \cite{12}. For totally wildly ramified representations, results are confined to a small number of old, but distinguished, papers. We briefly examine the relation between this paper and that historical context. 
\par 
Leaving aside the peripheral case of \cite{11}, the significant work concerns dimen\-sion $p$, in the context of proving the existence of the Langlands correspondence. The case $p = 2$ is in Kutzko \cite{32}, \cite{33} (as recounted in \cite{8}), $p=3$ is Henniart \cite{23} while $p\ge5$ is M\oe glin \cite{36}. 
\par 
The keystone of Kutkzo's work is the management of the case where, in the notation of the rest of the section, $\Psi_\alpha$ has a single jump. He proves that this is equivalent to $m\le 3w_\alpha$ (as we noted in 6.2 Example). He identifies the field we called $T_\sigma$ in 10.5: it is the splitting field of the polynomial $X^3-\roman{tr}(\alpha)X^2+\det(\alpha)$ (45.2 Theorem of \cite{8}). This approach is extended to odd $p$ in \cite{36} V.4 Proposition. A similar ``universal polynomial'' appears in \cite{11} 5.1 Theorem for epipelagic representations, i.e., those with Swan exponent $1$ in arbitrary dimension $p^r$. These results anticipate the more general 10.5 Theorem (1). 
\par 
Kutzko's construction of the Langlands correspondence has little to say about relating parameter fields $F[\alpha]$, $L_\xi$ on the  two sides. To define the correspondence, he relies on the Weil representation. That construction has remained resistant to further elucidation. 
\par 
M\oe glin's paper \cite{36}, for $p\ge5$, goes significantly further in that respect. It builds on Kutzko-Moy \cite{34} and Kutzko \cite{31} along with Carayol \cite{17}. It also relies on a number of working hypotheses that have since been verified, notably: 
\roster 
\item characterization of the Langlands correspondence via local constants of pairs (see \cite{24}); 
\item  compatibility of Kazhdan's lift \cite{30} and the Kutzko-Moy tame lift \cite{34} with Arthur-Clozel base change \cite{1} (see \cite{26}, \cite{4} respectively). 
\endroster 
All of those cited papers assume $F$ to be of characteristic zero. That restriction is removed in \cite{27} and \cite{28}. 
\par 
A feature of \cite{36} is the treatment of the relation between parameter fields. To re-arrange matters in accordance with the scheme here, we start with a simple stratum $[\frak a,m,0,\alpha]$ in $\M pF$ (as throughout) such that $m>w_\alpha$. Write $E = F[\alpha]$ and let $\theta\in \scr C^\dagger(\frak a,\alpha)$. Let $\chi_\theta$ be a character of $E^\times$ agreeing with $\theta$ on $U^1_E$. The representation $\sigma(\chi_\theta) = \Ind_{E/F}\,\chi_\theta$ is then irreducible, totally wild and of Carayol type. If $E/F$ is cyclic, then $\sigma(\chi_\theta)$ is absolutely wild. In this case, M\oe glin shows that the set of representations $\sigma(\chi_\theta)$, for $\theta\in \scr C^\dagger(\frak a,\alpha)$, is what we have called $\scr G^\dagger(\alpha)$. That is, the Langlands correspondence matches parameter fields. 
\par 
In general, the problem of describing parameter fields for H-singular representations seems to be of a different order. In the case of epipelagic representations \cite{11} (where $m=1$), the field $F[\alpha]$ is so ill-determined as to make the question meaningless without some qualification. 
 
\enddocument